%% file: AHirzebruchrevised.tex
\documentclass [10pt,reqno]{amsart}

\usepackage {hannah}
\usepackage {color,graphicx}
\usepackage {texdraw}
\usepackage [latin1]{inputenc}

\DeclareMathOperator {\val}{val}
\DeclareMathOperator {\trop}{trop}
\DeclareMathOperator {\irr}{irr}
\DeclareMathOperator {\ov}{ov}

\newcommand{\ZZ}{{\mathbb Z}}

\newcommand{\CC}{{\mathbb C}}
\newcommand{\RR}{{\mathbb R}}

\newcommand{\TT}{{\mathbb T}}

\renewcommand{\C}{{\mathcal C}}

\renewcommand{\S}{{\mathcal S}}
\renewcommand{\P}{{\mathcal P}}
\renewcommand{\L}{{\mathcal L}}
\renewcommand{\E}{{\mathcal E}}
\newcommand{\X}{{\mathcal X}}

\newcommand{\tg}{{\text{``}}}
\newcommand{\td}{{\text{''}}}
\newcommand{\Ve}{\text{Vert}}
\newcommand{\Ed}{\text{Edge}}

\newcommand{\ctor}{(\CC^*)^2}
\newcommand{\ctorn}{(\CC^*)^n} 
\newcommand{\Arg}{\operatorname{Arg}}

\newcommand{\sonen}{(S^1)^n}
\newcommand{\tj}{t_j}

\newcommand{\Log}{\operatorname{Log}}

\title [Deformation of tropical Hirzebruch surfaces]{Deformation of
  tropical Hirzebruch surfaces and enumerative geometry}

\author {Erwan Brugall\'{e} and Hannah Markwig}

\address{Erwan Brugall\'{e},  Universit\'e Pierre et Marie Curie,  Paris 6, 4 place Jussieu, 75 005 Paris, France}
\email{brugalle@math.jussieu.fr}
\address {Hannah Markwig, Universit\"at des Saarlandes\\ Fachrichtung Mathematik\\ Postfach 151150, 66041 Saarbr\"ucken\\ Germany}
\email {hannah@math.uni-sb.de}



\begin {document}

\begin {abstract}
We illustrate the use of tropical methods by
 generalizing a formula due to Abramovich and Bertram, extended later
 by Vakil. Namely,
we exhibit   relations between enumerative invariants  of the
Hirzebruch surfaces $\Sigma_n$ and 
$\Sigma_{n+2}$, obtained by
 deforming the first surface to the latter.

Our strategy involves
a tropical counterpart of 
deformations of  Hirzebruch
surfaces, and  tropical enumerative geometry on a tropical
surface in three-space.
\end {abstract}

\maketitle
\tableofcontents

\section{Introduction}

\subsection{Results}
In \cite{AB01}, Abramovich and Bertram related genus 0 enumerative invariants
of  $\Sigma_0=\CC P^1\times \CC P^1$ and the second
Hirzebruch surface $\Sigma_2$ (i.e. the quadratic cone $x^2+y^2+z^2=0$
in $\CC P^3$ blown up at the node). The strategy of their proof
is to understand how algebraic curves on $\Sigma_0$ behave when this
latter surface 
deforms
to  $\Sigma_2$.
 Later on this method was extended by Vakil in \cite{Vak00b} 
to non-rational
enumerative invariants of $\Sigma_0$ and $\Sigma_2$, and more generally
to relate enumerative invariants of an almost Fano complex surface and
any of its deformations. 

The goal of this paper is to 
illustrate the use of tropical  deformation techniques by
 generalizing  
this formula
to the case of $\Sigma_n$ and $\Sigma_{n+2}$.
Here, $n\ge 0$ and 
$\Sigma_n= \mathbb P(\mathcal O_{\CC P^1}(n)\oplus \CC )$
is the so called \emph{$n$-th Hirzebruch surface}.
Hirzebruch surfaces are toric.

We consider enumerative invariants $N_{\chi}(\delta)$ counting curves
in $\Sigma_n$ with prescribed Euler characteristic $\chi$, data of
intersections with the toric boundary, and enough point conditions (see
Definition \ref{def-toricrelinv}). We encode the intersections with
the toric boundary in terms of a \emph{Newton fan} $\delta$ (see
Definition \ref{def-Newtonfan}). A Newton fan is a multiset of vectors
$\delta=\{v_1^{m_1},\ldots, v_k^{m_k}\}$ (where the notation used here
indicates that the vector $v_i$ appears $m_i$ times) that 
refines the notion of the fan dual 
to a Newton polygon.
In particular, to any Newton fan $\delta\subset \mathbb{Z}^2$ we can naturally associate a
dual polygon $\Pi_\delta$.
The toric
surface $\Sigma_n$ in which we consider the curves to be counted is
determined once we fix a Newton fan
(see Notation \ref{deltazero}). 

In our formula, we also consider enumerative invariants that we call
the (1-1)-relative invariants of $\Sigma_n$ (see Definition
\ref{def:rel inv}). They count curves passing through prescribed
points and having prescribed tangent germs at 
a fixed point.
We refer
to these enumerative invariants as $\mathcal N(u,n,d,\alpha)$, they
depend on the integer $u$ prescribing the number of 
irreducible
components, the integer $n$ defining the Hirzebruch surface $\Sigma_n$
in question, and two vectors $d$ and $\alpha$ encoding the
intersection multiplicities with boundary divisors of the tangent
germs. 
For fixed $n$ and $u$, there exist only finitely many possibilities
for $d$ and $\alpha$, as explained (along with further details on the
computation of these invariants) in Remark \ref{rem:finitelt many
  N(u)}. 

Before stating our formula relating these enumerative invariants, we need to introduce some notation.

\begin{notation}\label{deltazero}
 Let $\delta_0$ be the Newton fan 
 $$\delta_0=\{(1,n)^a, (0,-1)^{an+b}, (-1,0)^a, (0,1)^b \}.$$ 
We write
 $\delta \vdash \delta_0$ if $\delta$ is a Newton fan satisfying
$$\begin{array}{ll}\delta= & 
\{(1,n+2)^m, (0,-1)^{a(n+1)+b}, (-1,0)^A, 
\\ &\quad \quad (-\alpha_1,\beta_1),\ldots,
(-\alpha_r,\beta_r), (0,\beta_{r+1}),\ldots,
(0,\beta_{r+s}), (0,1)^U \}
\end{array}$$
with
$$0< m\le a,\quad 0\leq A\leq \min\{m, b\},$$
$$\alpha_i,\beta_i>0 \quad \text{for}\quad i=1\ldots,r, \quad \text{and} 
\quad \beta_{r+1},\ldots, \beta_{r+s}>1.$$

Note that for every $\delta \vdash \delta_0$,
$$\Pi_\delta\subset Conv\{(0,0), (0,a), (b-a,a), \left(a(n+1)+b,0 \right)\},$$
where the latter is the polygon dual to $$\{(1,n+2)^a, (0,-1)^{a(n+1)+b}, (-1,0)^a,
(0,1)^{b-a} \}.$$ 
In particular
once $\delta_0$ is fixed, 
the choices of such $\delta$ are limited (see Figure \ref{fig delta}).

\begin{figure}[ht]
\centering
 \input{Deltadelta.pstex_t}
\caption{Finitely many fans $\delta\vdash \delta_0$.}
\label{fig delta}
\end{figure}
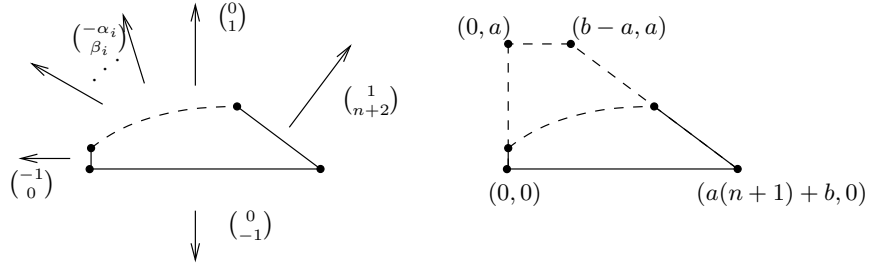

For a Newton fan $\delta\vdash\delta_0$, 
we define the following
quantities:
\begin{itemize}

\item $\chi'=\chi- 2(a+b-m-r-s-A-U )$;
\item $d= (\beta_1+\alpha_1,\ldots,\beta_r+\alpha_r, 
\beta_{r+1},\ldots, \beta_{r+s} , 1^{U+A-b})$;
\item $ \alpha= (\alpha_1, \ldots,\alpha_r )$.

\end{itemize}

\end{notation}

Our main result is a formula relating enumerative
invariants
of  $\Sigma_n$ and $\Sigma_{n+2}$ which we prove at the end of
Section \ref{sec:proof main}.
\begin{theorem}\label{main formula}
Let $n\ge 0$ and $\chi\in\ZZ$ be two integers, and let $\delta_0$ as in Notation \ref{deltazero} be the Newton fan 
 $$\delta_0=\{(1,n)^a, (0,-1)^{an+b}, (-1,0)^a, (0,1)^b \}.$$ Then we
have the following equation: 
\begin{align*} 
N_\chi(\delta_0)=\sum_{\delta\vdash \delta_0} \left( \binom{U}{b-A}
 \cdot\prod_{i=1}^r\gcd(\alpha_i,\beta_i) \cdot \prod_{i=1}^s\beta_{r+i}\cdot
   \mathcal N(a-m,n,d,\alpha) \ N_{\chi'}(\delta) \right).
\end{align*}

\end{theorem}

Note that the factor $N_{\chi'}(\delta)$ indeed counts curves on $\Sigma_{n+2}$
by Notation \ref{deltazero} and Example \ref{ex faninsigman}.
Note also that as pointed out in Remark \ref{rem:finitelt many N(u)}, 
 for a fixed $n$ the computations of all numbers 
$ \mathcal N(a-m,n,d,\alpha)$ reduces to the computation of
finitely many cases.

We use methods from tropical geometry to prove Theorem \ref{main formula}.
In the tropical world, algebraic curves are replaced by certain
balanced piece-wise linear graphs called tropical curves. Tropical
geometry has gained lots of attention recently. One of the interesting
results is that we can determine enumerative invariants of toric
surfaces by counting the corresponding tropical curves instead (Mikhalkin's
Correspondence Theorem, see \cite{Mi03}). This fact is at the base of
our tropical approach to enumerative invariants of Hirzebruch
surfaces. 

The proof of  Theorem \ref{main formula}
basically consists of the suitable Correspondence Theorems
(see Theorem \ref{thm-corres} and Theorem \ref{thm:corres 2}) relating
the above enumerative invariants to their tropical counterparts, and a
proof of the tropical version of the above formula (see Section
\ref{sec:proof main}).
The main idea underlying our tropical strategy is a deformation of
tropical Hirzebruch surfaces, which is a tropical analogue of
Kodaira's deformation of complex Hirzebruch surfaces.
We provide a more detailed description of our method  in
Section \ref{sec:method}.

The two Correspondence Theorems \ref{thm-corres} and \ref{thm:corres 2} and the formula of Theorem \ref{main formula} together with its tropical method of proof should be viewed as the main contributions of our paper.



As examples of applications, 
we specialize Theorem \ref{main formula} to the cases $n=0,1$ and $2$.

\begin{example}[The case $n=0$]\label{subsec-n0}
If $n=0$ our formula reduces to the Abramovich-Bertram Formula \cite{AB01}
relating enumerative invariants of $\Sigma_0$ and $\Sigma_2$. 
According to Example 
\ref{exa simple trop enum}, $\mathcal N(u,0,d,\alpha)\ne 0$ if and only if 
$d=(1^u)$ and $\alpha=0$, so the only Newton fans $\delta$ which contribute to 
$N_\chi(\delta_0)$ are of the form
$$\delta=\delta_m=\{(1,2)^m, (0,-1)^{a+b}, (-1,0)^m, (0,1)^{a+b-2m} \},$$
for which we have $r=s=0$, $A=m$, $U=a+b-2m$, and $\chi'=\chi$.
Eventually, Theorem \ref{main formula} reduces to  Abramovich-Bertram Formula 
$$N_\chi(\delta_0)=\sum_{m=0}^a \left(\begin{array}{c} a+b-2m  \\ a-m  \end{array}\right)
\cdot N_\chi\big(\delta_m).  $$

We point out that the above sum starts at $m=0$ although the original
Abramovich-Bertram formula starts at $m=1$. This difference comes
from the fact our formula involves \emph{reducible} enumerative
invariants,
 whereas the original involves \emph{irreducible} invariants. In this
latter situation one has  $ N^{irr}_\chi\big(\delta_0)=0$ which explains the
difference in the initial value in both summations.

In the special case $n=0$, it is immediate to adapt the proof of
Theorem \ref{main formula} to get a formula only involving 
 irreducible curves, and
of course one recovers the original  Abramovich-Bertram formula
$$N_\chi^{irr}(\delta_0)=\sum_{m=1}^a \left(\begin{array}{c} a+b-2m  \\ a-m  \end{array}\right)
\cdot N_\chi^{irr}\big(\delta_m).  $$

However in the general case, the most natural way of phrasing Theorem
\ref{main formula} is via reducible enumerative invariants.
See also Proposition \ref{prop:application}
for another application of our methods to irreducible enumerative
invariants.
\end{example}

\begin{example}[The case $n=1$]\label{subsec-n1}
We now set up a formula relating enumerative invariants of $\Sigma_1$ and
 $\Sigma_3$. 
This formula can also be 
deduced from
Ionel's and Parker's symplectic sum formula \cite{IP00} 
(note however
that this latter provides another, though equivalent,
formula than ours).
According to Example 
\ref{exa simple trop enum}  we have a nonzero contribution for a fan
$\delta$ only if 
$d=(2^{u-t_a},1^{2t_a})$   and
$\alpha=(1^{u-t_a-t_b})$.

By notation \ref{deltazero}, 
we thus have directions $(-\alpha_i,\beta_i)=(-1,1)$ and
$(0,\beta_{r+i})=(0,2)$
 only. Hence the
only Newton fans $\delta$ which contribute to 
$N_\chi(\delta_0)$ are of the form
$$\delta=\delta_{m,t_a,t_b}=\{(1,3)^m, (0,-1)^{2a+b},
(-1,0)^{A},(0,2)^{t_b},(-1,1)^R, (0,1)^{b-A+2t_a} \} $$ 
with $R=a-t_a-t_b-m$ and $A=2m+t_a+t_b-a$.
Theorem \ref{main formula} reduces to
$$N_\chi(\delta_0)=\sum_{m=0}^a\ \sum_{t_a+t_b=0}^{a-m} 2^{t_b}\cdot(2t_a-1)!!
 \cdot\left(\begin{array}{c} a+b+t_a-t_b-2m \\ 2t_a  \end{array}\right)
\cdot N_{\chi + 2t_a}\big(\delta_{m,t_a,t_b}).$$

\end{example}

\begin{example}[The case $n=2$]\label{subsec-n2}

According to Example 
\ref{exa simple trop enum}  we have a nonzero contribution only for  fans
 of the form
$$\begin{array}{cl}
\delta=\delta_{m,\underline{t}}= &\{(1,4)^m, (0,-1)^{3a+b},
(-1,0)^{A},(-2,1)^ {t_a}, (-1,2)^{t_b}, 

\\ &\quad  (-1,1)^{t_c},
(0,3)^{t_f},(0,2)^{t_e}, (0,1)^{b-A+3t_d+t_c+t_e} \} 
\end{array}$$ 
with 
$$\underline{t}=(t_a,t_b,t_c,t_d,t_e,t_f), $$  
$$t_a=a-t_c-t_b-t_f-t_e-t_d-m, \quad \mbox{and}\quad
 A=3m+2(t_f+t_e+t_d-a)+t_b+t_c.$$
Theorem \ref{main formula} reduces to
$$
N_\chi(\delta_0)=\sum_{m=0}^a\ \sum_{t_b+t_c+t_d+t_e+t_f=0}^{a-m} 
\mathcal D_{m,\underline{t}} 
 \cdot
 N_{\chi + 2t_c+2t_e+4t_d}\big(\delta_{m,r,s})$$
with
$$\mathcal D_{m,\underline{t}} = 3^{t_b+t_f} \cdot2^{t_e}\cdot 
  \binom{3t_d+t_c+t_e}{t_c, t_e,3,\ldots,3} \cdot 
\frac{t_c!t_e!}{t_d!} \cdot \left(\begin{array}{c} 2a+b+t_d-t_e-2t_f-3m -t_b \\ 3t_d+t_c+t_e  \end{array}\right) .$$
There are $t_d$ copies of $3$ in the above multinomial coefficient.
\end{example}



\subsection{Our method}\label{sec:method}

Inspired by Kodaira's deformation of Hirzebruch surfaces, we consider an analogous deformation in the tropical world (see Appendix \ref{hirz}).
Since the tropical world has a \emph{discrete} nature, quite often
one does not need to go to the limit to actually \emph{see} the
limit. This is also the case here: for our proof, it is not really necessary to fully formalize the point of view of a tropical Kodaira deformation, and for that reason we only present it in the Appendix \ref{hirz}, as motivation and background for our tropical deformation technique (see, in particular, Remark \ref{rem:justification}).
For our technique, it is enough to work on a tropical model of $\TT\Sigma_n$ close enough to the limit of the deformation.

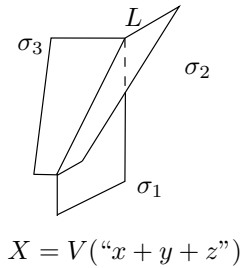
\begin{figure}[ht]
\centering
 \input{surfaceX.pstex_t}
\caption{The tropical surface $X$.}
\label{fig X2}
\end{figure}

We model this closeness to the limit as follows: the tropical model
of  (an open part of) $\TT\Sigma_n$ 
we use is the tropical surface $X$ in $\RR^3$ defined by the tropical
polynomial $\tg x+y+z \td$. It consists of three
half-planes $\sigma_1=\{x=y\geq z\}$, $\sigma_2=\{x=z\geq
y\}$ and $\sigma_3=\{y=z\geq x\}$ meeting along the line
$L=\{x=y=z\}=\RR(1,1,1)$, see Figure \ref{fig X2}.
 The configuration of points through which we require curves to pass
 through are chosen on the face $\sigma_1$ very far down from the line
 $L$.
In this setting, the tropical curves we are enumerating  naturally
contain two parts: one corresponding to the face $\sigma_1$, and one
corresponding to the two upper faces $\sigma_2$ and $\sigma_3$. The
former correspond to curves in
$\TT\Sigma_{n+2}$, and the latter to curves in $\TT\Sigma_{n}$. This behaviour is parallel to Kodaira's deformation in the complex world as we discuss in more detail in Section \ref{subsec-background}.

Note that this ambient tropical surface $X$, different from the
 ``usual'' $\RR^2$, is
 imposed
 by our strategy based on deformation of tropical Hirzebruch
 surfaces. Indeed the presence of a unique exceptional curve on
 $\Sigma_n$, with different self-intersection for different values
 of $n$,  is an obstruction for any deformation of $\Sigma_n$ to
 $\Sigma_{n+2}$ to be toric. As a consequence, to model this
 deformation tropically one needs to make use of  \emph{tropical
   modifications} (see \cite{Mi06}) of tropical Hirzebruch
 surfaces. In particular the tropical model of $\TT \Sigma_n$ involved
 in this deformation is no longer a tropical toric compactification of $\RR^2$.

Our two main ingredients for the proof of Theorem \ref{main formula} are the enumeration of tropical curves in the tropical surface
$X$ and its relation, via Correspondence Theorems \ref{thm-corres} and \ref{thm:corres 2}, to the enumeration of
algebraic curves in Hirzebruch surfaces. Both aspects have not been much explored in the literature yet. Most papers about tropical enumerative geometry deal with the case of tropical
curves in $\RR^n$. 
In the case of curves in $X$, the tropical inclusion
is more subtle  than
the set theoretic one, as it has
 already been noticed by several people (see
for example \cite{Vig1,Br12,Br17}). In other words a tropical curve 
$C$ might be set-theoretically contained in $X$ without being
tropically contained in $X$. One has to require extra conditions for
this latter inclusion to hold.
This phenomenon quite complicates the   enumerative geometry of
general tropical varieties. However, the surface $X$ we are interested
in here is simple enough so that only one extra condition,
the so-called
\emph{Riemann-Hurwitz condition},  suffices to rule
out parasitic tropical curves. The necessity of this condition
 has
been observed earlier by the first author and Mikhalkin (\cite{Br12}, see
also \cite{Br13}, \cite{ABBR13a},
 and \cite{ABBR13b}). 
In the other direction, we prove Correspondence Theorems 
for tropical curves in $X$ by
reducing to Correspondence Theorems for tropical curves in
$\RR^2$ with a multiple point at some fixed point on the toric
boundary. Thanks to this reduction, our proofs
goes by a mild 
adaptation of the proofs of Correspondence Theorems from \cite{Mi03}
and \cite{Shu12}, in the framework of phase-tropical morphisms given in
\cite{Br9}.

\subsection{Background, context, and difficulties}\label{subsec-background}


Kodaira proved that the surface $\Sigma_n$ can be 
deformed
to 
 $\Sigma_{n+2k}$ with $k\ge 0$.
As in \cite{AB01}, one may try to relate enumerative invariants of
$\Sigma_n$ and $\Sigma_{n+2k}$ by studying 
the limit of curves  when $\Sigma_n$ 
deforms
to 
 $\Sigma_{n+2k}$. However this
analysis  gets much more complicated when
$n>0$. The main reason for that is that  as soon as $n>0$  both surfaces 
 $\Sigma_n$ and $\Sigma_{n+2k}$ contain an exceptional curve, say $E_n$
and $E_{n+2k}$, and that curves in $\Sigma_n$ degenerate to curves
in $\Sigma_{n+2k}$ with 
singularities at $k$ fixed points on
 $E_{n+2k}$. 
Those latter points
may be
thought as the ``virtual intersection 
points
 of $E_n$ and $E_{n+2k}$''
defined by the chosen 
deformation
from $\Sigma_n$ to $\Sigma_{n+2k}$.

Let us explain in details the origin of the complications arising when $n>0$.

First let us
decompose
  the Kodaira 
deformation
of $\Sigma_n$ to $\Sigma_{n+2k}$ into two steps:
 a deformation of $\Sigma_n$ to the normal cone of a curve of
 bidegree $(1,k)$,
followed by the blow-down of the $\Sigma_n$ copy in the
special fiber.
More precisely, let $V$ be a non-singular curve of bidegree
$(1,k)$
  in
$\Sigma_n$\footnote{The presentation 
$\Sigma_n= \mathbb P(\mathcal O_{\CC P^1}(n)\oplus\CC)$ 
induces a projection 
$\pi:\Sigma_n\to\CC P^1$, and
 a curve of bidegree $(1,k)$ 
is characterized by the fact that it intersects in exactly $k$ points
the exceptional section
 and in one point any fiber of the map $\pi$; see Section \ref{def hirz} for details.}, and let
$\Sigma$ be the trivial family $\Sigma_n\times \CC $ blown-up along
the curve $V\times\{0\}$. The natural projection
$\pr:\Sigma \to\CC$ defines a 
flat degeneration 
of $\Sigma_n$ into the
reducible surface $\pr^{-1}(0)=\Sigma_n\cup \mathbb
P(\mathcal{N}_{V/\Sigma_n}\oplus \CC )$ intersecting transversally along
$V\subset\Sigma_n$ and $V_\infty\subset \mathbb
P(\mathcal{N}_{V/\Sigma_n}\oplus \CC )$. 
Since $V$ has self-intersection number 
$n+2k$ in $\Sigma_n$, we get that $\mathbb
P(\mathcal{N}_{V/\Sigma_n}\oplus \CC )=\Sigma_{n+2k}$ and
$V_\infty=E_{n+2k}$ is the exceptional curve in $\Sigma_{n+2k}$.

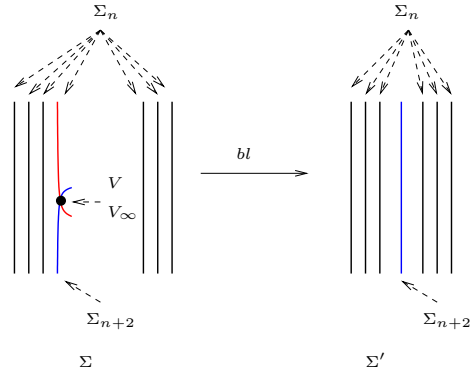
\begin{figure}[ht]
\centering
 \input{Def.pstex_t}
\caption{The $3$-folds $\Sigma$ and $\Sigma'$.}
\label{fig:deg sigma_n}
\end{figure}

It turns out that the $\Sigma_n$
copy in  $\pr^{-1}(0)$
can be contracted to $V$ by a 
blow down $\mbox{bl}:\Sigma \to \Sigma'$, and that the
induced projection $\pr':\Sigma'\to\CC$ is precisely the Kodaira deformation of
$\Sigma_n$ to 
$\Sigma_{n+2k}$ (see Figure \ref{fig:deg sigma_n}). 
With this picture in mind,
the  ``virtual intersection points of $E_n$ and $E_{n+2k}$'' 
we
mentioned above are now simply the intersection points of $V$ and
$E_{n}$ in the  $\Sigma_n$ copy  of the central fiber
$\pr^{-1}(0)\subset\Sigma$. 
Relations between enumerative invariants of $\Sigma_n$ and
$\Sigma_{n+2k}$ should now be derived
 from a careful analysis of how curves in
$\Sigma_n$ degenerate when this latter surface degenerates to 
$\Sigma_n\cup \Sigma_{n+2k}$.

At this point, 
the origin of the complications might appear more clearly with 
a symplectic point of view on the problem
and the methods.
On the level of the underlying symplectic manifolds (recall that
$\Sigma_n$ is a K\"ahler 
manifold), the above deformation $\Sigma$ of the reducible surface  
$\Sigma_n\cup \Sigma_{n+2k}$ to $\Sigma_n$ can be seen as a
symplectic sum of $\Sigma_n$ and $\Sigma_{n+2k}$ glued respectively
along $V$ and $E_{n+2k}$. There exist powerful symplectic sum formulas 
(see
\cite{IP00,LiRu01,EGH}, or \cite{l:adffgwi} for a degeneration formula in
the algebraic setting) relating Gromov-Witten invariants of a symplectic
sum with those of the symplectic summands, but these Gromov-Witten
invariants are \emph{relative} to some \emph{smooth} symplectic
hypersurface. Symplectic sum formulas provide an alternative proof of the
Abramovich-Bertram-Vakil results stated above, and more generally
express Gromov-Witten invariants of a symplectic $4$-manifold in term of its
Gromov-Witten invariants relative to an 
embedded symplectic
sphere $E$
with self-intersection $-l$.
Note however that as soon as $l\ge 2$, these formulas involve
enumeration of ramified coverings, which singularly complicates 
 actual computations when $l\ge 3$.
 Still, one can  deduce
in this way
a relation 
among
enumerative invariants of
$\Sigma_1$ and $\Sigma_3$
 equivalent to our Theorem \ref{main formula} specialized to
the case $n=1$ (see Example \ref{subsec-n1}).
However when $n>1$,  relating enumerative invariants of $\Sigma_n$ and
$\Sigma_{n+2k}$ requires to consider Gromov-Witten
invariants relative to a \emph{singular} symplectic curve.
Indeed, the complex structure on the algebraic surface $\Sigma_n$ is
not generic as soon as $n>1$ and enumerating algebraic curves on
$\Sigma_n$ is the same than computing Gromov-Witten invariants of the
underlying symplectic manifold  relative to the symplectic divisor
$E_n$. In particular,  
relating enumerative invariants of $\Sigma_n$ and $\Sigma_{n+2k}$ using
the deformation $\Sigma$ can be seen as expressing Gromov-Witten invariants of
$\Sigma_n$ relative to $E_n$ in terms of Gromov-Witten invariants of
$\Sigma_n$ relative to $E_n\cup V$ and of $\Sigma_{n+2k}$ relative to $E_{n+2k}$.
Since the curves $E_n$ and $V$ intersect in $k$ points, Gromov-Witten
invariants relative to a singular divisor show up naturally with the
method we intend to apply.
Up to our knowledge, Gromov-Witten invariants relative to a singular
divisor have been defined only recently (see \cite{Io12,Par11}, or \cite{GrSi13,AbCh11} in the  algebraic setting)  and
there does not exist a general symplectic sum/degeneration formula
for  those invariants yet.

We view our tropical approach as a powerful tool to overcome these problems.

If one 
translates the above strategy in the tropical setting, then
the family $\Sigma$ is replaced  by a single tropical surface $X$, and the
study of degenerations of holomorphic curves is replaced by the
enumeration of tropical curves in $X$. 
We then perform this enumeration in the special case of the degeneration of $\Sigma_n$  to 
$\Sigma_n\cup \Sigma_{n+2}$. We obtain in this way 
Theorem \ref{main formula}, which 
may be seen as such a symplectic sum/degeneration formula in some particular
instance of normal crossing divisor.
The case of the degeneration of $\Sigma_n$  to 
$\Sigma_n\cup \Sigma_{n+2k}$ 
should also be doable
tropically, but 
requires some additional efforts (see Section \ref{sec:rem}).

\medskip
The aim of
the above discussion is  to
replace
our work in the context of current mathematical developments
and to
explain where the difficulties we have to deal with come from.
Having said that, we  formulate our main result, Theorem 
\ref{main formula}, in the algebraic language without referring 
explicitly to
relative Gromov-Witten invariants.
Hence  our formula expresses (some) enumerative invariants 
of $\Sigma_n$ in
terms of (some) enumerative invariants of $\Sigma_{n+2}$  and 
some finitely many 
simple relative enumerative invariants of $\Sigma_n$,  both
of these latter invariants 
involving curves
with a prescribed very singular point.

\subsubsection*{Organization of the paper}
Our paper is organized as follows. In Section \ref{sec-formula} we
define Hirzebruch surfaces, 
and
the enumerative invariants we are
interested in.
We recall basic tropical definitions in Section \ref{sec:basic geotrop}, and state a
 Correspondence Theorem (Theorem \ref{thm:corres 2})
to compute $(1,1)$-relative invariants of
$\Sigma_n$ in Section \ref{sec:inter trop numbers}. 
Section \ref{sec-tropenum} is devoted to the study
 of 
basic
tropical enumerative geometry in $X$. In particular we state a
 second Correspondence Theorem (Theorem \ref{thm-corres})
 in this setting, and prove Theorem
 \ref{main formula}. Both Correspondence Theorems are proved in
 Section \ref{sec-corres}. Section \ref{sec:rem} contains some
  remarks about possible
 extensions of our present work. 
We end this paper by giving
 in Appendix \ref{hirz} the tropical counterpart of Kodaira deformation of
 Hirzebruch surfaces.

\subsubsection*{Acknowledgements}
Part of this work was accomplished at the Mathematical Sciences Research
Institute (MSRI) in Berkeley, CA, USA, during the semester program in
Fall 2009 on 
Tropical Geometry. The authors would like to thank the MSRI for
hospitality. E.B. is
partially supported by the ANR-09-BLAN-0039-01 and
ANR-09-JCJC-0097-01.
H.M. is partially supported by DFG-grant 4797/5-1 and by GIF-grant 1174-197.6/2011.
We are thankful to an anonymous referee for helpful remarks on an earlier version.

\section{Complex Hirzebruch surfaces}\label{sec-formula}

\subsection{Definition}\label{def hirz}
A \textit{Hirzebruch surface}, also called
 a rational geometrically ruled surface,
is a compact complex surface which admits a holomorphic fibration
to $\CC P^1$ with fiber
 $\CC P^1$.
The classification of Hirzebruch surfaces 
is well known (see for example \cite{Beau}): they are all isomorphic
to exactly one of the surfaces 
$\Sigma_n=\mathbb P(\mathcal O_{\CC P^1}(n)\oplus \CC)$
 with $n\ge 0$.
The surface $\Sigma_n$ is called the \emph{$n$th Hirzebruch surface}. For example $\Sigma_0=\CC P^1\times \CC P^1$.
Let us denote by $B_n$ (resp. $E_n$, and $F_n$) the section 
$\mathbb P(\mathcal O_{\CC P^1}(n)\oplus \{0\})$ (resp. the section
$\mathbb P(\{0\}\oplus  \CC)$, and a fiber). 
The curves $B_n$, $E_n$, and $F_n$ have self-intersections
$B_n^2=n$, $E_n^2=-n$, and  $F_n^2=0$. When $n\ge 1$, the curve $E_n$ itself
determines uniquely the 
Hirzebruch surface since it is the only reduced and irreducible
algebraic curve in $\Sigma_n$ with negative self-intersection.
For example $\Sigma_1$ is the projective plane blown up at a point,
and $\Sigma_2$ is the quadratic cone with equation $x^2+y^2+z^2$ is
$\CC P^3$ blown up at the node. In both cases the fibration
 is given by the 
extension  of the projection from the blown-up point to a 
line (if
$n=1$) or a hyperplane section (if $n=2$) 
which does not pass through the blown-up point.

The group 
$\mbox{Pic}(\Sigma_n)=\mbox{H}_2(\Sigma_n,\mathbb Z)$
is isomorphic to $\mathbb Z\times\mathbb Z$ and is generated by the
classes of $B_n$ and  $F_n$. Note that we have
$E_n=B_n-nF_n$ in $\mbox{H}_2(\Sigma_n,\mathbb Z).$ 
An algebraic curve $S$
in $\Sigma_n$ is 
said to be of \textit{bidegree} $(a,b)$ if it realizes the homology
class $aB_n+bF_n$ in $\mbox{H}_2(\Sigma_n,\mathbb Z)$.

The surface $\Sigma_n$ is a projective toric surface which
can be obtained by taking two copies of $\CC\times \CC P^1$
 glued
by 
the biholomorphism
$$ 
\begin{array}{ccc}
\CC^*\times \CC P^1 &\longrightarrow & \CC^*\times \CC P^1
\\ (x_1,y_1) &\longmapsto & (\frac{1}{x_1},\frac{y_1}{z_1^n})
\end{array}
$$

\begin{figure}[ht]
\centering
\begin{tabular}{ccc}
\includegraphics[width=5cm, angle=0]{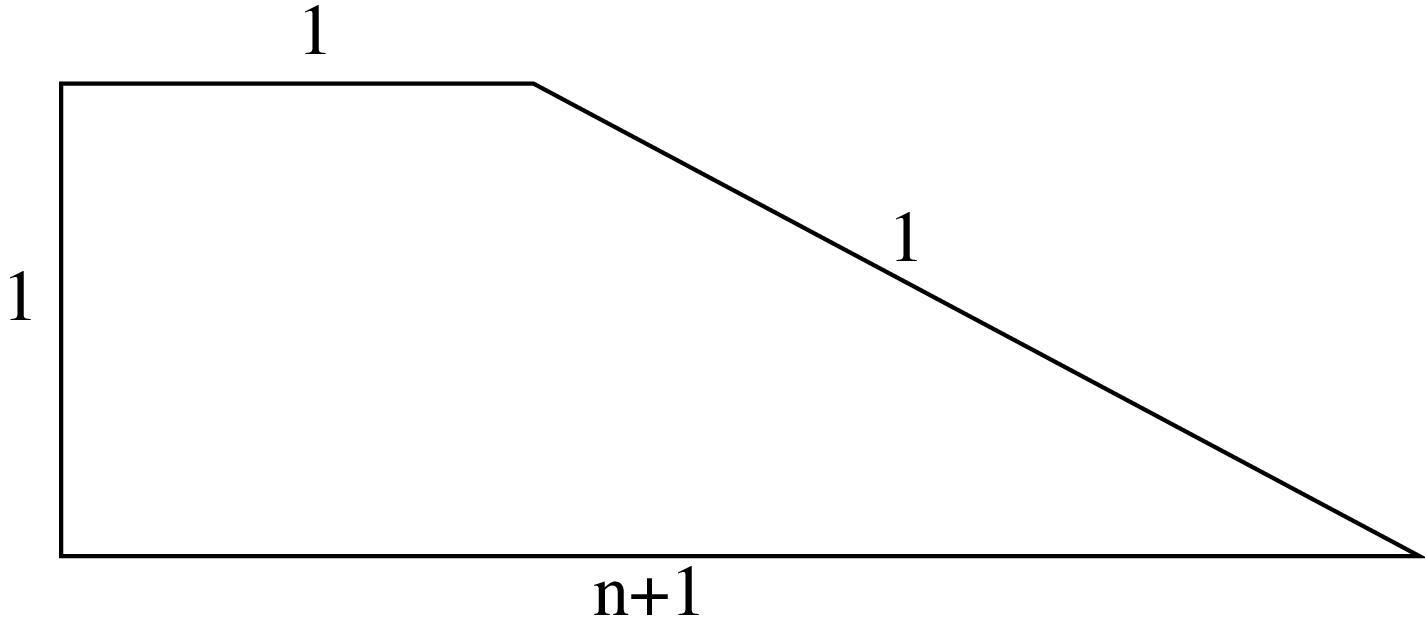}& \hspace{5ex} &
\includegraphics[width=5cm, angle=0]{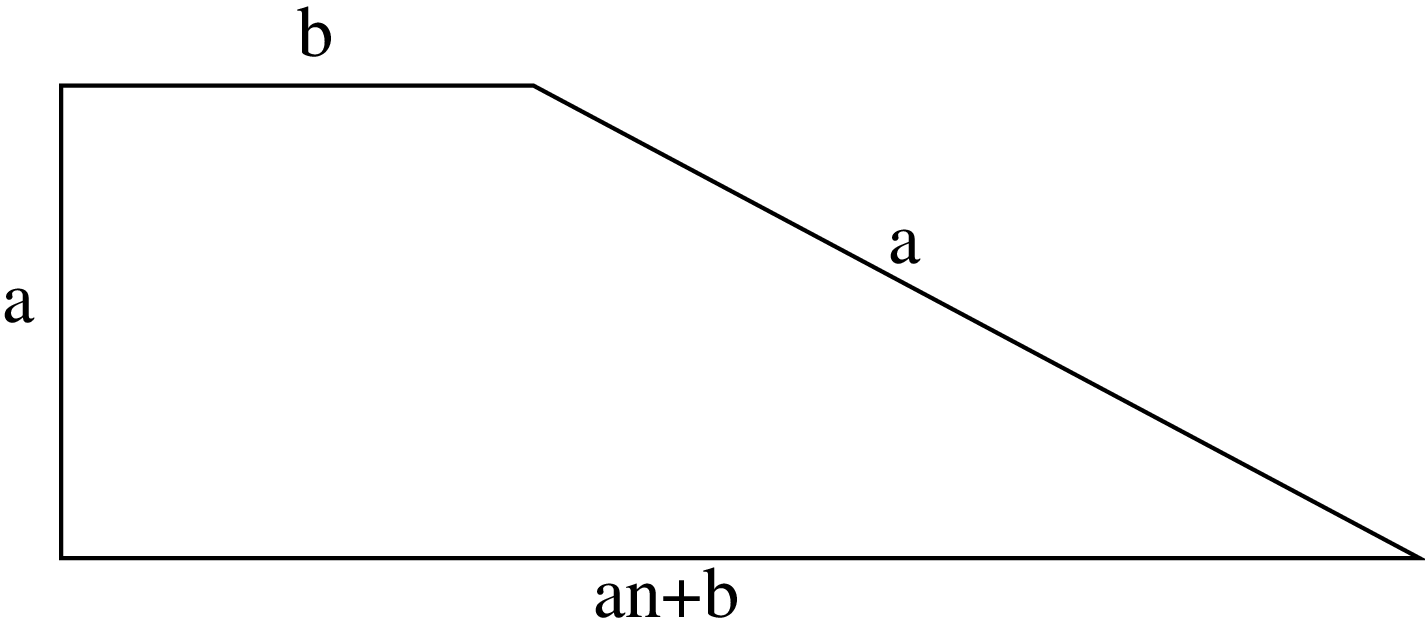}
\\
\\ a) A polygon defining $\Sigma_n$. && b) The Newton
    polygon of a curve \\ & & of bidegree $(a,b)$ in $\Sigma_n$. 
\end{tabular}
\caption{Polygons of Hirzebruch surfaces.}
\label{toric polygon Sn}
\end{figure}

The coordinate system $(x_1,y_1)$ 
in the
first chart
is called \textit{standard}.
 The surface $\Sigma_n$ is the toric surface defined by
the polygon depicted in Figure \ref{toric polygon Sn}a (the number
labeling an edge corresponds to its integer length).
If $S$ is a curve of bidegree $(a,b)$ in $\Sigma_n$ then 
its Newton polygon 
 in a standard
coordinate
system
lies inside the trapeze with vertices
$(0,0)$, $(0,a)$,
$(b,a)$, and
$(an+b,0)$ (see Figure \ref{toric polygon Sn}b), with equality
if $S$ is generic and $a,b\ge 0$.
A curve of class  $B_n$,  $E_n$, or $F_n$
is defined by the equation
$y_1=P(x_1)$, $y_1=+\infty$, or $x_1=c$, respectively, where $P(x)$ is a complex
polynomial of degree at most $n$ and $c\in\CC P^1$.

\subsection{Toric relative  invariants of Hirzebruch
  surfaces.}\label{subsec-enum-inv}
Since we work with curves in toric surfaces having possibly
non-transversal intersection with toric divisors, it is more suitable to
deal with
Newton fans rather than Newton polygons.
We will also use Newton fans in $\ZZ^3$ in Sections \ref{sec-tropenum}
and \ref{sec-corres}.

\begin{definition}[Newton fans]\label{def-Newtonfan}
A \textit{Newton fan} is a multiset $\delta=\{v_1,\ldots,v_k\}$ of vectors $v_i\in \ZZ^r$ satisfying $$\sum_{i=1}^k v_i=0.$$ 
The positive integer $w_i=\gcd(v_{i1}, \ldots,v_{ir})$ (resp. the vector
$\frac{1}{w_i}v_i$) is called the \textit{weight} (resp. the
\textit{primitive direction}) of $v_i$.
We will use the notation $$\delta=\{v_1^{m_1},\ldots, v_k^{m_k}\}$$ to indicate
that the vector $v_i$ appears $m_i$ times in $\delta$.
If $\delta= \{(\alpha_1,\beta_1),\ldots,(\alpha_k,\beta_k)\}$ is a Newton fan of vectors in $\ZZ^2$, one can
construct the \emph{dual polygon} $\Pi_\delta$ in $\RR^2$ in the following
way: for each primitive integer direction $(\alpha,\beta)$ in $\delta$,
we consider the vector $w(-\beta,\alpha)$, where $w$ is the sum of the
weights of all vectors in $\delta$ with primitive integer direction
$(\alpha,\beta)$;
$\Pi_\delta$ is the unique (up to translation) polygon whose
oriented edges (the orientation is induced by the usual orientation of
$\RR^2$) are exactly the vectors $w(-\beta,\alpha)$.

\end{definition}

To any complex algebraic curve $S$ in $(\CC^*)^n$, we may associate its
Newton fan $\delta_S$ as follows: consider the toric compactification
$\mbox{Tor}(\Pi_S)$ of $(\CC^*)^n$ given by a polytope $\Pi_S$ such
that $S$ does not 
intersect boundary components of $\mbox{Tor}(\Pi_S)$ of codimension two or more. Then each puncture $p$ of $S$ corresponds to a facet
$\gamma$ of
 $\Pi_S$. 
 We associate to $p$ the element
$w_pv_p$ where $v_p$ is the primitive normal vector to $\gamma$ oriented
outward $\Pi_S$, and $w_p$ is the order of contact at $p$ of $S$
with the toric divisor corresponding to $\gamma$ in the toric variety
 $\mbox{Tor}(\Pi_S)$. 
The choice of $\Pi_S$ is clearly not unique
however $\delta_S$ does not depend on this choice.
Note that if $n=2$, one may choose $\Pi_S$ equal
to $\Pi_{\delta_S}$.

The enumerative invariants of $\Sigma_n$ we are mainly interested in here
 are defined by counting curves with
prescribed intersection profile with the toric divisors of
$\mbox{Tor}(\Sigma_n)$, and hence correspond to  \emph{toric relative
  invariants} of the surface $\Sigma_n$.
\begin{definition}[Toric relative invariants of $\Sigma_n$]\label{def-toricrelinv}
Given a Newton fan $\delta$ in $\ZZ^2$ and an integer $\chi\in\ZZ$, the number
of 
algebraic curves (resp. irreducible
algebraic curves) in $\mbox{Tor}(\Pi_\delta)$ 
 with Newton fan $\delta$, whose
normalization has Euler
characteristic $\chi$, and  passing
through a generic configuration $\omega$ of $\# \delta -\frac{\chi}{2}$
points  does not depend on $\omega$; we denote this
number by $N_\chi(\delta)$ (resp. $N^{irr}_\chi(\delta)$).

Alternatively $N_\chi(\delta)$ (resp. $N^{irr}_\chi(\delta)$)
is the number of algebraic curves (resp. irreducible
algebraic curves) in $(\CC^*)^2$ 
 with Newton fan $\delta$, whose
normalization has Euler
characteristic $\chi-\#\delta$, and  passing
through a generic configuration  of $\# \delta -\frac{\chi}{2}$
points in $(\CC^*)^2$.

\end{definition}

Remember that if one defines the genus of a compact non-singular 
reducible complex algebraic
curve $S=S_1\cup\ldots\cup S_k$, where $S_i$ is irreducible, as 
$g(S)=g(S_1)+\ldots g(S_k) +1 -k$, then one has
$\chi(S)=2-2g(S)$.

\begin{example}\label{ex faninsigman}
In this text we  consider complex curves in 
$(\CC^*)^2$ with Newton fans of
the form
$$\delta=\{(1,n)^a, (0,-1)^{an+b},  (-\alpha_1,\beta_1),\ldots, (-\alpha_k,\beta_k) \} $$
with $\alpha_i,\beta_i\ge 0$ for
all $i$. A complex curve with such a Newton fan can naturally be seen
as a curve of bidegree $(a,b)$ in $\Sigma_n$ with a singularity at the
point $(0,\infty)$ in the standard coordinates corresponding to $\Pi_\delta$. In particular, a generic
 algebraic curve in $\Sigma_n$ of bidegree $(a,b)$
 will have the Newton fan
$$\delta(a,b,n)=\{(1,n)^a, (0,-1)^{an+b}, (-1,0)^a, (0,1)^b \} $$
in  standard  coordinates.

We give below  some well known enumerative invariants of Hirzebruch surfaces:
$$N_2(\delta(1,0,n))= N_2(\delta(0,1,n))=1\quad \mbox{and}\quad
 N_2(\delta(2,2,0))=12.$$
\end{example}

\subsection{$(1,1)$-relative  invariants of Hirzebruch
  surfaces.}\label{subsec-enum-inv-rel}

Theorem \ref{main formula} involves some additional relative
enumerative
invariants of Hirzebruch surfaces $\Sigma_n$. We give their
definition below, and refer to Sections \ref{sec:inter trop numbers} 
and \ref{sec:inter trop numbers corres}
for more details about
those invariants and their tropical computation.
Given two algebraic curves $S$ and $S'$ in an algebraic surface
intersecting in finitely many points,
 we denote  by $S\circ S'$ the
 intersection number of $S$ and $S'$, and 
by $(S\circ S')_p$ the intersection multiplicity of $S$ and $S'$ at a point $p$.

Let
$S_0$ be a curve of
bidegree $(1,1)$ in $\Sigma_n$ and let $p_0\in S_0\setminus E_n$ be a point. We
denote by $F_0$ the unique curve of bidegree $(0,1)$ passing through
$p_0$ and we choose a (non-standard) 
local system of coordinates $(x_0,y_0)$ on $\Sigma_n$ at $p_0$ such
that 
$S_0$ has local equation $y_0=0$ and $F_0$ has local equation $x_0=0$. 
Given two integer numbers $d_1$ and $d_2$, a \emph{$(d_1,d_2)$-germ}
at $p_0$ is a curve in $\Sigma_n$ with local equation $x_0^{d'_2}
+cy_0^{d'_1}=0 $ with $c\in\CC^*$, $d'_1=\frac{d_1}{\gcd(d_1,d_2)}$, and
$d'_2=\frac{d_2}{\gcd(d_1,d_2)}$.
Let  $D$ be a
local branch at $p_0$ of a reduced algebraic curve in $\Sigma_n$ 
containing neither $S_0$ nor $F_0$ as an irreducible component.
Denote by $d_{S_0}$ (resp. $d_{F_0}$)  the local intersection
multiplicity of $D$ with $S_0$ (resp. $F_0$) at $p_0$.
Then  there exists a unique 
$(d_{S_0},d_{F_0})$-germ at $p_0$ whose intersection multiplicity with
$D$ at $p_0$ is maximal.  We call this curve the
\emph{tangent germ} of $D$ at $p_0$.

\begin{definition}[(1-1)-relative enumerative invariants of
    $\Sigma_n$]\label{def:rel inv}
Let $n\ge 0,$ and  $u, \alpha_1,\ldots  \alpha_r,
d_{1},\ldots d_{r+s}>0 $ be integer numbers such that
$\alpha_i\le d_i-1$ and $\sum_{i=1}^{r+s}d_i=u(n+1)$,  and let $S_0$, $F_0$, and
$p_0$ as above. 
Choose a configuration $\omega=\{ G_1,\ldots, G_r, p_1,\ldots,p_s\}$
of $s$ distinct point $p_1,\ldots, p_s$ in $S_0\setminus \left(E_n\cup\{p_0\} \right)$, and 
$r$ distinct germs at $p_0$  such that $G_i$ is a
$(d_i+\alpha_i(n+1),\alpha_i)$-germ at 
$p_0$. 
We denote by $\mathcal S_{S_0,p_0}(\omega)$ the set of all algebraic curves
$S$ in
$\Sigma_n$ of bidegree $(a,0)$ such that
\begin{itemize}
\item $S$ has a smooth branch tangent to $S_0$ at $p_i$ with
  intersection multiplicity $d_{r+i}$ for $i=1,\ldots, s$;
\item $S$ has exactly $r$ branches $D_1,\ldots,D_r$ at $p_0$, $(D_i\circ S_0)_{p_0}=d_i+\alpha_i(n+1)$ and $(D_i\circ
  F_0)_{p_0}=\alpha_i$; 
\item $G_i$ is the tangent germ of $D_i$ at $p_0$;
\item $S$ has $u$ connected components, whose normalization are 
all rational;
\item each connected component of $S$ intersect $F_0\setminus\{p_0\}$
  in a single point, with intersection multiplicity $1$.
\end{itemize}

For $\omega$ generic, the set $\mathcal S_{S_0,p_0}(\omega)$ is finite
and its cardinal is independent of $S_0$, $p_0$, and $\omega$. We denote it by
$\mathcal N(u,n,d,\alpha)$. Note that it is independent
of the ordering of  the pairs $(d_{1},\alpha_1),\ldots,(d_{r},\alpha_r)$ and of
$d_{r+1},\ldots, d_{r+s}$. If $r=0$, we write $\alpha=0$.
\end{definition}

\begin{remark}\label{rem:finitelt many N(u)}
Given  fixed $n$ and $u$, there clearly exist finitely many choices for
$d$ and $\alpha$. Furthermore the computation of 
 $\mathcal N(u,n,d,\alpha)$ reduces to the computation of all
finitely many possible $\mathcal N(1,n,d',\alpha')$, and to the
combinatorial enumeration 
of how elements of $\omega$ can be
distributed among the $u$ irreducible components of $S$.
\end{remark}

\begin{lemma}\label{lem:sum d_i}
 We use the notation of Definition \ref{def:rel inv}. 
If $\widetilde S$ is an irreducible component of an element $S$
of $\mathcal S_{S_0,p_0}(\omega)$, then 
$$\sum_{D_i\ \mbox{branch of}\ \widetilde S} d_i = n+1 $$
and all intersection points of $S$ with $S_0$ are exactly the points
$p_0,p_1,\ldots,p_s$.
\end{lemma}
\begin{proof}
The curve $\widetilde S$ has bidegree $(\widetilde a,0)$ in
$\Sigma_n$, so $\widetilde S\circ S_0 = \widetilde a(n+1)$ and
$\widetilde S\circ F_0 = \widetilde a$. By the hypothesis we
have
$\widetilde S\circ S_0 \ge \sum d_i + (n+1)\sum \alpha_i$ and
$\widetilde S\circ F_0=\sum \alpha_i +1 $, where the sums are taken
over integers $i$ such that $D_i$ is a branch of $\widetilde S$.
So  we deduce that $\sum d_i \le n+1$. But since $\sum_{i=1}^{r+s} d_i
= u(n+1)$, the latter inequality is in fact an equality.
\end{proof}

We show in Corollary \ref{cor:other def of N} in 
Section \ref{sec:inter trop numbers corres}
that $\mathcal 
N(u,n,d,\alpha)$ has an equivalent definition in terms of the enumeration
of algebraic
curves in $\ctor$  with a given Newton fan $\delta$ and a fixed multiple point
on the toric boundary of $\mbox{Tor}(\Pi_\delta)$.
We give below the values of $\mathcal N(u,n,d,\alpha)$
in the cases $n=0,1,2$,
for which we will
specialize Theorem \ref{main formula}.
Recall that $k!!=k\cdot (k-2)\cdot (k-4)\cdot \ldots 3\cdot 1$ if $k$
is odd, and that  multinomial coefficients 
are defined by
$$ \binom{t}{t_1,\ldots,t_k}=\binom{t}{t_k}\binom{t-t_k}{t_1,\ldots,t_{k-1}}. $$

\begin{example}\label{exa simple trop enum}
We 
suppose that $d_{1}\ge\ldots\ge d_{r}$
and $d_{r+1}\ge\ldots\ge d_{r+s}$. 
The following values for $n=0,1$ and $2$ are computed in Lemmas
\ref{lem:compute0},  \ref{lem:compute1}, and \ref{lem:compute2}: 
$$\mathcal N(u, 0,d,\alpha)= \left\{
\begin{array}{cl}
  1& \text{if} \ d=(1^u) \text{ and }  \alpha=0 
\\  0 & \text{otherwise,}
\end{array}\right. \quad \quad \quad \quad \quad \quad \quad \quad \quad$$

$$\mathcal N (u, 1, d,\alpha)=\left\{
\begin{array}{cl}
(2t_a-1)!!& \mbox{if } d=(2^{u-t_a},1^{2t_a})  \mbox{ and
}\alpha=(1^{u-t_a-t_b}) 
\\   0 & \text{otherwise,}
\end{array}\right. $$

$$\mathcal N (u, 2, d,\alpha)=\left\{
\begin{array}{l}
 3^{t_b}\cdot
  \binom{3t_d+t_c+t_e}{t_c, t_e, 3, \ldots,3}\frac{t_c!t_e!}{t_d!}

\\ \\ \quad   \mbox{if }  d=(3^{t_a+t_b},2^{t_c},3^{t_f},2^{t_e}, 1^{3t_d+t_c+t_e}) 
\mbox{ and }\alpha=(2^{t_a},1^{t_c+t_b})
\\ \quad   \mbox{with } t_a+t_b+t_c+t_f+t_e+t_d=u
\\ \\
\\   0 \quad \text{otherwise}
\end{array}\right. $$
where there are $t_d$ copies of $3$ in the above multinomial coefficient.

\end{example}

\section{Tropical enumerative geometry of Hirzebruch surfaces}\label{sec:basic geotrop}

\subsection{Basics in tropical geometry}
Here we briefly recall standard definitions 
 in tropical geometry in order to fix notations we
 use
in the rest of paper. For a more comprehensive introduction to the
subject the reader may refer for example 
to 
 \cite{BIT},  \cite{Gath1}, \cite{Ma08},
  \cite{Mi03}, \cite{Mi06},   \cite{St2} and references therein.

\subsubsection{Tropical curves and morphisms}
We start by defining tropical curves and their morphisms to $\RR^n$.

Given a (non-necessarily compact) topologically complete metric graph $C$, we denote
by $\Ed(C)$ the set of edges of $C$, by $\Ve(C)$ the set of vertices
of $C$, and by $\Ed^\infty(C)$ the set of non-compact 
edges of $C$. Elements of $\Ed^\infty(C)$ are called \emph{ends}
of $C$, and elements of $\Ed^0(C)=\Ed(C)\setminus\Ed^\infty(C)$ are called
\emph{bounded edges}. The valency of a vertex $V$ is the number of
edges of $C$ adjacent to $V$, and is denoted  by $\val(V)$.

A \emph{degenerate metric graph} is the data of a graph $C$ equipped
with a complete metric on $$C\setminus \left(\bigcup_{e\in E} e \right),$$
called the non-degenerate part of $C$,
where $E$ is some subset of $\Ed^0(C)$. We write $l(e)\in \RR_{>0}$ for the length of the bounded non-degenerate edge $e$.
Edges in $E$ are called
\emph{degenerate}, and their length $l(e)$ is defined to be 0. A connected
component of the union of degenerate edges of $C$ is called a
\emph{degenerate component} of $C$. The metric graph $C$ is said to be
\emph{tamely degenerate} if $C$ does not have
a non-degenerate edge whose adjacent vertices are both contained in
 the same degenerate component.

\begin{definition}[Abstract tropical curves]
An \emph{abstract (punctured) tropical curve} is a tamely degenerate metric
graph without degenerate loop-edges, equipped with  a
function 
$$\begin{array}{ccc}
\Ve(C)&\longrightarrow & \ZZ_{\geq 0}
\\ V &\longmapsto & g_V
\end{array} $$
 satisfying the stability condition. 
The integer $g_V\in \ZZ_{\geq 0}$ is called  the
\emph{genus} of the vertex $V$. 
 The stability condition states that vertices of genus $0$ are at
 least $3$-valent, and vertices of genus $1$ are at least $1$-valent.
 Two abstract tropical curves
are isomorphic (and will from now on be identified) if there exists an
homeomorphism between them, restricting on an
isometry on the non-degenerate parts, that respects the genus function on the vertices.

The tropical curve $C$ is said to be \emph{irreducible} if it is connected.

The \emph{topological Euler characteristic} of $C$ is
$\chi(C)=b_0(C)-b_1(C)$ where $b_i(C)$ is the $i$th Betti number of $C$.

The \emph{tropical Euler characteristic} of $C$ is
$\chi_{\trop}(C)=\chi(C) -  \sum_{V \in\Ve(C)} g_V$.

The \emph{genus} of an abstract tropical curve $C$ is defined to be 
$g(C)=1-\chi_{\trop}(C)$.

We say that an abstract tropical curve is \emph{explicit} if $g_V=0$ for all vertices $V$.

The \emph{combinatorial type} of an abstract tropical curve is the
homeomorphism class of $C$ together with the genus function and the
degenerate edges, i.e.\ we drop the information about the positive lengths.
\end{definition}

\begin{remark}
The consideration of tropical curves with degenerate edges has already
been suggested in the literature (see for example {\cite[Definition
    6.2]{Br9}}). However 
to our knowledge
 they were not rigorously introduced and used
 up to now since 
tropical curves without degenerate edges suffice
for many applications of tropical geometry, for example for the
enumeration of curves in toric varieties.
Degenerate edges turn out to be necessary when considering tropical
morphisms to tropical varieties different for $\RR^n$, 
where it is essential to
distinguish  the set-theoretic inclusion from the tropical-theoretic
inclusion.
We refer to Section \ref{sec:phase trop} for more details.
\end{remark}

\begin{remark}
We recover the classical formula
$\chi(C)=\#\Ve(C) - \#\Ed(C)  $ only in the case when $C$ is
compact. When $C$ is not compact, we have 
$$\chi(C)=\#\Ve(C) - \#\Ed^0(C)=\#\Ve(\overline C) - \#\Ed(\overline C),  $$
where $\overline C$ is the compact tropical curve obtained from $C$
by adding one vertex to each of its ends.

Note that $\chi(C)=\chi_{\trop}(C)$ if and only if $C$ is explicit,
and  that
$g(C)=b_1(C) +  \sum_{V \in\Ve(C)} g_V + 1-b_0(C)$.
\end{remark}

\begin{definition}[Tropical morphisms of curves to $\RR^n$]\label{def:trop morph}
A  tropical morphism is a continuous map $h:C\to \RR^n$ from an abstract tropical curve $C$ satisfying the
following conditions:

\begin{itemize}
 
\item On each non-degenerate edge $e$ the map $h$ is integer affine
  linear, i.e.\ of the 
    form $ h_{|e}(t) = a + t \cdot v $ with $ a \in \RR^n $ and $ v \in
    \ZZ^n $ once we have identify $e$ with an interval of $\RR$ by an
    isometry.
If
    $ V \in \partial e $ and we parameterize the edge $e$ starting at $V$, the
    vector $v$ in the above equation will be denoted $ v(V,e) $ and called the
    \emph {direction vector} of $e$ starting at $V$. If $V$ is understood from
    the context (e.g.\ in case $e$ is an end, having only one adjacent vertex)
    we will also write $ v(e) $ instead of $ v(V,e) $. The greatest
    common divisor of the entries of $v(e)$ is called the
    \emph{weight} of $e$ for $h$.

\item each degenerate edge $e$  adjacent to the vertices $V$ and
  $V'$ is equipped with two vectors $ v(V,e) =-v(V',e)$ in $\ZZ^n$;
moreover $e$ is mapped to a point by $h$
  (i.e. the restriction of $h$ on $e$ can be thought as an
  \emph{infinitesimal} integer affine linear map).

\item At each vertex $V$ the \emph{balancing condition} is satisfied, i.e.\ $$\sum_{e\in\Ed(C); \;V \in \partial e} v(V,e)=0.$$

\end{itemize}
Note that degenerate edges of $C$ contribute to the balancing condition. 
   Two tropical morphisms $h:C\to\RR^n$ and $\tilde
    h:\tilde C\to \RR^n$  are called
    \emph{isomorphic} (and will from now on be identified) if there is
    an isomorphism of tropical curves
    $ \phi: C \to \tilde C $  such that $ \tilde h \circ \phi = h $.

The Newton fan $\delta$ of a tropical morphism $h:C\to \RR^n$
 is the multiset $ \delta = \{ v(e) \}_{e\in\Ed^\infty(C)} $.

The \emph{combinatorial type} of $h:C\to \RR^n$ 
is the combinatorial type of the underlying abstract curve $C$ together
with the direction vectors of all edges.
\end{definition}

\subsubsection{Tropical morphisms of curves to a tropical surface}

The proof of Theorem \ref{main formula} works by enumerating of
 tropical curves in a
tropical surface $X\subset \RR^3$ (see Remark \ref{rem:justification}
for a justification). We first recall the definition of a
tropical hypersurface of $\RR^n$, and then  define 
tropical morphisms $h:C\to X$ for a particular surface $X$.

Recall that the tropical operations on $\RR$, denoted by $\tg +\td$
and $\tg\times \td$, are defined by
$$\tg a+b \td=\max(a,b)\quad \text{and}\quad \tg a\times b\td=a+b. $$
A tropical polynomial is then a piecewise-affine function
$$\begin{array}{ccccc}
P: & \RR^n&\longrightarrow & \RR
\\ & x&\longmapsto & \tg \sum_i a_ix^i\td&= \max_i(a_i + \langle x,i\rangle) ,
\end{array} $$
where $\langle , \rangle$ denotes the standard Euclidean product on
$\RR^n$.

\begin{definition}[Tropical hypersurfaces]\label{def trop hyper}
The tropical hypersurface $Z(P)$ defined by a tropical polynomial 
$P(x)=\tg \sum_i a_ix^i\td$ in $n$
variables is the subset of $\RR^n$ where the value of $P(x)$ is given by (at least)
two distinct monomials of $P$.
\end{definition}
In other words, a points $x_0$ is in $Z(P)$ if and only if there exist
two distinct $i$ and $j$ in $\ZZ^n$ such that
$P(x_0)=a_i + \langle x,i\rangle=a_j + \langle x,j\rangle$.

\begin{remark}
We work with a simplified definition of a tropical hypersurface: we neglect the weights that one usually assigns to the facets of the polyhedral complex $Z(P)$.
We do not discuss these weights here, since for
the particular tropical surface $X$ considered in the rest of this
paper they do not play any role.
\end{remark}

\begin{example}\label{const-surface}
In the whole paper, we denote by $P_X$ the tropical polynomial $P_X(x,y,z)=\tg x+y+z\td$ and by $X$ the tropical surface $Z(P_X)$ in $\RR^3$. The surface $X$  consists of three
2-dimensional cells, $\sigma_1=\{x=y\geq z\}$, $\sigma_2=\{x=z\geq
y\}$ and $\sigma_3=\{y=z\geq x\}$, that meet in the line
$L=\{x=y=z\}=\RR(1,1,1)$ (see Figure \ref{fig X2}).

\end{example}

Now we define tropical morphisms $h:C\to X$. As it has
 already been mentioned in the introduction, the tropical inclusion
is more subtle  than
the set theoretic inclusion. 
In the simple situation we deal with in
this paper, i.e. $X$ is just
made of three faces meeting along the line $L$,  there are only two
extra conditions we have to impose on $h$ to be tropically contained
in $X$: the direction vectors of degenerate edges should be contained in
the union of the linear spans of the faces of $X$, and $h$ should
satisfy 
the so-called
Riemann-Hurwitz condition.
 The necessity of this latter condition has first
been observed by the first author and Mikhalkin (\cite{Br12}, see
also \cite{Br13}, \cite{ABBR13a}, and 
\cite{ABBR13b}). More  conditions that have to be required
for a tropical morphism $h:C\to \widetilde X$ to a more general tropical surface $\widetilde X$
can be found in \cite{Br17} and \cite{GathSchW}.

The above-mentioned Riemann-Hurwitz condition is based on tropical
intersection theory. We first  recall  basic facts and definitions
in tropical intersection theory restricting to the
information we need for our purposes.
We use notations from Example \ref{const-surface} in next definitions.

We define $\tilde \sigma_1=\{x=y\}$, $\tilde \sigma_2=\{x=z\}$,
$\tilde \sigma_3=\{y=z\}$, and 
$T_LX=\tilde \sigma_1 \cup\tilde \sigma_2 \cup \tilde \sigma_3$.

\begin{definition}[Tropical premorphisms to $X$]\label{def:premorphism}
A \emph{tropical premorphism} $h:C\to X$ is a tropical morphism
$h:C\to\RR^3$ such that 
\begin{itemize}
\item $h(C)\subset X$;
\item for any edge $e$ adjacent to the vertex $V$, one has
  $v(V,e)\in T_LX$; moreover if  $v(V,e)\in \tilde \sigma_i\setminus X $
   one has $v(V,e')\in \tilde \sigma_i$ for all edges $e'$
  adjacent to $V$.

\end{itemize}

An edge $e$ of $C$ is said to
be \emph{tropically mapped to $L$} if $h(e)\subset L$ and $v(V,e)$ is
parallel to $(1,1,1)$.
\end{definition}
Note that the second condition above is non-empty only for degenerate
edges. 
Obviously an edge $e$ is tropically
mapped to $L$ if $e$ is non-degenerate and $h(e)\subset L$. On the
other hand, if $e$ is degenerate then the fact that $f(e)\subset L$
does not imply that $v(V,e)$ is parallel to $(1,1,1)$.

\begin{definition}[Intersection multiplicities] \label{def int}
Let $h:C\to X$ be a tropical premorphism, and let $V$ be a vertex of
$C$.
If $h$ 
does not map a neighborhood of $V$ entirely in some plane $\tilde\sigma_i$, the \emph{intersection multiplicity} of 
$C$ with $L$ at $V$, denoted by $d_V$, is defined as
$$d_V
=\sum_{V \in \partial e\; \text{and}\;h(e)\subset\sigma_i} |P_X(v(V,e)) + P_X(-v(V,e))|$$
for any choice of face $\sigma_i$ of $X$. 
Otherwise we set $d_V=0$. 
 
A vertex $V$ of
$C$ is said to be \emph{tropically mapped to $L$} if either $d_V>0$ or
each edge
adjacent to $V$ is tropically mapped to $L$. 

The \emph{overvalency} of $V$ is defined by 
\begin{displaymath}
 \ov_V:=k_V- d_V-2+2g_V,
\end{displaymath}
 where $k_V$ is the number of edges of $C$ adjacent to $V$ and not
tropically
mapped 
to $L$  by $h$.

The \emph{total intersection number} 
of $C$ with $L$, denoted by $C\circ L$, 
 is defined as
$$C\circ L=\sum_{V\in\Ve(C) }d_V. $$ 
\end{definition}

\begin{remark}\label{total int}
The intersection multiplicity  $d_V$ does not depend on the choice of the face $\sigma_i$
(see \cite{AR07}, \cite{Shaw1}). More precisely, if we choose
$\sigma_i=\sigma_1$, then 
$$v(V,e)=(x,x,z) \quad\text{and}\quad |P_X(v(V,e)) + P_X(-v(V,e))|=|x-z|.$$

The total intersection number of $C$ and $L$
 depends only on the Newton fan of $h:C\to X$: 
the Newton fan equals the recession fan of the (embedded) tropical curve $h(C)$,
i.e.\ the fan obtained by shrinking all bounded edges of $h(C)$
to length $0$; the total intersection number of $C$ with $L$ equals
the tropical intersection of $h(C)$ with $X$ in $\RR^3$; since $h(C)$
is rationally equivalent to its Newton fan, this intersection depends
only on the Newton fan. 
\end{remark}

\begin{definition}[Tropical morphisms to $X$]\label{def:morphX}
A tropical premorphism $h:C\to X$ is a tropical morphism to $X$ if
 $\ov_V\geq 0$ for any vertex $V$ of $C$ with $d_V>0$.
\end{definition}
This extra condition for a premorphism to be a morphism is a consequence of the
Riemann-Hurwitz formula in complex geometry (see Section
\ref{sec:phase trop} for more details) and is usually referred to as the
\emph{Riemann-Hurwitz condition}.

\begin{example}
Two examples of tropical premorphisms $h_1:C_1\to X$ and $h_2:C_2\to X$, 
where $C_{1}$ and $C_2$ are two
tropical curves made of one vertex and 3 ends, are depicted in
Figure \ref{exmor}a. 
The integer close the the vertex of $C_i$ denotes its genus, and the
vector close to the image of an edge $e$ is the vector $v(e)$.
In each case, $C_i \circ L=2$.
The tropical curve $C_1$ has genus $0$, and $C_1$
has genus 1, so $h_1$ is not a tropical morphism to $X$, while $h_2$
is.

Figure \ref{exmor}b shows an example of a tropical morphism
$h_3:C_3\to X$, where $C_3$ is a
tropical curve of genus 0 made of one vertex and 4 ends.
In the whole paper we use the same convention
to encode both an abstract tropical
curve and its image by a tropical premorphism  in the same picture: 
 the two edges we draw close to each other in the
$(-1,0,0)$-direction are distinct edges in $C_3$ which have the same
 image by $h_3$. 
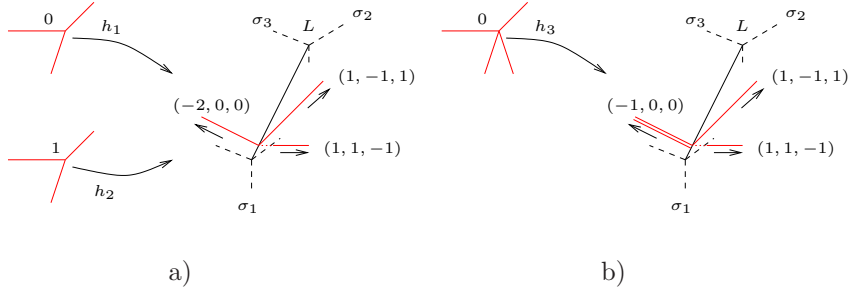
\begin{figure}[ht]
\centering
\begin{tabular}{ccc}
\input{Ex_mor.pstex_t} & \hspace{3ex} &
\input{Ex_mor2.pstex_t}
\\
\\ a)  && b) 
\end{tabular}
\caption{Examples of tropical (pre-)morphisms.}
\label{exmor}
\end{figure}

\end{example}

The \emph{combinatorial type} of a tropical morphism $h:C\to X$ is the
combinatorial type of the underlying abstract tropical  curve together
with the direction vectors for all edges, and together with the set of
vertices and edges that are 
tropically
mapped into $L$.

\begin{notation}\label{not-ci}
Given $h:C\to X$ a tropical morphism, 
 we define $C_i=h^{-1}(\sigma_i)$ 
and $h_i=h_{|C_i}$ for $i=1,2,3$. Note that the projection $(x,y,z)\mapsto (x,y)$
identifies
each face $\sigma_2$ and $\sigma_3$ of $X$ with a half-space in
$\RR^2$, and maps the lattice $\sigma_i\cap \ZZ^3$ isomorphically 
to the lattice $\ZZ^2\subset\RR^2$. Hence using this projection and 
 extending edges of $C$ 
that meet $L$ to ends, 
the map $h_{i}:C_{i}\to\RR^2$ can be understood as a tropical morphism to
the plane for $i=1,2$.
In the same way, the projection $(x,y,z)\mapsto (y,z)$
identifies the face  $\sigma_1$ of $X$ with a half-space in
$\RR^2$, and using this projection the map $h_1:C_1\to\RR^2$ can be
understood as a tropical morphism to 
the plane.
\end{notation}

\subsection{Tropical analogues of $(1,1)$-Relative invariants of $\Sigma_n$}\label{sec:inter trop numbers}
We now define tropical analogues of $(1,1)$-relative invariants and state their correspondence to the numbers 
$\mathcal
N(u,n,d,\alpha)$ of Definition \ref{def:rel inv}
 in Theorem
\ref{thm:corres 2}. We end this section by
computing explicitly those invariants in the cases $n=0,1,$ and $2$.
The proof of Theorem \ref{thm:corres 2}.
is postponed to Section \ref{sec:inter trop numbers corres}.

\subsubsection{Tropical analogues of $\mathcal N(u,n,d,\alpha)$}
If $C_1$ and $C_2$ are two tropical curves in the plane $\RR^2$, we
denote by $C_1\circ C_2$ their tropical intersection number, and by 
$(C_1\circ C_2)_p$  the stable intersection multiplicity of $C_1$ and
$C_2$ at the point $p$ (see for example \cite{St2}, Section 4).

\begin{definition}\label{def:trop invmult}
Fix the following data:
\begin{itemize}
 \item two 
integers  $u\ge 1$ and $n\ge 0$,
\item a tuple $d=(d_1,\ldots,d_{r+s})$ of positive integers satisfying
  $d_i\le n+1$ and $\sum_{i=1}^{r+s} d_i=u(n+1)$,
\item a tuple $\alpha=(\alpha_1,\ldots,\alpha_r)$ of positive
 integers satisfying $\alpha_i\leq d_i -1$ (we write $\alpha=0$ if $r=0$).
\end{itemize}

Suppose in addition that 
$$\frac{d_i}{\alpha_i}\le \frac{d_{i+1}}{\alpha_{i+1}}\quad \forall
i=1,\ldots r-1.$$
Furthermore, 
choose a collection of $r+s$
 points $\omega=(p_1,\ldots, p_{r+s})$ on the line $L=\{x=y\}$ in $\RR^2$, in such a way that the $x$-coordinates of the
 points $p_i$ increase as $i$ increases.

Let $\TT\mathcal S'(u,n,d,\alpha,\omega)$ be the set of tropical morphisms $h:C\to
\RR^2$ passing through $\omega$ (i.e.\ satisfying $\omega\subset h(C)$) such that
\begin{itemize}
\item $C$ is the disjoint union of of $u$ rational tropical curves;

\item for each $i=1,\ldots,r+s$, the curve $C$ has exactly one vertex
  $V_i$ with $h(V_i)=p_i$ such that the intersection multiplicity of $C$ with $L$ at $V_i$  equals
 $d_i$. Furthermore, each vertex $V_i$ is adjacent to $k_i$ ends of direction $(-1,0)$ and to $l_i:=d_i-k_i$ ends of direction $(0,1)$.
In addition, these vertices satisfy:
\begin{itemize}
 \item \emph{for $i=1,\ldots,r$:} $V_i$ is $d_i+2$-valent, in addition
   to the ends described above it is adjacent to one end of direction
   $(\alpha_i, \alpha_i)$. In the following, we sum up all these
   requirements by saying that $V_i$ is a vertex \emph{of type 1} 
(see Figure \ref{fig specvert2}a).
\item \emph{for $i=r+1,\ldots,r+s$:} $V_i$ is $d_i+1$-valent. In the
  following, we sum up all these requirements by saying that $V_i$ is
  a vertex \emph{of type 2} 
(see Figure \ref{fig specvert2}b).
\end{itemize}
\item each connected component of $C$ contains exactly one
 end with direction vector $(1,-n)$;
\item $C$ has exactly $u(n+2) + r$ ends (in particular they all have
  been described above).
\end{itemize}
\begin{figure}[ht]
\centering
 \begin{tabular}{ccc}
\input{specvert1.pstex_t}& \hspace{3ex} &
\input{specvert2.pstex_t}
\\ a) Vertex of type 1 & & b) Vertex of type 2
\end{tabular}
 \caption{Vertex types on $L$ of tropical morphisms in $\TT\mathcal S'(u,n,d,\alpha,\omega)$.}
\label{fig specvert2}
\end{figure}
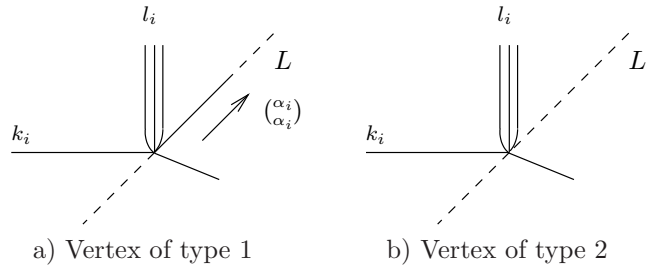

Let $h:C\to \RR^2$ be an element of $\TT\mathcal S'(u,n,d,\alpha,\omega)$.
For each vertex $V$ of $C$ not mapped to $L$, we choose two adjacent
vectors $e_{V,1}$ and $e_{V,2}$, and we define the multiplicity of $h:C\to \RR^2$
as 
$$\mu_h=\prod_{i=1}^{r+s}\binom{k_i+l_i}{k_i}\prod_{V\in\Ve(C),\ V\ne V_i}|\det(v(e_{V,1}),v(e_{V,2}))|. $$

Finally we define the number
$$\TT\mathcal N(u,n,d,\alpha)=\left(\prod_{i=1}^r \frac{1}{\gcd(d_i,\alpha_i)}\right)\cdot
\left(\prod_{i=1}^s \frac{1}{d_{r+i}}\right)\cdot
\sum_{h\in \TT\mathcal S'(u,n,d,\alpha,\omega)}\mu_h.$$
We also set
$$\TT\mathcal N(0,n,d,\alpha)=1 \mbox{ if }d=\alpha=0,\ \mbox{ and } 0
\mbox{ otherwise}.$$ 
\end{definition}
\begin{remark}\label{rem: obs int inv}
It follows from the definition that 
the total intersection number of a component $\widetilde C$ of $C$
with $L$ equals
$\sum_{V_i\in\widetilde C } d_i =n+1$.

Note also that since $0<\alpha_i\le d_i-1$, we have $d_i\ge 2$ for
$i\le r$.
\end{remark}

The next lemma states in particular 
that $\TT\mathcal N(u,n,d,\alpha)$ is indeed finite. The proof follows
 exactly along the lines of \cite{GM051}, so we do not include it here.
Denote by $\mathcal U$ the subset of $L^{r+s}$ consisting of tuple
$(p_1,\ldots, p_{r+s})$ such that $p_i$ has a strictly smaller
$x$-coordinate than $p_j$ if $i<j$, for all $i,j=1\ldots r+s$.
\begin{lemma}
There exists a dense open subset $U$ of $\mathcal U$ such that for any choice
of collection $\omega$ in $U$
satisfying the requirement of Definition
\ref{def:trop invmult} one has
\begin{itemize}
\item the set 
$\TT\mathcal S'(u,n,d,\alpha,\omega)$ is finite;
\item for an element $h:C\to \RR^2$ in $\TT\mathcal S'(u,n,d,\alpha,\omega)$, any
  vertex $V$ of $C$ not mapped to $L$ is 3-valent; furthermore $h$ is
  an embedding in a neighborhood of $V$;
\item the number $\TT\mathcal N(u,n,d,\alpha)$ does not depend on the choice
  of $\omega$ in $U$.
\end{itemize}
\end{lemma}

The independence of the number $\TT\mathcal N(u,n,d,\alpha)$ of $\omega$ is also an immediate corollary of the next theorem,
whose proof is postponed to Section \ref{sec:inter trop numbers corres}.

\begin{theorem}\label{thm:corres 2}
For any $u,n,d,$ and $\alpha$ satisfying the hypothesis of Definition
\ref{def:trop invmult} one has
$$\TT\mathcal N(u,n,d,\alpha)= \mathcal N(u,n,d,\alpha). $$
\end{theorem}

\subsubsection{Examples}
We end this section by computing the numbers $\mathcal N(u,n,d,\alpha)$
for $n=0,1$ and $2$. The cases $n=0$ and $1$ could easily by done
without  the whole machinery of tropical geometry. However we perform
their computation tropically as a warm-up for less straightforward computations.
We first show that there are only two kinds of vertices of type 2, and that
$\mu_h$ only gets non-trivial contribution from connected components of
$C$ containing a vertex of type 1.
We now fix a tuple $\omega=(p_1,\ldots, p_{r+s})$ in $U$, and an element
$h:C\to\RR^2$ of $\TT\mathcal S'(u,n,d,\alpha,\omega)$.

\begin{figure}[ht]
\centering
\begin{tabular}{c}
\input{Simple.pstex_t} 
\end{tabular}
\caption{The element of $\TT \mathcal S'(1,n,d,0,\omega)$.}
\label{fig:easy 2nd kind}
\end{figure}
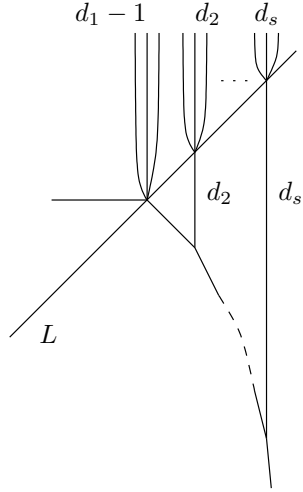

\begin{lemma}\label{lem:easy 2nd kind}
Let $\widetilde C$ be a connected component of $C$ containing an end $e$
with direction vector $(-1,0)$
adjacent to a vertex $V_{r+j}$ of type 2. 
Then $\widetilde C$ does not contain any vertex of type 1,
and $e$ is the only end of $\widetilde C$
with direction vector $(-1,0)$. 
In particular one has
$$\prod_{i\ | \ V_{i}\in \widetilde    C}\frac{1}{d_{i}}\cdot
\prod_{i\ | \ V_{i}\in \widetilde
    C}\binom{k_{i}+l_{i}}{k_{i}}\cdot\prod_{V\in\Ve(\widetilde C),\ V\ne
  V_i}|\det(v(e_1),v(e_2))| = 1. $$
\end{lemma}
If $\alpha=0$ 
(i.e.\ no vertices of type 1) and $u=1$ (i.e.\ $C$ is irreducible)
the curve $C$ is mapped to $\RR^2$ as depicted in Figure
\ref{fig:easy 2nd kind}. 

\begin{proof}
Let $T$ be the smallest subtree of $\widetilde C$ containing the end $e'$
 of $\widetilde C$ with direction vector $(1,-n)$, 
as well as all the vertices of type 2 contained in $\widetilde C$ (see
Figure \ref{fig:orientation} where we depict only what happens in the
half-plane $x>y$).
 Orient the edges in $\widetilde C$ such
 that they point towards the end $e'$.
Assume that $T\neq \widetilde C$.
For each vertex $V_i$ of $\widetilde C$ we denote by $e_i$ the unique
edge of $\widetilde C$ adjacent to $V_i$ and mapped to the half-plane $x>y$.
By definition of the set
 $\TT\mathcal S'(u,n,d,\alpha,\omega)$, the direction vector of the edge
$e_{r+i}$ has $x$-coordinate $k_{r+i}$. 
Any edge of $ \widetilde C$ adjacent to $T$ but not in $T$ connects $T$ to some of the points 
 $p_1,\ldots,p_r$, which have a strictly less
$x$-coordinate than the points $p_{r+1},\ldots,p_{r+s}$. Thus any such edge has a positive
$x$-coordinate. Since the sum of the $x$-coordinates of all such edges and of all edges
$e_{r+i}$ in $ \widetilde C$  equals 
 $1$ we can conclude that there is in fact only one nonzero summand
 which equals one. By assumption, we have $k_{r+j}>1$, so we must have
 $k_{r+j}=1$ and $T=\widetilde C$. The claim follows. The claim about the contribution of such a component to $\TT\mathcal N(u,n,d,\alpha)$ follows since the vertices below $L$ contribute $\prod_{i=2}^s d_i$, the binomial factors are all one except for $V_1$ where we get $d_1$, and we divide by $\prod_{i=1}^s d_i$. 
\begin{figure}[ht]
\centering
\begin{tabular}{c}
\input{tildeC.pstex_t} 
\end{tabular}
\caption{A component $\tilde C$, together with the tree $T$
  in thick.}
\label{fig:orientation}
\end{figure}
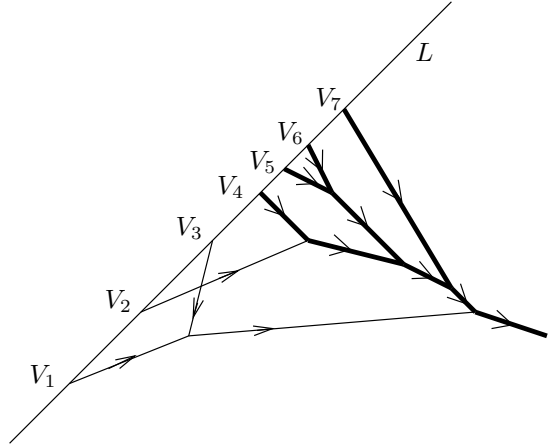
\end{proof}

\begin{lemma}\label{lem:compute0}
We have 
$$\mathcal N(u, 0, d,\alpha)= \left\{
\begin{array}{cl}
  1& \text{if} \ d=(1^u)\ \text{and}\ \alpha=0
\\  0 & \text{otherwise}.
\end{array}\right.$$
In particular, for any morphism in $ \TT\mathcal S'(u,0,(1^u),0,\omega)$
every connected component is mapped to a horizontal line as in Figure
\ref{fig:n=0,1}a. 
\end{lemma}
\begin{proof}
It is enough to prove the lemma in the case $u=1$.
Let $h:C\to \RR^2$ be in $\TT\mathcal S'(1,0,d,\alpha,\omega)$. Remark 
\ref{rem: obs int inv} implies that
$\sum d_i=1$ and $r=0$.  So $s=1$, which
completes the proof by Theorem \ref{thm:corres 2}.
\end{proof}

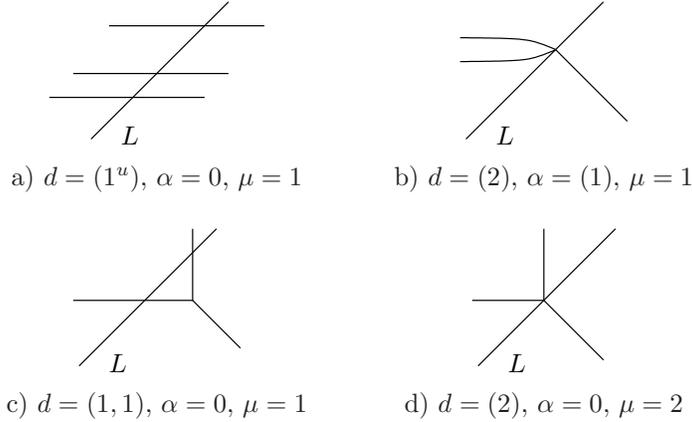
\begin{figure}[ht]
\centering
\begin{tabular}{ccc}
\input{Compute0.pstex_t}& \hspace{3ex}  &
\input{Compute1a.pstex_t}
\\ a) $d=(1^u)$,  $\alpha=0$,  $\mu=1$ &&b)  $d=(2)$, $\alpha=(1)$, $\mu=1$
\\ 
\\
\input{Compute1b.pstex_t}&  \hspace{3ex}  &
\input{Compute1c.pstex_t}
\\ c) $d=(1,1)$, $\alpha=0$, $\mu=1$ &&d) $d=(2)$, $\alpha=0$, $\mu=2$

\end{tabular}
\caption{$\mathcal N (u, n, d,\alpha)$ for $n=0$ and $1$.}
\label{fig:n=0,1}
\end{figure}

\begin{lemma}\label{lem:compute1}
Suppose that $d_{r+1}\ge\ldots\ge d_{r+s}$. 
We have
$$\mathcal N (u, 1, d,\alpha)=\left\{
\begin{array}{cl}
 (2t_a-1)!!& \mbox{if }  d=(2^{u-t_a},1^{2t_a}) 
\mbox{ and }\alpha=(1^{u-t_a-t_b})
\\   0 & \text{otherwise.}
\end{array}\right. $$
\end{lemma}
\begin{proof}
Let $h:C\to \RR^2$ be an element of $\TT\mathcal S'(u,1,d,\alpha,\omega)$.
We study the connected components of $C$, i.e.\ we assume first
that $u=1$, which
implies $\sum d_i =n+1=2$. 

If $r>0$, then by Remark \ref{rem: obs int inv} we have $d=(2)$ and
$\alpha=(1)$. This implies immediately that
 $C$ has only one vertex, and that the only possibility for
the Newton fan of $h$ is $\{(1,-1),(1,1), (-1,0)^2\}$. In this case
$C$ is mapped to $\mathbb R^2$ as depicted in Figure \ref{fig:n=0,1}b
and $\mu_h=1$.

If $r=0$ according to Lemma \ref{lem:easy 2nd kind} the 
Newton fan of $h$ is  $\{(-1,0),(0,1),(1,-1)\}$, the curve 
$C$ is mapped to $\mathbb R^2$ as depicted in Figure \ref{fig:n=0,1}c
or d, and $\mu_h$ equals respectively $1$ and $2$.

\vspace{1ex}
In the general case, every connected component of $C$ is mapped to
$\mathbb R^2$ as depicted in Figures \ref{fig:n=0,1}b, c, and d.
Thus we must have $\alpha=(1^{r})$, and $r$ equals the number of components as in Figure \ref{fig:n=0,1}b. Denote by $t_b$ the number of vertices of type 2 with $d_i=2$ (i.e.\ the number of components as in Figure \ref{fig:n=0,1}d) and by $2t_a$ the number of vertices of type 2 with $d_i=1$ (i.e.\ $t_a$ is the number of components as in Figure \ref{fig:n=0,1}c). Since $u$ is the total number of components we have $r+t_a+t_b=u$. We can thus express $d=(2^r,2^{t_b},1^{2t_a})$ as $d=(2^{u-t_a},1^{2t_a})$ and $\alpha=(1^r)=(1^{u-t_a-t_b})$.

There is a one to one correspondence between $\TT\mathcal
S'(u,1,d,\alpha,\omega)$ and the  partitions into pairs of the set of $p_i$
such that $d_i=1$.
Hence we have $\# \TT\mathcal S'(u,1,d,\alpha,\omega)= (2t_a-1)!!$. 
The multiplicities of the components as in Figure \ref{fig:n=0,1}d exactly cancels with the global factor $\prod \frac{1}{d_{r+i}}$ in the definition of $\TT \mathcal N (u, 1, d,\alpha)$,
thus we have $\TT \mathcal N (u, 1, d,\alpha)=\mathcal N (u, 1, d,\alpha)=  (2t_a-1)!!$ by Theorem \ref{thm:corres 2}. 
\end{proof}

\begin{figure}[ht]
\centering
\begin{tabular}{ccc}
\input{Compute2a.pstex_t}&  \hspace{3ex}  &
\input{Compute2b.pstex_t}
 \\ a) $d=(3)$,  $\alpha=(2)$,  $\mu=1$ &&b)  $d=(3)$, $\alpha=(1)$,
$\mu=3$
\\ \\ \input{Compute2c.pstex_t}&  \hspace{3ex}  &
\input{Compute2f.pstex_t} 
\\  c) $d=(2,1)$,  $\alpha=(1)$,  $\mu=1$ 
 &&d) $d=(1^3)$, $\alpha=0$, $\mu=1$
\\ \\
\input{Compute2e.pstex_t}&  \hspace{3ex}  &
\input{Compute2d.pstex_t}
\\ e) $d=(2,1)$, $\alpha=0$, $\mu=2$ &&f) $d=(3)$, $\alpha=0$, $\mu=3$

\end{tabular}
\caption{Components of curves contributing to $\mathcal N (u, 2, d,\alpha)$.}
\label{fig:n=2}
\end{figure}
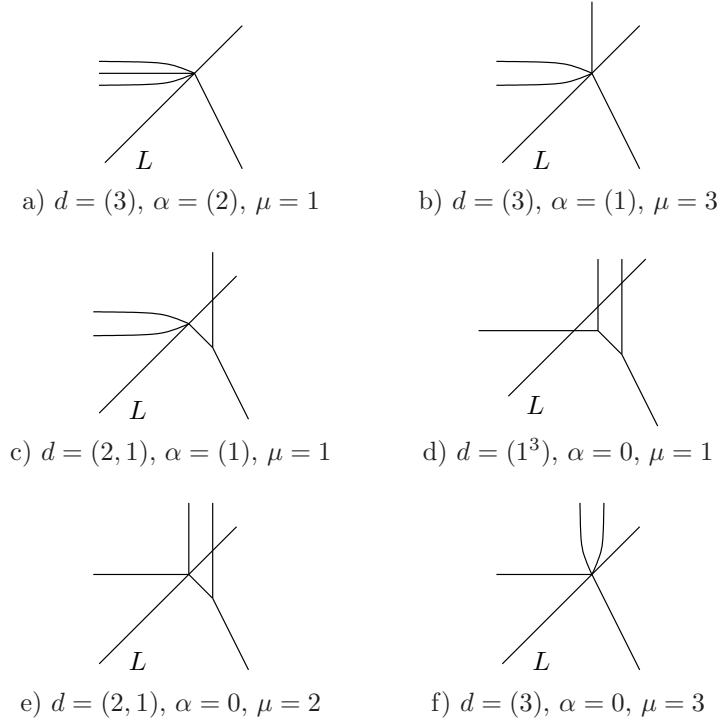

\begin{lemma}\label{lem:compute2}
Suppose that $d_{1}\ge\ldots\ge d_{r}$
and $d_{r+1}\ge\ldots\ge d_{r+s}$. 
Then we have
$$\mathcal N (u, 2, d,\alpha)=\left\{
\begin{array}{l}
3^{t_b} \cdot \binom{3t_d+t_c+t_e}{t_c, t_e,3,\ldots,3}\frac{t_c!t_e!}{t_d!}
\\ \\ \quad   \mbox{if }  d=(3^{t_a+t_b},2^{t_c},3^{t_f},2^{t_e}, 1^{3t_d+t_c+t_e}) 
\mbox{ and }\alpha=(2^{t_a},1^{t_c+t_b})
\\ \quad   \mbox{with } t_a+t_b+t_c+t_f+t_e+t_d=u
\\ \\
\\   0 \quad \text{otherwise,}
\end{array}\right. $$
where there are $t_d$ copies of $3$ in the above multinomial coefficient.
\end{lemma}
\begin{proof}
Since the technique is exactly as in Lemmas \ref{lem:compute0}
and \ref{lem:compute0}, we briefly sketch the proof.
Each component of an element of $\TT\mathcal S'(u,1,d,\alpha,\omega)$
has to be mapped to $\mathbb R^2$ as depicted in Figure
\ref{fig:n=2}, and the result follows from considering all possible
partitions of the points $p_1,\ldots,p_{r+s}$. Here, $t_i$ denotes the number of components as in Figure
\ref{fig:n=2}i for all i=a,b,c,d,e,f. 
\end{proof}

\section{Tropical enumerative geometry in the plane $X$}\label{sec-tropenum}
This section is the core of the present paper.
We start by setting up
a tropical enumerative problem
in the tropical surface $X$  (see Example \ref{const-surface}) that we relate to the enumerative
geometry of complex Hirzebruch surfaces via a Correspondence Theorem
(Section \ref{sec:corr}). This latter is an adaptation of
Mikhalkin's Correspondence Theorem in \cite{Mi03}, and will be proved in
Section \ref{sec-corres}.
In Section
\ref{sec:proof main}, we prove our main result, Theorem \ref{main formula}. 
We also apply our new method to deduce a formula enumerating \emph{irreducible} curves for some cases. This formula generalizes a result by Abramovich and Bertram (see \cite{Vak00b}).

\vspace{2ex}

\subsection{Basic tropical enumerative geometry in $X$}\label{enum general}
Here we describe a particular kind of
tropical enumerative problems 
concerning tropical 
morphisms through point conditions in $X$, and describe 
properties of the
tropical morphisms that are solutions.
Recall that the tropical surface $X$ is made of three
2-dimensional cells, $\sigma_1=\{x=y\geq z\}$, $\sigma_2=\{x=z\geq
y\}$ and $\sigma_3=\{y=z\geq x\}$, 
meeting along the line $L=\RR(1,1,1)$.

\begin{notation}\label{enum pb}
Let $\Delta$ be a Newton fan only containing vectors in $X$, but no
vectors in $L$, and let $\chi\in \ZZ$.
We denote by  $\Delta_i$ the set of directions of $\Delta$ 
in $\sigma_i$ for $i=1,2,3$, and by $d$ 
the intersection multiplicity of $L$ with the Newton fan $\Delta$ 
(see Remark \ref{total int}). 
For the rest of this section, we assume that any direction in $\Delta_3$
has tropical intersection multiplicity 1 with $L$, i.e. $d=\#
\Delta_3$.

Given a 
configuration $\omega$ 
of $\# \Delta - \chi - d$ points in
$\sigma_1\cup \sigma_2$, we denote by $\TT\mathcal S(\omega)$ 
the set of all tropical
  morphisms $h:C\to X$  with Newton fan 
$\Delta$, with
  $\chi_{\trop}(C)=\chi$, and  
passing through all points in $\omega$. 

\end{notation}
In the proof of Theorem \ref{main formula} we
use a configuration
$\omega\subset \sigma_1$, however almost all results from this section
still hold without
this assumption.
The rest of this section is  devoted to prove the following
proposition describing the elements of $\TT\mathcal S(\omega)$.
Recall that given a tropical morphism $h:C\to X$, we define
$C_i=h^{-1}(\sigma_i)$ (see Notation \ref{not-ci}).
\begin{proposition}\label{prop-curves}
For a generic  configuration $\omega$,
the set $\TT\mathcal S(\omega)$ 
is finite. Moreover any tropical morphism $h:C\to
X$ in $\TT\mathcal S(\omega)$ satisfies the following properties:
\begin{enumerate}
\item there is no edge $e$ of $C$ with $h(e)\subset L$ set theoretically;

\item any vertex $V$ of $C$ such that $h(V)\notin L$ is 3-valent; 
furthermore $h$ is
  an embedding in a neighborhood of $V$;

\item the tropical curve $C$ is explicit;

\item for any point $p\in\omega$, the set $h^{-1}(p)$ consists of a
  single point which is in the interior of an edge of $C$; 
\item $C_3$ is a union of $d$ ends of $C$;
\item for any vertex $V$ of $C$ such that $h(V)\in L$ one has
  $ov(V)=0$, and $V$ is adjacent to exactly 
one edge of $C_1$ and $C_2$;

\end{enumerate}

If $\omega \subset \sigma_1$ we have in addition:
\begin{enumerate}

\item[(7)] $C_2$ is a union of $\#\Delta_2$ trees.
\end{enumerate}

\end{proposition}
Note that it follows from the conditions above
that $C$ contains no degenerate edges, and that
$C_3$ has exactly $d_V$ ends adjacent to each vertex $V$ with $h(V)\in L$.
The proof of Proposition \ref{prop-curves} will follow from several
technical lemmas.

We first recall some basic facts about the structure of the space
$M^\alpha$ of all
tropical morphisms to $\RR^2$ or $X$ 
with a given combinatorial type $\alpha$. The space $M^\alpha$ is also
called the \emph{(rigid) space of deformations} of $\alpha$ (or of a tropical morphism $h$ of type $\alpha$).
We can naturally identify the space $M^\alpha$  with an unbounded open
polyhedron in $\RR^{2+\#\Ed^0(C)}$: a tropical morphism $h:C\to \RR^2$
or $X$
 in $M^\alpha$
is entirely determined by the
coordinates of a root vertex of $C$  and the lengths of all bounded
edges of $C$.
 We cannot vary the
lengths 
of non-degenerate edges
independently however, since they have to satisfy the
equations that the cycles close up, and that certain vertices and
edges are mapped to $L$ if $h$ is a morphism to $X$. 
We call the dimension of $M^\alpha$ also the \emph{dimension of $\alpha$}.

We say that a plane tropical morphism $h:C\to X$  (resp. $h:C\to
\RR^2$) \emph{contracts a
  cycle $\gamma$ of $C$} if  the set-theoretic intersection of $X$
with the (classical)
affine span in $\RR^3$ of $h(\gamma)$ (resp. if $h(\gamma)$) 
has dimension at most 1.  

\begin{example}
The tropical morphisms to $\RR^2$ depicted in Figure \ref{fig:ex contracted
  cycle}a and b do not contract the cycle, however the one to $\RR^2$
and  to $X$  respectively depicted in
Figure \ref{fig:ex contracted cycle}c and d do.
\begin{figure}[ht]
\centering
\begin{tabular}{ccccccc}
\input{Cycle3.pstex_t}& \hspace{2ex} &
\input{Cycle4.pstex_t}  & \hspace{2ex} &
\input{Cycle2.pstex_t}  & \hspace{2ex} &
\input{Cycle1.pstex_t}  
\\  a) && b) && c) && d)
\end{tabular}
\caption{Examples of non-contracted and contracted cycles.}
\label{fig:ex contracted cycle}
\end{figure}
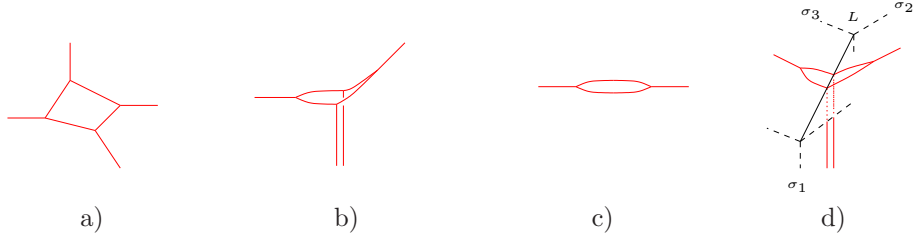
\end{example}

\begin{lemma}[see e.g. {\cite[Proposition 3.9]{GM051}}]\label{rem-dimplane}
The dimension of a combinatorial type $\alpha$ of tropical
morphisms to $\RR^2$ with Newton fan $\delta$ and topological Euler characteristic $\chi$ 
which does not contract any cycle is at most
$\#\delta-\chi$.
Moreover 
 $\dim(M^\alpha)=\#\delta-\chi$
 if and only if every vertex of $\alpha$ is
$3$-valent 
and $\alpha$ does not contain any degenerate edge.
\end{lemma}

The dimension of a combinatorial 
type of tropical morphisms to $X$ is much harder to determine.
 We compute it only in a special situation which will be sufficient
 for our purposes. 
Recall that given  a tropical 
morphism $h:C\to X$ , an edge $e$ of $C$ is 
 tropically mapped to $L$ if $h(e)\subset L$ and $v(V,e)$ is
parallel to $(1,1,1)$, and a vertex $V$ of $C$ is 
tropically mapped to $L$ if either $d_V>0$ or each edge  adjacent to
$V$ is tropically mapped to $L$.

\begin{notation}\label{not1}
Let $h:C\to X$ be a tropical morphism 
of type $\alpha$ with no vertex $V$ tropically mapped to $L$
 with $d_V=0$.
 We define $$\widetilde C=C\setminus\{e\in\Ed(C)\; |\;
h(e)\text{ is tropically mapped to } L\}$$ and $\widetilde
h=h_{|\widetilde C}$.
Note that $\chi(\widetilde C)=\chi(C) +l$, where $l$ is the number
of bounded
edges of $C$ tropically mapped to $L$, and that the rigid space of
deformations of $\widetilde h$ is naturally identified with
$M^\alpha$.
We define $\widetilde C_i=\widetilde h^{-1}(\sigma_i)$ 
and $\widetilde h_i=\widetilde h_{|\widetilde C_i}$. As explained
in
Notation \ref{not-ci}, we may 
think of  $\widetilde h_i:\widetilde C_i\to
\sigma_i\subset \RR^2$ as a
tropical morphism to the plane $\RR^2$. 
Let   $V_1,\ldots,V_k$ be the vertices of $\widetilde C$
 tropically mapped to $L$, and
 $y_{ji}$ for $j=1,\ldots,k$ and $i=1,2,3$ the number of edges of
 $\widetilde C_i$ 
adjacent to $V_j$. We set $y_i=\sum_{j=1}^k y_{ji}$.
We denote by
 $\chi_i$ the topological Euler characteristic of $\widetilde C_i$,
and by $v_j$ the valency of $V_j$ in $\widetilde C$.
By definition we have
$$v_j=\sum_{i=1}^3 y_{ji}, \quad \quad \sum_{j=1}^k v_j=y_1+y_2+y_3$$
and 
\begin{equation}\label{eq-euler}
\chi(\widetilde C)= \chi_1+\chi_2+\chi_3-\sum_{j=1}^k (v_j-1)=
\chi_1+\chi_2+\chi_3-y_1-y_2-y_3+k.
\end{equation}

\end{notation}

\begin{lemma}\label{lem-dim}
Consider a combinatorial type $\alpha$ of 
tropical morphisms $h:C\to
X$  without any
contracted cycle, with no vertex $V$ tropically mapped to $L$ with $d_V=0$,
and which maps $l$ bounded edges of $C$ tropically to $L$.
Then the dimension of $\alpha$ is less than or equal to
$\#\Delta-\chi(C)-l -d$.
\end{lemma}
\begin{proof}
According to Lemma \ref{rem-dimplane}, the space of deformations
of each morphism
$\widetilde h_i$ has  dimension  at most $\#\Delta_i+y_i-\chi_i$.
We cannot vary $\widetilde h_1$, $\widetilde h_2$ and $\widetilde h_3$ 
independently however, since we want to glue them to a single map
 $\widetilde h$.
First of all, this imposes some conditions on $\widetilde h_1$ itself: for each
vertex $V_j$, we have to require that the $y_{j1}$ adjacent edges meet
the same point on $L$. This yields $y_{j1}-1$ conditions for each
vertex $V_j$. 
These conditions are all independent since
we assumed that $C$ has no contracted cycle, hence
 we get $y_1-k$ conditions altogether for $\widetilde h_1$. 
Since all $y_{j2}$ edges of $\widetilde C_2$
adjacent to $V_j$ meet the point  $\widetilde h_1(V_j)$, 
we get $y_2$ conditions for $\widetilde h_2$ altogether. 
Also those are independent because we do not have any
contracted cycle. Analogously, we get $y_3$ independent conditions for
$\widetilde h_3$. 
Thus the dimension of $\alpha$ is less than or equal to 
\begin{align*}
  &\#\Delta_1+y_1-\chi_1+\#\Delta_2+y_2-\chi_2+\#\Delta_3+y_3-\chi_3-y_1-y_2-y_3+k\\=&
  \#\Delta-\chi(\widetilde C)-y_1-y_2-y_3+2k\\=& \#\Delta
  -\chi(C)-l -\sum_{j=1}^k v_j+2k\\ =& \#\Delta -\chi(C) -l
  -\sum_{j=1}^k (\ov_{V_j}+d_{V_j}+2)+2k\\ \leq &\#\Delta-\chi(C) -l
  -\sum_{j=1}^k (d_{V_j}+2)+2k
\\=& \#\Delta -\chi(C) -l-d-2k+2k\\=& \#\Delta-\chi(C) -l-d\end{align*}
where the first equality follows from Equation (\ref{eq-euler}),
and the inequality holds since $\ov_{V_j}\geq 0$ for $j=1,\ldots,k$. 
\end{proof}

\begin{lemma}\label{lem-ineq} 
Consider a combinatorial type $\alpha$ of tropical morphisms $h:C\to
X$ as in Notation \ref{enum pb} and \ref{not1}.
For each $i=1,2,3$ we have
\begin{displaymath} 
\#\Delta_i\geq \chi_i \mbox{ and } y_i\geq \chi_i.
\end{displaymath}
The first inequality is an equality if and only if 
$\widetilde C_i$ is a union of $\#\Delta_i$ trees, the second if and
only if $\widetilde C_i$ is a union of $y_i$ trees.  
\end{lemma}
\begin{proof}
First note that there cannot be components of $C$ whose image lies in
$\sigma_i$. If that was the case, it follows from the balancing
condition that such a component can only have ends with directions in
$L$ which contradicts our general assumptions on $\Delta$.
 Indeed, if such a
component has an end of another direction in $\sigma_i$, then it must
meet the line $L$ 
and hence it contains also parts in the
other faces of $X$. 
Let $k_i$ denote the number of connected components of $\widetilde C_i$.
 Since there are no components whose image lies in $\sigma_i$,
 we have $k_i\leq y_i$. Since every
connected component must contain at least
 one end, we have $\# \Delta_i\geq k_i$.
Any connected component of $\widetilde C_i$ has 
Euler characteristic less than or equal to 1, so 
$\chi_i\le k_i$.
Moreover we have the equality
$\chi_i=k_i$ if and only if 
$\widetilde C_i$ is a union of $k_i$
trees. 
\end{proof}

\begin{lemma}\label{lem-explicit}
 Let $h:C\to X$ be an element of $\TT\mathcal S(\omega)$ 
with no vertex $V$ tropically mapped to $L$ with $d_V=0$,
 and which
does not   contract any cycle. 
Then $h$ satisfies properties 
$(3)$ and $(4)$ of Proposition \ref{prop-curves}.
Furthermore one has the following:
\begin{itemize}
\item no bounded edge of $C$ is tropically mapped to $L$;
\item any vertex  of $C$ tropically mapped to
 $ L$ 
 is adjacent to exactly 
one edge of $C_1$ and $C_2$;

\item  every vertex $V$ of $C_1$ or $C_2$ not tropically mapped to $L$ is
$3$-valent and $h$ is
  an embedding in a neighborhood of $V$;

\item any degenerate edge of $C$ contained in $C_1$ or $C_2$ is
  adjacent to a vertex of $C$ tropically mapped to $L$.
\end{itemize}
\end{lemma}
\begin{proof}
We use the notations introduced in \ref{enum pb} and \ref{not1}.
 According to Lemma
\ref{rem-dimplane} 
the dimension of the space of deformations of $\tilde{h_i}$ is at most 
 $\#\Delta_i+y_i-\chi_i$. 
Assume that $r$ of the $\#\Delta-\chi-d$ point conditions lie in $\sigma_1$, the other $\#\Delta-\chi-d-r$ points in $\sigma_2$.
Then the morphism $\widetilde h_1:\widetilde C_1\to \RR^2$ passes
through a generic configuration of $r$ 
 points in the plane, and satisfies moreover the gluing
 conditions that the $y_{1j}$ ends adjacent to $V_j$ in $\widetilde
 C_1$
 must meet
 the same point. As in the proof of \ref{lem-dim} this gives $b_1=y_1-k\geq 0$ extra
 independent conditions 
 on $\widetilde h_1$.
Hence we can conclude that $\#\Delta_1+y_1-\chi_1\geq r+b_1$.
In the same way we have
$\#\Delta_2+y_2-\chi_2\geq \#\Delta-\chi-d-r +b_2 $ where $b_2=y_2-k\geq 0$.
Finally, we have $k$ extra independent gluing conditions for the
morphisms $h_1$ and $h_2$ to match along $L$. 
Altogether we have
\begin{align*}
 &0 \geq r+b_1-\#\Delta_1-y_1+\chi_1+ \#\Delta-\chi-d-r+b_2-\#\Delta_2-y_2+\chi_2+k \\ =& 
\#\Delta_3+ \sum_V g_V -\chi(\widetilde C) +l  -d +k-y_1+\chi_1 +b_1 -y_2+\chi_2+b_2 \\ =&
\#\Delta_3+ \sum_V g_V - \chi_3 +y_3-k +l -d +k +b_1+b_2 \\=&
\sum_V g_V+ (y_3-\chi_3)+l+b_1+b_2
\end{align*}
where the second equality follows from Equation (\ref{eq-euler}) and
the third equality follows since $d=\#\Delta_3$ by the assumption
made in Notation \ref{enum pb}.
The five summands above are all non-negative. For the second
one, this follows from Lemma \ref{lem-ineq}, for the others it is obvious.
We can thus conclude that each summand is zero. It follows that $l=0$,
the curve $C$ is explicit, and that $y_1=y_2=k$. Hence the tropical morphism $h$ satisfies property $(3)$.
Using Lemma \ref{lem-ineq}, it follows that  $C_3$ is a union of $y_3$ trees.

Also, it follows that the dimension of the type of $ h_i: C_i\to
\RR^2$ equals $\#\Delta_i+y_i-\chi_i$ for $i=1,2$. By Lemma
\ref{rem-dimplane}, every vertex of $C_1$ and $C_2$ 
not tropically mapped to $L$
is $3$-valent, and every bounded edge of $C_1$ or $C_2$ is non-degenerate. 

Suppose now that $C$ has a $3$-valent vertex $V_0$ not tropically
mapped to $L$
in the neighborhood of which $h$ is not an
embedding. This means that two edges of $C$ have the same (primitive) direction
vector from $V_0$. Since $h_i$ does not contract any cycle and
$y_1=y_2=k$, 
by gluing
edges which have the same image by $h_i$ one can produce a tropical
morphism $h':C'\to X$ in $\TT\mathcal S(\omega)$ 
with no vertex $V$ tropically mapped to $L$ with $d_V=0$,
 which
does not   contract any cycle, and with a vertex $V_1$  not tropically
mapped to $L$
which is not
$3$-valent. But we just showed above that
this is impossible, hence a contradiction.

We showed that the rigid space of deformation of $h$ has dimension $\#\Delta
-\chi -d$, which is exactly the number of independent conditions
imposed by
$\omega$. In particular $h$ cannot satisfy any further independent
condition,
like having a vertex  or two points of
$C$ mapped to a point of $\omega$, i.e. $h$
also satisfies $(4)$ 
\end{proof}

\begin{remark}\label{rem-c2}
 If we assume that our point conditions lie in $\sigma_1$, a modification of the proof of Lemma \ref{lem-explicit} above shows in addition that $C_2$ is a union of $\#\Delta_2$ trees: In this case, $\#\Delta_1+y_1-\chi_1\geq \#\Delta-\chi-d+b_1$ and 
\begin{align*}
 &0 \geq
 \#\Delta-\chi-d+b_1-\#\Delta_1-y_1+\chi_1 \\ =& 
 (\#\Delta_2-\chi_2) + \sum_V g_V + (y_2-k) + (y_3-\chi_3) + b_1 +l,
\end{align*}
Here, the first summand is nonnegative due to \ref{lem-ineq} and the third summand is nonnegative since there must be an edge inside $\sigma_2$ adjacent to each of the $k$ vertices in $L$. As before it follows that all summands are zero, and in addition to the results of Lemma \ref{lem-explicit}, we can conclude that in this situation $C_2$ is a union of $\#\Delta_2$ trees.
\end{remark}

\begin{corollary}\label{cor-nocyc}
 If $h:C\to X$ is an element of $\TT\mathcal S(\omega)$, 
then  $h$ does not contract any cycle  and 
no bounded edge of $C$ is tropically 
mapped to $L$. In particular, $C$ does 
not contain any vertex $V$ tropically mapped to $L$ with $d_V=0$.
\end{corollary}
\begin{proof}
Suppose that $h:C\to X$ is an element of $\TT\mathcal S(\omega)$
which  contains 
a vertex $V$ tropically mapped to $L$ with $d_V=0$. 
We denote by $C_L$ the union of vertices and 
edges of $C$ which are tropically mapped
to $L$.

On each connected component of $h^{-1}(L)\cap C_L$ we identify 
 all points having the same image in $L$. In this way, we produce  
 a tropical morphism $h':C'\to X$ with the same image
 as $h$,  with no vertex $V$ tropically mapped to $L$ with $d_V=0$, and with 
$b_1(C')\le b_1(C)$. Note
 that any vertex of $V$ of $C'$ obviously satisfies $\ov_V\ge 0$. 
If $b_1(C')< b_1(C)$ we increase the genus of an arbitrary vertex of $C'$
in such a way that $g(C')=g(C)$.
If $h'$ does not contract any cycle, then it follows from
 Lemma \ref{lem-explicit} that no bounded edge of $C'$ is tropically mapped
  to $L$. So $C=C'$, $h=h'$, and
 the corollary is proved in this case.

So we are left to prove the corollary when $h$ contracts a cycle and 
$C$ does not contain
any vertex $V$ tropically mapped to $L$ with $d_V=0$.
Let $\gamma$ be a contracted cycle of $h:C \to X$. Because of the
above, we may assume that $\gamma\cap C_L$ is finite and
$h(\gamma\cap C_L)\cap L$ is either empty or a single point. There is a  continuous involution on 
$\gamma\setminus C_L$ 
which exists on any morphism with the same combinatorial type as $h$.
We can quotient
 $\gamma$ by this involution producing a tropical premorphism
$h':C'\to X$, where
 $b_1(C')=b_1(C)-1$. We want to construct a tropical morphism
 $h'':C''\to X$ of genus $g(C)$ out of $h':C'\to X$. To do so, we have
 to increase the genus of one vertex of $C'$
by one, making
 sure that the overvalencies at vertices tropically 
mapped to $L$ are all nonnegative.

If 
$h(\gamma\cap C_L)\cap L=\emptyset$
then the overvalency of the
vertices of $C'$ mapped to $L$ are  the same as those of 
$C$, and we increase by one 
the genus of an arbitrary
vertex of $C'$.

If $h(\gamma\cap C_L)\cap L\neq \emptyset$ then
 it is a
point. Let $V_1,\ldots,V_s$ be 
the vertices of $\gamma$ tropically mapped to $L$.
The involution identifies the vertices $V_1,\ldots,V_s$ of $C$ with a single
vertex $V$ of $C'$,
with overvalency
$\ov_V=\sum_{i=1}^s \ov_{V_i}-s'+2(s-1)$, where $s'$ is the number of
pairs of edges in a neighborhood of the $V_i$ 
that are glued together by the involution. Since $s'=s$, 
 the overvalency can only become negative if $s=1$. This is the case
when the involution identifies two edges adjacent to a single vertex $V_1$
on $L$. In this situation we increase the genus of $V_1$ by one (the
corrected overvalency is now positive) thus again producing a morphism
$h'':C''\to X$ of the same genus. Otherwise, we increase the genus of
an arbitrary vertex of $C'$ by one.

 We repeat this for any contracted cycle. Eventually, we obtain a
 tropical morphism $h'':C''\to X$  in $\TT\mathcal{S}(\omega)$ of genus
 $g$, without any contracted cycle, and
 with $b_1(C'')<b_1(C)$. It 
 follows from Lemma \ref{lem-explicit} that $C''$ is explicit, a contradiction.
 \end{proof}

\begin{lemma}\label{lem-ov0}
If $h:C\to X$ is an element of $\TT\mathcal S(\omega)$, 
then $\ov(V)=0$ for any vertex $V$ tropically mapped to $L$,  and
$C_3$ is a union of $d$ ends. 
\end{lemma}
\begin{proof}
By Corollary
\ref{cor-nocyc} the morphism $h$ does not contract any cycle and does
not tropically map
any bounded edge of of $C$  to $L$.
 By Lemma  \ref{lem-explicit}, the curve $C$ is
explicit, and $y_1=y_2=k$. Hence considering the sum $0\leq \sum_V\ov(V)$
over all vertices with $h(V)\in L$ we obtain
$$
\begin{array}{ccl}
 \sum_V\ov(V)  &=& \sum_V (\val(V)- d_V-2+2g_V) 
\\ &=&
-d+ \sum_V \val(V) - 2k
\\ &=& -d+y_1+y_2+y_3-2k
\\ &=& -d+y_3
\\ &\leq 0
\end{array}
$$
since by definition of $d=C\circ L$ we have  $d\geq y_3$. Hence $ \sum_V\ov(V)
=0$ and $y_3=d$ as claimed.
\end{proof}

\begin{proof}[Proof of Proposition \ref{prop-curves}]
There are
finitely many combinatorial types of tropical morphisms in $X$ with
Newton fan  $ \Delta$ and of Euler characteristic $\chi$. 
If $h:C\to X$ is an element of $\TT\mathcal S(\omega)$ of combinatorial
type $\alpha$, 
it follows from Lemmas \ref{lem-dim}, \ref{lem-explicit} and Corollary \ref{cor-nocyc} that $\dim(M^\alpha)\leq \#\Delta-\chi-d$.
Since we fix $\#\Delta-\chi-d$ independent 
linear conditions, it follows that we have equality and that there is
a unique tropical morphism of combinatorial type $\alpha$ in
$\TT\mathcal S(\omega)$. This proves that  $\TT\mathcal S(\omega)$ is
finite.

 Corollary \ref{cor-nocyc} implies that the assumption of 
 Lemma \ref{lem-explicit} are satisfied by all elements of
 $\TT\mathcal S(\omega)$, in particular  $(3)$ and $(4)$ are satisfied.
Suppose that $C$ has a degenerate edge $e$. By  Lemma
\ref{lem-explicit}, any degenerate component of $C$ is a tree and any edge $e$ of this tree
is adjacent to a vertex $V$ of $C$ tropically mapped to $L$. The other
vertex $V'$ adjacent to $e$ is not tropically mapped to $L$ since
otherwise $e$ would be tropically mapped to $L$ contradicting
Corollary \ref{cor-nocyc}. Hence $V'$ is a trivalent vertex of
$C$. Now two possibilities can occur:
\begin{enumerate}
\item[(a)] $V'$ is adjacent to only one degenerate edge of $C$; in
  this case, by contracting $e$ we produce a new tropical morphism
  $h':C'\to X$ in $\TT\mathcal S(\omega)$ with a vertex $V$
  satisfying $ov(V)>0$.

\item[(b)] $V'$ is adjacent to two degenerate edges; in this case by
  gluing these two edges as in the proof of Corollary \ref{cor-nocyc}
we produce a new tropical morphism
  $h':C'\to X$ in $\TT\mathcal S(\omega)$ with a vertex $V$
  satisfying $ov(V)>0$.
\end{enumerate}
Hence both cases contradict  Lemma \ref{lem-ov0}, and $C$ does not
have any degenerate edge. In particular any element of $\TT\mathcal
S(\omega)$ also satisfies $(1)$.

Lemmas \ref{lem-ov0}  and  \ref{lem-explicit} 
give  $(5)$ and $(6)$. 
It follows from $(5)$ that there are no vertices of $C$ mapped to
in $\sigma_3\setminus L$. Since we know from Lemma \ref{lem-explicit} that $(2)$ is satisfied for any vertex in $\sigma_1\cup \sigma_2\setminus L$, $(2)$ follows.

Remark \ref{rem-c2} proves $(7)$ for the case when $\omega\subset \sigma_1$.
\end{proof}

\subsection{Relation with enumerative geometry of Hirzebruch surfaces}\label{sec:corr}

One can compute tropically enumerative invariants of Hirzebruch
surfaces by enumerating tropical curves in $\RR^2$ with Newton fan 
(see Example \ref{ex faninsigman}):
$$\delta_0=\{(1,n)^a, (0,-1)^{an+b}, (-1,0)^a, (0,1)^b \}. $$ 
However to obtain the equation of Theorem \ref{main formula}, one has to 
enumerate tropical curves in $X$: Theorem \ref{main formula} is
based on the deformation of $\Sigma_{n+2}$ to $\Sigma_n$ which can be
described tropically by the tropical surface $X$ (see Appendix \ref{hirz}
for more detailed explanations). Hence we are interested in tropical
curves with Newton fan
$$\{(1,n,1)^a, (0,-1,1)^{an+b}, (-1,0,0)^a,
(0,1,0)^b,(0,0,-1)^{a(n+1)+b} \}. $$

For convenience later in the formula,
we apply the
transformation $(x,y,z)\mapsto (x,-y,z)$, i.e. we fix the following Newton
fan
$$\Delta= \{(1,-n,1)^a, (0,1,1)^{an+b}, (-1,0,0)^a, (0,-1,0)^b,
(0,0,-1)^{a(n+1)+b} \}.$$
Using Notation \ref{enum pb}, we have 
\begin{align*}\Delta_1&=\{(0,0,-1)^{a(n+1)+b}\},\\ \Delta_2&=\{(1,-n,1)^a,(0,-1,0)^b\}\mbox{ and}\\ \Delta_3&=\{ (0,1,1)^{an+b}, (-1,0,0)^a\}.\end{align*}
Note that here 
$d=(n+1)a+b$, and
$\#\Delta_3=d$. In particular we are in the situation covered by
Section \ref{enum general}.
As in Section \ref{enum general}, let us choose an integer
$\chi\in\ZZ$, and a generic configuration $\omega$ of 
$\#\Delta-\chi-d$ points in $\sigma_1\cup \sigma_2$. 

Following Notation \ref{enum pb}, we denote by $\TT\mathcal S(\omega)$
the set of all tropical morphisms $h:C\to X$ passing through $\omega$,
with Newton fan $\Delta$, and 
with $C$ a tropical curve such that $\chi_{\trop}(C)=\chi$.
Recall that any element of $\TT\mathcal S(\omega)$ satisfies the
properties $(1)-(6)$ given in Proposition \ref{prop-curves}.

Given an element $h:C\to X$  of $\TT\mathcal S(\omega)$, we denote by
 $\Ve_L(C)$ (resp. $\Ve_{\sigma_i}(C)$) the set of vertices of $C$
 mapped to $L$ (resp. $\sigma_i\setminus L$).
If  $V\in\Ve_{\sigma_1}(C)\cup\Ve_{\sigma_2}(C)$, then it follows from Proposition \ref{prop-curves} that $\val(V)=3$. 

Given $V\in\Ve_{\sigma_i}(C)$, we choose any two of its adjacent edges
$e_{V,1}$ and $e_{V,2}$.
Note that 
we have 
$v(e_{V,j})=(a_{V,j},a_{V,j}, b_{V,j})$ if $i=1$, and
 $v(e_{V,j})=(a_{V,j},b_{V,j},a_{V,j})$ if $i=2$ for some $a_{V,j}$ and $ b_{V,j}$.

A vertex $V\in\Ve_{L}(C)$ is  adjacent to $d_V$ ends  mapped to $\sigma_3$,
 say $k_V$ ends with direction
$(-1,0,0)$ and $l_V$ ends with direction $(0,1,1)$ (pointing away from
 $L$). Note that $k_V+l_V=d_V$.

\begin{definition}[Multiplicity of a tropical morphism in $\TT\mathcal S(\omega)$]\label{def:mult}
Let $h:C\to X$ be an element of $\TT\mathcal S(\omega)$.

We define the multiplicity of a vertex $V\in\Ve_{\sigma_i}(C)$ as

$$\mu_V= \left|\det\left( \begin{array}{cc}a_{V,1}&a_{V,2}\\ b_{V,1}& b_{V,2} \end{array}  \right)\right|. $$

We define the multiplicity of a vertex $V\in\Ve_{L}(C)$ as
$$\mu_V=\binom{k_V+l_V}{k_V}. $$

We define the multiplicity of $h$ as
$$\mu_h=\prod_{V\in\Ve(C)}\mu_V.$$
\end{definition}

We also define the two following numbers
\begin{displaymath}
 \TT N_{\chi}(\Delta,\omega)=\sum_{h\in\TT\mathcal S(\omega)}\mu_h.
\end{displaymath}
and 
\begin{displaymath}
 \TT N^{\irr}_{\chi}(\Delta,\omega)=\sum_{h\in\TT\widetilde{\mathcal S}(\omega)}\mu_h.
\end{displaymath}
where $\TT\widetilde{\mathcal S}(\omega)$ is the set of 
tropical morphisms $h:C\to X$ in $\TT\mathcal S(\omega)$ from an irreducible
 tropical curve $C$. 

The next Theorem is one of  the main results of this paper.
\begin{theorem}[Correspondence Theorem]\label{thm-corres}
Let $\Delta,\chi$, and $\omega$ be as explained in the beginning of
Section \ref{sec:corr}, and let as in Notation \ref{deltazero} 
$$\delta_0=\{(1,n)^a, (0,-1)^{an+b}, (-1,0)^a, (0,1)^b\}.$$ 
Then we have
$$\TT N_\chi(\Delta,\omega)=N_{2\chi}(\delta_0)\quad
\text{and}\quad  \TT N^{irr}_{\chi}(\Delta,\omega)=N^{irr}_{2\chi}(\delta_0).$$
\end{theorem}
We 
prove this theorem  in Section \ref{sec-corres} in the case of
irreducible curves, from which the case of reducible curves follows immediately.
A consequence of the Correspondence Theorem \ref{thm-corres} is that the
numbers  $\TT N_\chi(\Delta,\omega)$ and 
$\TT N^{irr}_\chi(\Delta,\omega)$
do not depend on the choice of $\omega$, as long as $\omega\subset
\sigma_1\cup\sigma_2$ is generic. 
We will thus also write
$\TT N_\chi(\Delta)$ and 
$\TT N^{irr}_\chi(\Delta)$ 
for $\omega$
satisfying the requirements.

\subsection{Proof of Theorem \ref{main formula}}\label{sec:proof main}

We still fix $\Delta$, $\chi$, and $\omega$ 
 as in the beginning of
Section \ref{sec:corr}.
\begin{notation}\label{not-pointslow}
For the rest of the paper, we suppose that 
  $\omega\subset \sigma_1$ and that the points in $\omega$
have very low $z$-coordinate compared to the $x$ and
  $y$-coordinates.
\end{notation}

To prove Theorem \ref{main formula}, we need the following lemma,
which describes how tropical morphisms in $\TT\mathcal S(\omega)$ can
meet the line $L$. We still use Notation \ref{not-ci}. Remember that any element of $\TT\mathcal S(\omega)$ satisfies the
properties $(1)-(7)$ of Proposition \ref{prop-curves}.

\begin{figure}[ht]
\centering
\begin{tabular}{ccc}
\input{vertextypeskleinA.pstex_t}
& \hspace{6ex} &
\input{vertextypeskleinB.pstex_t} 
\\ a) && b)
\\ \\
\input{vertextypeskleinC.pstex_t}
& \hspace{6ex} &
\input{vertextypeskleinD.pstex_t} 
\\ c) &&  d) $\begin{array}{l} n\geq \alpha\geq 0 \\ 0<k+l\leq n+1
\\ \beta=k+l-\alpha>0\end{array}$
\end{tabular}
\caption{Four ways to hit $L$.}
\label{fig:4types}
\end{figure}
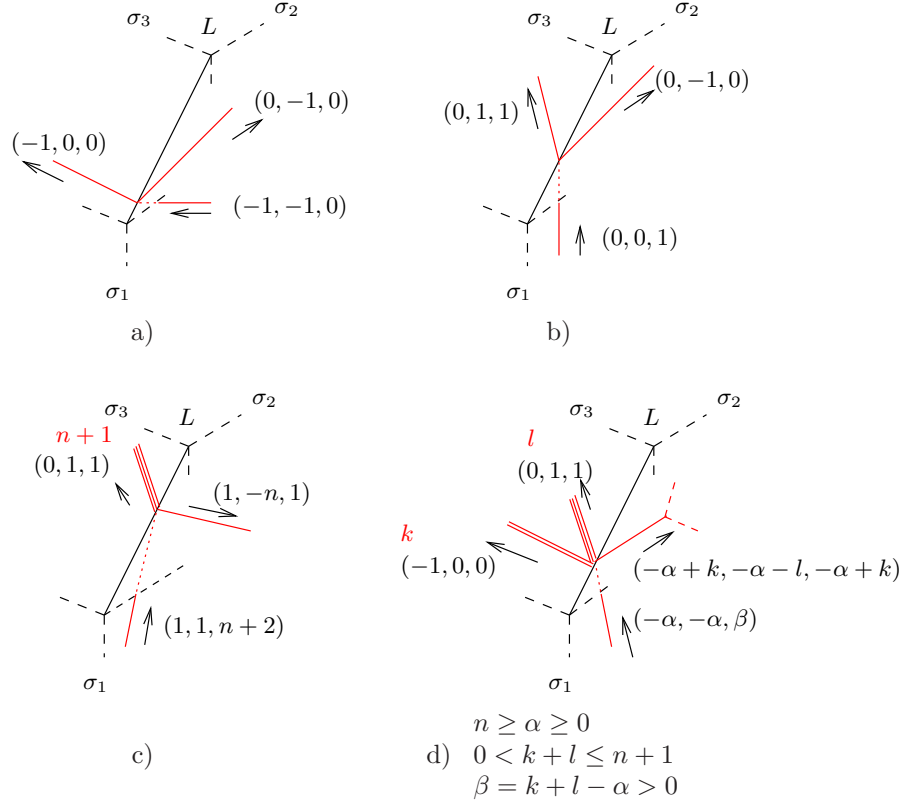

\begin{lemma}\label{lem-vertextypes}
Let $h:C\to X$ be an element of $\TT\mathcal S(\omega)$, and
 $V$ be a vertex of $C$ mapped to the line $L$. Then 
all
possibilities of how $h$ can look like in a neighborhood of $V$ are
depicted in Figure \ref{fig:4types}.
\end{lemma}

Recall that  $ C_3$ only consists of  ends of $C$ (Proposition \ref{prop-curves}(5)), and that any
connected component of $C_2$ contains exactly one end of $C$ by \ref{prop-curves}(7). 
In cases a, b and c, the connected component of 
$ C_2$ containing $V$ also just consists of one end. 
In case d, this component might contain other vertices than $V$,
i.e. the part in $\sigma_2$ in picture $d$ can
continue and contain more vertices on $L$. In particular,
$(-\alpha+k,-\alpha-l,-\alpha+k)$ does not need to be the direction of an
end of $C$.

\begin{proof}
By property $(6)$ in Proposition \ref{prop-curves},
every vertex $V_i$ of $C$ 
mapped into $L$ has exactly one adjacent edge mapped in $\sigma_1$. We
denote by $(-\alpha_i,-\alpha_i,\beta_i)$ the direction of this edge 
(pointing towards $L$), with $-\alpha_i< \beta_i$. Then by \ref{prop-curves}(5) some ends of
$C$  adjacent to $V_i$ are mapped to $\sigma_3$,
 say $k_i$ ends with direction
$(-1,0,0)$ and $l_i$ ends with direction $(0,1,1)$ (pointing away from
 $L$). 
Finally, by \ref{prop-curves}(6) again, exactly one edge adjacent to $V_i$ is mapped to
$\sigma_2$. By the balancing condition
this edge has direction $(-\alpha_i+k_i, -\alpha_i-l_i, \beta_i-l_i)$
(pointing away from $L$). Since  
$k_i+l_i=d_{V_i}=\beta_i+\alpha_i$, we get 
$-\alpha_i+k_i =\beta_i-l_i$.

Let us consider one connected component of $ C_2$, 
call it $\tilde{C}$ and 
assume it meets $L$ at the vertices $V_1,\ldots,V_r$.
Note that the sum of the intersection multiplicities of $\tilde{C}$ at
the $V_i$ equals the intersection multiplicity of the end of
$\tilde{C}$ with $L$. 
If the end of $\tilde{C}$ is of direction $(0,-1,0)$, then it
intersects $L$ with multiplicity $1$, if it is of direction
$(1,-n,1)$, then it intersects with multiplicity $n+1$. So we have
$\sum_{i=1}^r(k_i+l_i)$ equals $1$ in the first case, or $n+1$ in the
second case. 

First, let us consider the case when the end of $\tilde{C}$ is of
direction $(0,-1,0)$.  
Since all the $k_i$ and $l_i$ are nonnegative numbers, and $(k_i+l_i)>0$ for each $i$, it follows that $r=1$ and $\tilde{C}$ is in fact just an end, so $(-\alpha_1+k_1,-\alpha_1-l_1,-\alpha_1+k_1)=(0,-1,0)$.
There are two possibilities how the vertex $V_1$ can look like:
Either $k_1=1$, then $l_1=0$, $\alpha_1=1$ and $\beta_1=0$; 
thus $h$ is as depicted in Figure \ref{fig:4types}a in a neighborhood of $V$.
Or $k_1=0$, then $l_1=1$, $\alpha_1=0$, and $\beta_1=1$; 
thus $h$ is as depicted in Figure \ref{fig:4types}b in a neighborhood of $V$.

Now consider the case when the end of $\tilde{C}$ is of direction $(1,-n,1)$. By the balancing condition, we then have 
 $\big(\sum (-\alpha_i+k_i), \sum (-\alpha_i-l_i), \sum
(-\alpha_i+k_i)\big) = (1,-n,1)$.  
Since we choose the configuration $\omega\subset\sigma_1$ 
very far down from $L$ (Notation \ref{not-pointslow}), we may
 assume that the ``slopes relative to $L$'' of the edges of
 $ C_1$ meeting
 $L$, i.e.\ the values $\frac{-\alpha_i+\beta_i}{-\beta_i-\alpha_i}$,
 decrease from left to right. If we assume the vertices are ordered from
 left to right, i.e.\ $V_1$ is the vertex most left on $L$, $V_r$ the
 most right, then we have
 $\frac{-\alpha_i+\beta_i}{-\beta_i-\alpha_i}\geq
 \frac{-\alpha_j+\beta_j}{-\beta_j-\alpha_j}$ for $i<j$. 
If the relative slopes are not ordered like this for two edges, then
the infinite continuations of these two edges intersect, and change
their order after intersecting. Since we assume that the points are so
far down, we can assume that all these changes of orders happen before
the edges meet $L$, and thus the slopes are ordered as above. 
An example is depicted in Figure \ref{fig:ordering}, where the slopes
relative to $L$ from the most left edge meeting $L$ to the most right are
$1,1,0,-1,-1,-3,-3$. 

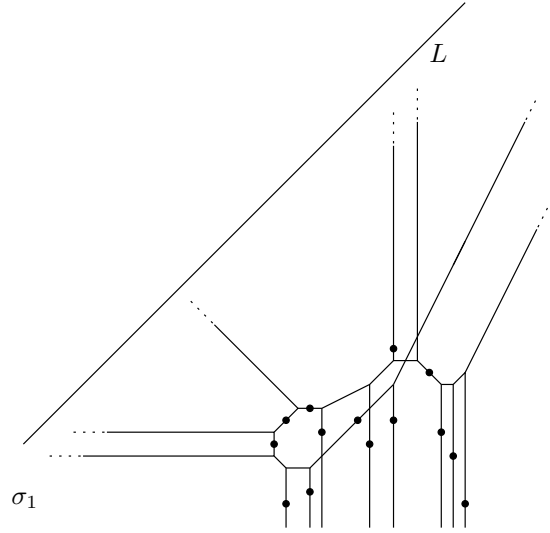
\begin{figure}[ht]
\centering
\begin{tabular}{c}
\input{slopes.pstex_t} 
\end{tabular}
\caption{The slopes
relative to $L$ from the left to the right are
$1,1,0,-1,-1,-3,-3$.}
\label{fig:ordering}
\end{figure}

Note that since $\beta_i+\alpha_i>0$, we have 
$$\frac{-\alpha_i+\beta_i}{-\beta_i-\alpha_i}<1 \Longleftrightarrow
\beta_i>0, \quad \text{and}\quad  
\frac{-\alpha_i+\beta_i}{-\beta_i-\alpha_i}<-1 \Longleftrightarrow
\alpha_i<0.$$
In particular, if $\alpha_i<0$ (resp. $\beta_i>0$)
for some $i$ then also $\alpha_j<0$ (resp. $\beta_j>0$) for all $j>i$. 

Suppose that $\alpha_r<0$.
 Denote by $T$ the smallest subtree of $\tilde{C}$ containing the end
 of $\tilde C$
 and all the vertices $V_j$, $j\geq i_0$, where $i_0$ is the smallest
 index such that $\alpha_{i_0}<0$. Orient the edges in $\tilde{C}$ such
 that they point towards the end of $\tilde C$ (see
Figure \ref{fig:orientation} where the
subtree $T$ is marked with thick edges).   
If $T=\tilde{C}$, then $i_0=1$ and all the $-\alpha_j+k_j$, $j=1\ldots
r$,
 are positive.
If $T\neq \tilde{C}$, there are two-valent vertices of $T$, i.e. in
$\tilde{C}$ two of the adjacent edges belong to $T$ and the third does
not. Since the end of $\tilde C$ belongs to $T$, the third edge then
has to point 
toward the vertex.  
Thus, it needs to connect to some vertices $V_j$, $j<i_0$, behind, and
since these vertices lie to the left of the $V_j$, $j\geq i$, it needs
to have a positive $x$-coordinate (in the orientation as
before). Hence any edge in $T$ and any edge
adjacent to $T$ has a positive $x$-coordinate.  
Every edge of $\tilde{C}$ is of direction $\sum_{m\in
  I}(-\alpha_m+k_m)$ for some subset $I\subset \{1,\ldots,r\}$.  
Edges adjacent to $T$ correspond to disjoint subsets of $\{1,\ldots,i_0-1\}$
whose union equals $\{1,\ldots,i_0-1\}$.  
These edges thus hand us a way to group the summands of
$\sum_{j=1}^{i_0-1}(-\alpha_j+k_j)$ in such a way that the sum of each
group is positive, even though the single summands $-\alpha_j+k_j$
might be negative. The summands $-\alpha_j+k_j$ for $j\geq i$ are
positive, too. 
For example in Figure \ref{fig:orientation} we have $ (-\alpha_3+k_3)<0$,
but $(-\alpha_1+k_1)+ (-\alpha_3+k_3)>0$.

Thus we grouped the whole sum $\sum_{j=1}^r (-\alpha_j+k_j)$ into positive summands, and in total we get $1$. 
Since by assumption we have $-\alpha_r+k_r>0$, we deduce that   
$r=1$, $\alpha_1=-1$, $k_1=0$  and $\tilde{C}$ is just a end, namely the end of direction $(1,-n,1)$. Thus $-\alpha_1-l_1=-n$, $l_1=n+1$ and $\beta_1=-\alpha_1+k_1+l_1=n+2$.
In $C_1$, we thus have an edge meeting $L$ with direction
$(1,1,n+2)$. 
It follows that $h$ is as depicted in Figure \ref{fig:4types}c in a neighborhood of $V$.

From now on we  assume that $\alpha_i\geq 0$ for all vertices $V_i$ in
$\tilde{C}$. 
The ordering of the slopes also shows us that for the left-most
vertex, $V_1$, we must have $-\alpha_1+k_1=\beta_1-l_1>0$,
since otherwise
there would be no way to connect $V_1$ to the end of $\tilde C$.
Since $l_1\geq 0$ it follows that
$\beta_1>0$, and so that $\beta_i>0$ for $i=1\ldots r$.  
Consider the $y$-coordinate $-\alpha_i-l_i$ of the edge of $\tilde{C}$
adjacent to $V_i$. Since $\alpha_i\geq 0$ and $l_i\geq 0$ it is
non-positive. The sum $\sum_{i=1}^r(-\alpha_i-l_i)$ equals 
$-n$. Since each summand is non-positive, we can conclude that each summand is bigger or equal to $-n$ and thus also $0\leq \alpha_i\leq n$.
From the above, we know that $\sum (k_i+l_i)=n+1$, so in particular
each summand $k_i+l_i$ has to be smaller or equal to $n+1$. 
Thus $h$ is as depicted in Figure \ref{fig:4types}d in a neighborhood of $V$.
\end{proof}

\begin{definition}\label{def:roof}
For an element $h:C\rightarrow X$ in $\TT\mathcal{S}(\omega)$, denote by
$C'\subset C$ the union of all connected components of $C\cap (\sigma_2\cup \sigma_3)$ containing a vertex of type d
 as in Lemma
\ref{lem-vertextypes}.
Let $h':C'\to \RR^2$ be the
composition $h'=h\circ \pi$, where $\pi$ denotes the projection to the
first two coordinates of $\RR^3$. The map $h'$ is a tropical morphism,
and we call it 
the \emph{roof} of $h:C\rightarrow X$. 
\end{definition}

Note that it follows immediately from Lemma \ref{lem-vertextypes} that
the
 roof of any element $h:C\rightarrow X$ in $\TT\mathcal{S}(\omega)$ is in 
$\TT\mathcal S'(a-m,n,d,\alpha,\omega')$ (see Definition \ref{def:trop
   invmult}) for 
values $m$, $n$, $d$ and $\alpha$ coming from a fan $\delta \vdash
\delta_0$ (see \ref{deltazero}), and $\omega'$ determined by $h_1:C_1\to\RR^2$.

\begin{proof}[Proof of Theorem \ref{main formula}]
As before (see Notation \ref{deltazero} and Section \ref{sec:corr}), let $$\delta_0=\{(1,n)^a, (0,-1)^{an+b}, (-1,0)^a, (0,1)^b \}\mbox{, and}$$ 
$$\Delta= \{(1,-n,1)^a, (0,1,1)^{an+b}, (-1,0,0)^a, (0,-1,0)^b, (0,0,-1)^{a(n+1)+b} \}.$$
It follows from the Correspondence Theorem \ref{thm-corres} that 
$$ N_{2\chi}(\delta_0)=\TT N_\chi(\Delta)=\TT N_\chi(\Delta,\omega)$$ 
where $\omega\subset X$ is a configuration of points as before,
i.e.\ in general position in $\sigma_1$ and very far down
 from $L$ (see Notation \ref{not-pointslow}).
Let $h:C\to X$ be in $\TT\mathcal S(\omega)$, i.e.\ a morphism to $X$ contributing to $\TT N_\chi(\Delta,\omega)$. 
We first want to show how we can split the information of $h$ into data as required by the right hand side of the equation.
We keep using Notation \ref{not-ci} and \ref{enum pb}. 
By abuse of notation, we often do not distinguish between the tropical morphism $h_1:C_1\to\mathbb \RR^2$ to the plane and the restriction of $h_1:C_1\to \sigma_1$: if we speak about directions of edges, we use two coordinates $x$ and $z$ in both cases. In the first case, these denote the two coordinates of $\RR^2$, in the second case, this is a shortcut for the direction $(x,x,z)$ in $\sigma_1$.

Lemma \ref{lem-vertextypes} tells us what Newton fan $\delta$ of 
$h_1:C_1\to\mathbb \RR^2$, viewed as a plane curve by prolonging the
edges that meet $L$, has:  
\begin{itemize}
 \item $a(n+1)+b$ ends of direction $(0,-1)$, since these are just the ends of $C$ in $\sigma_1$,
\item  $m\leq a$ ends of direction $(1,n+2)$ which become  vertices as
  in Figure \ref{fig:4types}c when meeting $L$,
\item  $A$ ends of direction $(-1,0)$ which become vertices as
  in Figure \ref{fig:4types}a when meeting $L$,
\item $U$ ends of direction $(0,1)$, $B$ of them become vertices as in Figure \ref{fig:4types}b when meeting $L$, $U-B$ become vertices as in Figure \ref{fig:4types}d. 
\item  ends of direction $(-\alpha_i,\beta_i)$ satisfying $n\geq \alpha_i\geq 0$, $0<\beta_i$ and $1<\beta_i$ if $\alpha_i=0$ which  become vertices as
  in Figure \ref{fig:4types}d.
\end{itemize}
Hence $h_1:C_1\to\RR^2$ has Newton fan 
$$\begin{array}{ll}\delta= & 
\{(1,n+2)^m, (0,-1)^{a(n+1)+b}, (-1,0)^A, 
\\ &\quad \quad (-\alpha_1,\beta_1),\ldots,
(-\alpha_r,\beta_r), (0,\beta_{r+1}),\ldots,
(0,\beta_{r+s}), (0,1)^U \}
\end{array}$$
i.e.\ with $\delta \vdash \delta_0$. 

We must have $A+B=b$ since this is the total number of ends with
direction $(0,-1,0)$ of $h$.
Since the total number of ends of direction $(0,1)$ in $\delta$ is
$U$, and since $B$ of these become vertices as in Figure \ref{fig:4types}b, we have $U-B=U+A-b$ ends of direction
  $(0,1)$ that are  adjacent to the roof $C'$ of $C$.
We clearly have $\chi(C)=\chi(C_1)+\chi(C')-\#C_1\cap C'$, so
since any irreducible component of  $C'$ is a tree we get
$$\chi(C_1)=\chi(C) - (a-m) + (r+s+U+A-b)=\chi-(a+b-m-r-s-A-U).$$

By Lemma \ref{lem-vertextypes}, 
the roof $h':C'\to \RR^2$ of $h:C\to X$ is an element of the set $\TT \mathcal
S'(a-m,n,d,\alpha,\omega')$ (see Definition \ref{def:trop invmult}) with 
 $$d=
(\beta_1+\alpha_1,\ldots,\beta_r+\alpha_r,\beta_{r+1},\ldots,\beta_{r+s} ,
1^{U+A-b})
\mbox{ and }
 \alpha= (\alpha_1, \ldots,\alpha_r )$$
and where $\omega'$ is determined by $h_1:C_1\to\RR^2$ and the choice of $B=b-A$ ends of
$C_1$ which become vertices as in Figure \ref{fig:4types}b when meeting $L$ and accordingly are not part of the roof.
We define
$$\mu_{h_1}=\prod_{V\in\Ve_{\sigma_1}(C)}\mu_V $$
where $\mu_V$ has been defined in Definition \ref{def:mult}. So by
definition
we have 
$$\mu_h=\prod_{i=1}^r\gcd(\alpha_i,\beta_i) \cdot \prod_{i=1}^s\beta_{r+i}\cdot \mu_{h_1}\mu_{h'}.$$

Conversely given a Newton fan $\delta$ as above, we denote by
$\TT\mathcal S''(\omega)$ the set of all tropical morphisms $h_1:C_1\to
\RR^2$ with Newton 
fan
$\delta$, Euler characteristic
$\chi'=\chi-(a+b-m-r-s-A-U) $, and 
passing through all points in $\omega$. According to the
Correspondence Theorems in \cite{Mi03} and \cite{Shu12} 
(see also \cite{GM052}), we have
$$\sum_{h_1\in \TT\mathcal S''(\omega)} \mu_{h_1}= N_{2\chi'}(\delta).$$
Now choose any element $h_1:C_1\to\RR^2$ of $\TT\mathcal S''(\omega)$,
and any set of $b-A$ ends of
$C_1$ with direction $(0,1)$. These $b-A$ ends define a set  $\TT\mathcal
S'(a-m,n,d,\alpha,\omega')$, and given any of its element $h':C'\to\RR^2$, one
can reconstruct a unique tropical morphism $h:C\to X$ in $\mathcal
C(\omega)$ reversing the construction above.
\end{proof}

In the following, we demonstrate that one can, with a little more care, also use our methods to count
irreducible curves. We generalize to any genus a
formula previously proved by Abramovich and Bertram for rational
curves (see \cite{Vak00b}).

Given a tuple $(l_1,\ldots,l_k)$ of positive integers, we denote by
$ \tau_{l_1,\ldots,l_k}$ the number of its symmetries, i.e. if 
$(l_1,\ldots,l_k)$ contains exactly $n_i$ times the entry $i$ then
$$\tau_{l_1,\ldots,l_k}=\prod_i n_i!. $$

\begin{proposition}\label{prop:application}
Let $n,b,g\ge 0$ be integers, and 
$$\delta_{n,b}=\{(1,n)^2,(0,-1)^{2n+b},(-1,0)^2,(0,1)^b\}.$$
Then
\begin{align*} &N^{\irr}_{2-2g}(\delta_{n,b})=
N^{\irr}_{2-2g}(\delta_{n+2,b-2})+\\&
\sum_{l_1+\ldots +l_{g+1}=g+1}^{n+1} \binom{2(n+b)+3}{n+1-\sum l_i
}\cdot \frac{(b+g-1)!\cdot\prod_{i=2}^{g+1} l_i^2}{(b-1)!}
\cdot
\\& \hspace{3cm}\cdot\left(\frac{(b+g)l_1^2}{\tau_{l_1,\ldots,l_{g+1}}}+
\frac{1}{\tau_{l_2,\ldots,l_{g+1}}}\cdot\binom{l_1}{2}\right).\end{align*}
\end{proposition}
\begin{proof}
To prove this equation, we mainly apply the techniques used in the
proof
of  Theorem \ref{main formula}, taking into account
irreducibility issues. 

First of all, it follows from the Correspondence Theorem
\ref{thm-corres}
 that $N^{\irr}_{2-2g}(\delta_{n,b})=\TT N^{\irr}_{1-g}(\Delta)$  where
$$\Delta= \{(1,-n,1)^2, (0,1,1)^{2n+b}, (-1,0,0)^2, (0,-1,0)^b,
 (0,0,-1)^{2(n+1)+b} \}.$$ 
Just as in Notation \ref{not-pointslow}, we choose a generic configuration of
$2(n+b) +3 +g$ points $\omega$ 
very far down in $\sigma_1$. 
As in the proof of Theorem \ref{main
  formula}, for any element $h:C\to X$ of 
$\TT \mathcal S^{\irr}(\omega)$, we have to understand the 
 contributions of $h_1:C_1\to\RR^2$ and the roof $h':C'\to\RR^2$. 
Let $\delta$ be the Newton fan of $h_1$.
Since $a=2$ in our case, there are not many possibilities
for $\delta$, the Newton fan of $h_1$:
$$\begin{array}{ll}\delta= & 
\{(1,n+2)^m, (0,-1)^{2(n+1)+b}, (-1,0)^A, 
\\ &\quad \quad (-\alpha_1,\beta_1),\ldots,
(-\alpha_r,\beta_r), (0,\beta_{r+1}),\ldots,
(0,\beta_{r+s}), (0,1)^U \}
\end{array}$$
with $m\leq 2$ (recall that $\delta\vdash\delta_{n,b}$, see Notation
\ref{deltazero}).  

\vspace{2ex}
Let us assume first that $m=2$. Recall that $\mathcal N(0,n,d,\alpha)\ne 0$
if and only if $d=\alpha=0$,  in which case it is equal to 1.
This implies that  $r=s=0$, which in its turn gives $A=2$ and
$U=b-2$.
Hence the contribution of all elements of $\TT \mathcal
  S^{\irr}(\omega)$ such that $m=2$ is equal to $N^{\irr}_{2-2g}(\delta_{n+2,b-2})$.

\vspace{2ex}
Let us assume  that $m=0$. In this case, the Newton polygon of $h_1$
is the segment with endpoints $(0,0)$ and $(0,2n+2+b)$. In particular the
space of deformations of $h_1$ has dimension $2n+2+b<2(n+b) +3 +g$, so
$h_1(C_1)$ 
cannot pass through all points in $\omega$ since this
latter configuration is generic.

\vspace{2ex}

Hence we are left to study the case $m=1$.
Let us first assume $A=1$, which is equivalent to $r=0$. So the roof
$h'$ looks like in Figure \ref{fig:easy 2nd kind}. 
Because of its Newton polygon  (see Figure \ref{fig:NP2}a), 
 the tropical
morphism $h_1$ has one irreducible rational component $h_0:C_0\to\RR^2$ 
with Newton polygon
depicted in Figure \ref{fig:NP2}b, and all the other irreducible
components have the segment $[0,1]$ as Newton polygon, i.e.\ they are vertical lines of weight one. Since $C$ is
irreducible and of genus $g$, the roof of $C$ has to connect $g+1$ vertical
ends of $C_0$ with all the other irreducible
components of $C_1$. We depict, in a floor diagram style (see
\cite{BM08}, \cite{BM}), in Figure \ref{fig:NP2}c 
how the morphism $h_1:C_1\to\RR^2$ looks like.
\begin{figure}[ht]
\centering
\begin{tabular}{ccc}
\input{NP2a.pstex_t}& 
\input{NP2b.pstex_t}&
 \input{FD.pstex_t}
\\ a) & b) &c)
\end{tabular}
\caption{The case $m=A=1$.}
\label{fig:NP2}
\end{figure}
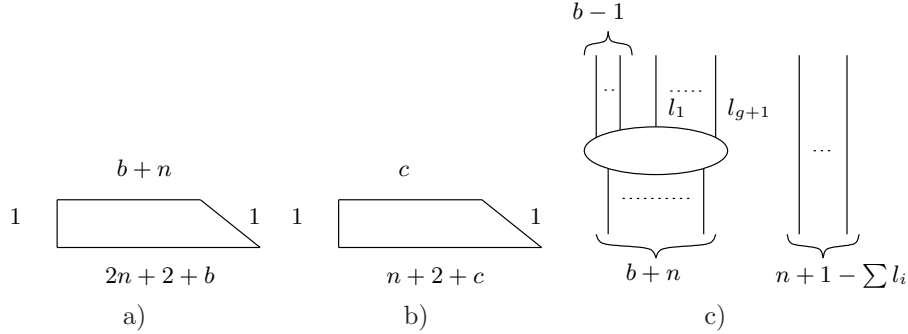
The roof of $C$ connects $g+1$ vertical ends of $C_0$ with
weights $l_1,\ldots,l_{g+1}$. $C_0$ contains $b-A=b-1$ ends of direction $(0,1)$ which become vertices of type b as in Figure \ref{fig:4types}b. Therefore we have $c=b-1+\sum l_i$, and there are $n+b-c=n+1-\sum l_i$ irreducible components which are vertical lines.
There are $\binom{2(n+b)+3}{n+1-\sum l_i}$ distinct choices for the
points of $\omega$ through which pass the vertical lines. 
The number of tropical morphisms $h_0:C_0\to \RR^2$ passing through the remaining points times the number of ways to choose the ends $l_i$ attached to the roof equals 
$$\frac{(b+g)!}{(b-1)!\tau_{l_1,\ldots,l_{g+1}}}.$$
In any case we have $\mu_{h_1}=\prod l_i$ and the roof $h'$ contributes $1$, so altogether the
contribution of all tropical morphisms in $\TT\mathcal
  S^{\irr}(\omega)$ such that $m=A=1$ is equal to 
$$\sum_{l_1+\ldots+ l_{g+1}=g+1}^{n+1} \binom{2(n+b)+3}{n+1-\sum l_i
}\cdot \frac{(b+g)!}{(b-1)!\tau_{l_1,\ldots,l_{g+1}}}
\cdot\prod_{i=1}^{g+1} l_i^2.$$

\vspace{1ex}
Finally, it remains to consider summands with $m=1$ and $A=0$,
i.e. $r=1$.
We denote 
$(-\alpha_1,\beta_1)=(-1,l_1-1)$, so the roof
$h'$ looks like in Figure \ref{fig:NP3}. 
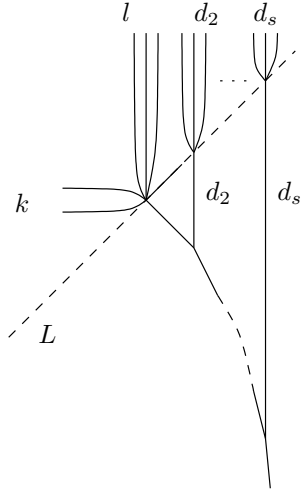
\begin{figure}[ht]
\centering
\input{NP3.pstex_t}
\caption{The case $m=1$, $A=0$.}
\label{fig:NP3}
\end{figure}
Since the end of $C_2$ has direction $(1,-n-2)$, we have $k=2$.
Just as before 
the roof of $C$ connects $g$ vertical ends of $C_0$ with
weights $l_2,\ldots,l_{g+1}$, and all the other irreducible components of
$C_1$ are just vertical lines of weight one. As before there are $n+1-\sum l_i$ such components.
There are $\binom{2(n+b)+3}{n+1-\sum l_i}$ distinct choices for the
points of $\omega$ through which pass the vertical lines. 
The number of tropical morphisms $h_0:C_0\to \RR^2$ passing through the remaining points times the number of ways to choose the ends $l_i$ attached to the roof equals 
$$\frac{(b+g-1)!}{(b-1)!\tau_{l_2,\ldots,l_{g+1}}}.$$

We have $\mu_{h_1}= \prod_{i=2}^{g+1} l_i$ and $\mu_{h'}=\binom{l_1}{2}$, so altogether the
contribution of all tropical morphisms in $\TT\widetilde{\mathcal
  S}(\omega)$ such that $m=1$ and $A=0$ is equal to 
$$\sum_{l_1+\ldots +l_{g+1}=g+1}^{n+1} \binom{2(n+b)+3}{n+1-\sum l_i
}\cdot \frac{(b+g-1)!}{(b-1)!\tau_{l_2,\ldots,l_{g+1}}}
\cdot\prod_{i=2}^{g+1} l_i^2\cdot\binom{l_1}{2}.$$
This completes the proof.
\end{proof}

\begin{example}
In the case $g=0$, Proposition \ref{prop:application} reduces to
\begin{align*} &N^{\irr}_{2}(\delta_0)=
N^{\irr}_{2}(\delta)+
\sum_{l=1}^{n+1} \binom{2(n+b)+3}{n+1-l}
\cdot\left(b\cdot l^2+
\binom{l}{2}\right).\end{align*}
which has been first proved by Abramovich and Bertram (see \cite{Vak00b}).
\end{example}

\section{Proofs of  Correspondence theorems}\label{sec-corres}

\input{proof.tex}

\appendix
\section{Hirzebruch surfaces and their deformations}\label{hirz}
In this appendix we translate to the tropical setting Kodaira
deformation of Hirzebruch surfaces. 
We include it since
it provides a justification for the strategy
of the proof of Theorem \ref{main formula}, however the
present paper is formally independent
from this appendix. As a consequence, we only sketch the following proofs, and we assume that the reader is well acquainted with tropical geometry.

We first recall Kodaira deformation of Hirzebruch surfaces before
turning to tropical deformations. 

\subsection{Kodaira deformation of Hirzebruch surfaces and
  deformation to the normal cone}\label{subsec-deform1}
As explained in Section \ref{def hirz}, 
 two Hirzebruch surfaces $\Sigma_n$ and $\Sigma_{n'}$ are not
biholomorphic if $n \ne n'$. However, if $n$ and $n'$ have the same
parity, one can deform one of the two surfaces to the other one.

\begin{theorem}[Kodaira, see \cite{Kod86}]\label{def cplx}
Let $n,k\ge 0$ be two integer numbers. There exists a complex manifold
$X_{n,k}$ of dimension 3 equipped with a submersion 
$\phi_{n,k}:X_{n,k}\to\CC$
such that
$$\forall t\ne 0,\  \phi_{n,k}^{-1}(t)=\Sigma_n, \quad \text{and}\quad
\phi_{n,k}^{-1}(0)=\Sigma_{n+2k}.$$  
\end{theorem}
Note that this implies that $\Sigma_n$ and $\Sigma_{n+2k}$ are
diffeomorphic. We will prove Theorem \ref{def cplx} in Section
\ref{trop hirz}
in the tropical language. However, our proof translates literally
to the complex setting.

\vspace{1ex}
Kodaira deformation of Hirzebruch surfaces can actually be reduced
to a standard procedure in both complex (deformation to the normal
cone, see \cite{F}) and symplectic (symplectic sum, see \cite{Ler95}
and \cite{IP00}, stretching the neck, see \cite{EGH})
geometries. 

\begin{proposition}[see \cite{F}]\label{def normal}
Let $S$ be a nonsingular complex surface, let  $C$ be a nonsingular algebraic
curve in $S$, and let  $X'_{S,C}$ be the blow-up of the curve
$C\times\{0\}$ in the complex 3-fold $S\times\CC$. Then the exceptional
divisor of $X'_{S,C}$ is isomorphic to $\mathbb P(\mathcal
N_{C/S}\oplus\CC)$, where $\mathcal N_{C/S}$ is the normal bundle 
of $C$ in
$S$.
\end{proposition}

The first Chern class of $\mathcal N_{C/S}$ is the self-intersection of
$C$ in $S$. In particular, if $C$ is rational of self-intersection $m$, 
 then $\mathbb P(\mathcal
N_{C/S}\oplus\CC)=\Sigma_{|m|}$.

We denote by $\phi'_{(S,C)}:X'_{S,C}\to\CC$ the obvious
(holomorphic) projection on the $\CC$ factor. Then Proposition
\ref{def normal} implies that
$\phi'_{(S,C)}$ is a submersion over $\CC^*$,  that
$$\forall t\ne 0,\  \phi_{(S,C)}'^{-1}(t)=S,$$
and that
$\phi_{(S,C)}'^{-1}(0)$ is the union of $S$ and $\mathbb P(\mathcal
N_{C/S}\oplus\CC)$
 intersecting transversely along $C$.

Suppose now that $S=\Sigma_n$ and $C$ is a smooth rational curve of
bidegree $(1,k)$. Since $C$ has self-intersection $n+2k$ in
$\Sigma_n$, we get from Proposition \ref{def normal} that
$$\phi_{(\Sigma_n,C)}'^{-1}(0)=\Sigma_n\cup \Sigma_{n+2k}. $$
In this case, it turns out that the two complex 3-folds, $X_{n,k}$ from Theorem
\ref{def cplx} and $X'_{\Sigma_n,k}$ from Proposition
\ref{def normal}, are related by a blow-up: one can contract the copy
of $\Sigma_n$ in $\phi_{(\Sigma_n,C)}'^{-1}(0)$ to obtain $X_{n,k}$.

\begin{proposition}\label{blowup cplx}
The complex manifold $X'_{\Sigma_n,C}$ is the blow-up $\mbox{bl}$ of
$X_{n,k}$ along the exceptional curve $E_{n+2}$ of $\Sigma_{n+2k}$, and 
$$\phi_{(\Sigma_n,C)}'=\phi_{n,k}\circ\mbox{bl}  .$$
\end{proposition} 
We prove the tropical analogue of Proposition \ref{blowup cplx} in Section \ref{trop
  hirz}. Once again, our tropical proof translates literally to the
complex setting.

\subsection{Tropical Hirzebruch surfaces}\label{trop hirz}
The construction of any non-singular toric variety can be performed exactly
in the same way in tropical and algebraic geometry. 
In particular, the tropical  Hirzebruch
surface of degree $n\ge 0$, denoted by $\TT \Sigma_n$ is constructed  
by taking two copies of $\TT\times \TT P^1$
 glued
along $\TT^*\times \TT P^1$
via the tropical isomorphism
$$ 
\begin{array}{ccccc}
\psi :&\TT^*\times \TT P^1 &\longrightarrow & \TT^*\times \TT P^1
\\ &(x_1,y_1) &\longmapsto & (\tg \frac{1}{x_1}\td ,\tg
\frac{y_1}{x_1^n}\td)&=  (-x_1 ,y_1 -nx_1)
\end{array}
$$
As in the complex setting, $\TT \Sigma_0=\TT P^1\times \TT P^1$, and 
$\TT \Sigma_1$ is $\TT P^2$ blown up at $[-\infty:0:-\infty]$.

\begin{figure}[ht]
\centering
\begin{tabular}{ccc}
\includegraphics[height=1.5cm, angle=0]{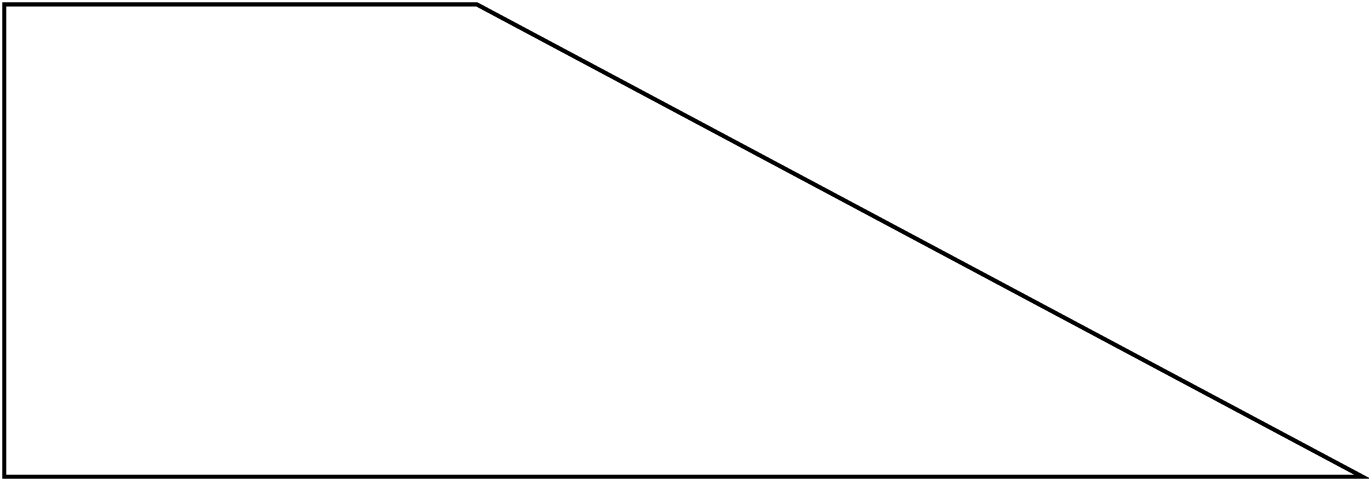} &
\includegraphics[height=2cm, angle=0]{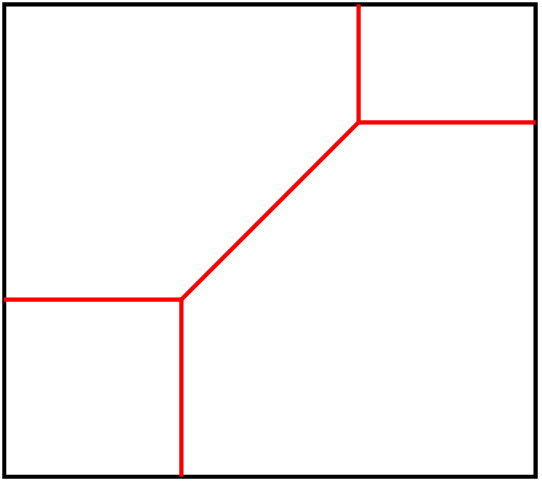} &
\includegraphics[height=2cm, angle=0]{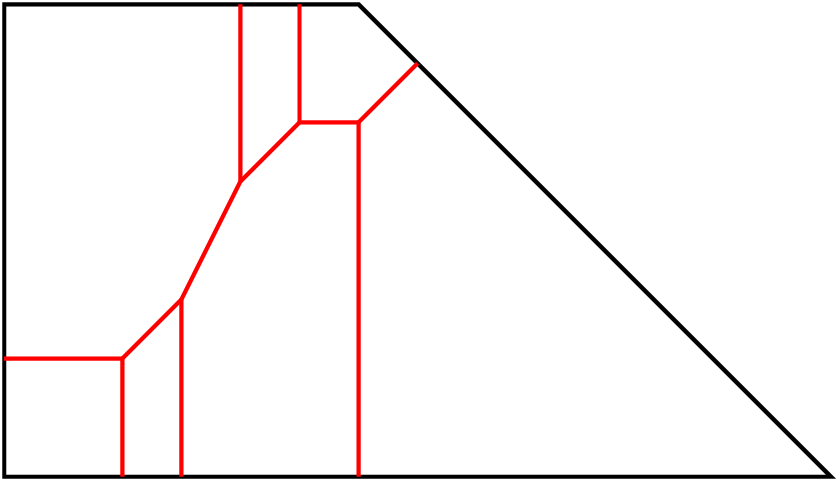}
\\
\\ a)  $\TT \Sigma_n$ & b) A curve of bidegree 
& c) A curve of bidegree 
\\ & $(1,1)$ in
$\TT\Sigma_0$   & $(1,2)$ in $\TT\Sigma_1$ 
\end{tabular}
\caption{Tropical toric surfaces and embedded curves.}
\label{TSn}
\end{figure}

The map $\psi$ 
sends
the vector $(1,n)$ to the vector $(-1,0)$. For
this reason, the tropical surface $\TT \Sigma_n$ is usually
represented by a quadrangle with two horizontal edges, one vertical
edge, and one edge of slope $-\frac{1}{n}$, see Figure \ref{TSn}a. More
generally, once a linear system is fixed, 
the tropical moment map provides an homeomorphism from any
 non-singular toric tropical variety to any Newton polygon
 corresponding to the linear system. This
 homeomorphism is also given by the Veronese embedding corresponding
 to the chosen linear system.

Let us denote by $B_n$ the tropical curve in $\TT \Sigma_n$ defined by the
tropical polynomial $\tg y_1\td$, and by  $F_n$ the curve defined by the
polynomial $\tg x_1\td$ (in the coordinate system
defined above on $\TT\Sigma_n$). That is to say,
the curve $B_n$
is the lowest horizontal edge, and  $F_n$ 
is the left vertical edge. 
 Note that $B_n^2=n$ and $F_n^2=0$.
The tropical Picard group of $\TT\Sigma_n$ is the free
abelian group of rank two generated by $B_n$ and $F_n$, and a tropical
1-cycle $C$ in $\TT\Sigma_n$ is said to have bidegree $(a,b)$ if 
it is linearly equivalent to $aB_n+bF_n$. Equivalently, $C$ is of bidegree
$(a,b)$ if and only if
$C\circ B_n=an +b$ and $C\circ F_n=a.$
\begin{example}
We depicted in Figure \ref{TSn}b a tropical curve of bidegree $(1,1)$ in
$\TT\Sigma_0$, and a tropical curve of bidegree $(1,2)$ in
$\TT\Sigma_1$ in Figure \ref{TSn}c.
\end{example}

The exceptional divisor $E_n$ of $\TT\Sigma_n$ is the upper horizontal
edge, defined by the rational function $\tg \frac{1}{y_1}\td$, and
represents the class $B_n-nF_n$ in $\mbox{Pic}(\TT\Sigma_n)$. In particular, 
one
has $E_n^2=-n$.

\subsection{Deformation of tropical Hirzebruch
  surfaces}\label{subsec-tropdeform} 

Let $C$ be a non-singular tropical curve in $\TT\Sigma_n$ of bidegree
$(1,k)$. By the genus formula $C$ is rational. Moreover we have
$C^2=n+2k$.

We start by  describing the deformation of $\TT\Sigma_n$ 
to the normal cone of $C$. Recall that the sign ``$=$'' between two
tropical varieties means ``isomorphic up to tropical modifications''
(see \cite{Mi06}). Note that we have to introduce tropical
modifications at this point since Kodaira deformation of Hirzebruch
surfaces in non-toric. Indeed the exceptional divisor is a toric
divisor, so 
any Hirzebruch surface is torically rigid.

\begin{theorem}\label{def normal trop}
There exists a non-singular tropical variety $\TT X'$ of dimension 3
and a tropical morphism $\Phi':\TT X'\to \TT P^1$ such that
$$\forall t\ne -\infty,\  \Phi'^{-1}(t)=\TT\Sigma_n, 
\quad \text{and}\quad
\Phi'^{-1}(-\infty)=\TT\Sigma_{n}\cup\TT\Sigma_{n+2k}.$$
Moreover, the intersection curve of
the two latter surfaces is 
 $C$ in $\TT\Sigma_{n}$, and the exceptional section $E_{n+2k}$ 
in $\TT\Sigma_{n+2k}$.
\end{theorem}
\begin{proof}
Let $\Pi'$ be the polytope in $\RR^3$ with vertices (see Figure \ref{X'})
$$(0,0,0),\ (2n+2k,0,0),\ (2k,2,0),\ (0,2,0), \ (0,0,1),\ 
(n+k,0,1),\ (1,k,1),\ (0,1,1).$$
The polytope $\Pi'$ defines a non-singular tropical toric variety
$\mbox{Tor}(\Pi')$ of
dimension 3. If $(x,y,z)$ are the
coordinates in the dense $\RR^3$-orbit of $\mbox{Tor}(\Pi')$, then the map
$(x,y,z)\to (x,y)$
induces a tropical morphism 
$\pi :\mbox{Tor}(\Pi')\to \TT\Sigma_n$ whose fibers are $\TT P^1$.
\begin{figure}[ht]
\centering
\begin{tabular}{c}
\includegraphics[height=3cm, angle=0]{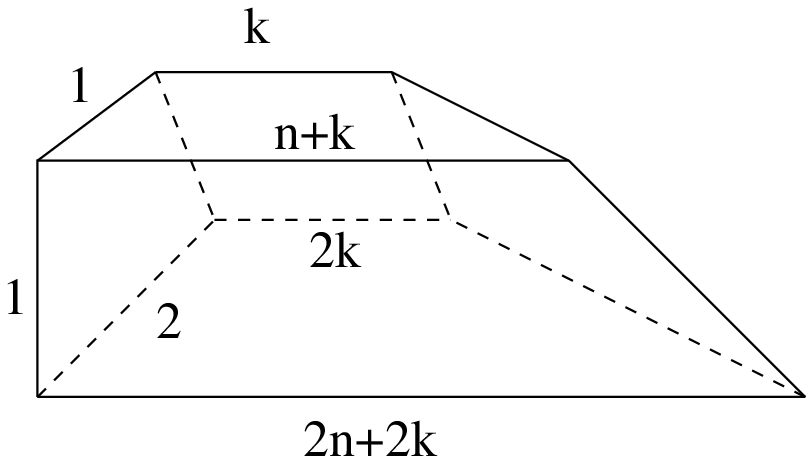} 
\\  
\end{tabular}
\caption{The polytope  $\Pi'$.}
\label{X'}
\end{figure}

We fix a tropical polynomial $p(x,y)$ defining the curve $C$ in 
$\TT\Sigma_n$, and we
define $\TT X'$ as the hypersurface in $\TT P^1\times \mbox{Tor}(\Pi') $
defined by the tropical polynomial 
$$\tg tz + p(x,y) \td$$
where $t$ is the coordinate in $\TT P^1$.

The tropical variety $\TT X'$ is non-singular, and there is a natural
tropical morphism $\Phi': \TT X'\to  \TT P^1$ whose fiber over $t_0$ is
the tropical hypersurface in $\mbox{Tor}(\Pi')$ defined by the tropical
polynomial
$\tg t_0z + p(x,y) \td$ (see Figure \ref{Modif S0}a
 in the case $n=0$ and $k=1$).
If $t_0\ne -\infty$, then $\Phi'^{-1}(t_0)$ is nonsingular, and the
morphisms
$\pi_{|\Phi'^{-1}(t_0)}: \Phi'^{-1}(t_0)\to \TT\Sigma_n$ is a tropical
modification of $\TT\Sigma_n$ along $C$. 
In particular $\Phi'^{-1}(t_0)= \TT\Sigma_n$.
The hypersurface $\Phi'^{-1}(-\infty)$ is the union of the tropical
surface $S_0$ defined by $p(x,y)$ with the surface $S_1$ in $\mbox{Tor}(\Pi')$
defined by the rational function $\tg \frac{1}{z}\td$ (i.e. the upper
horizontal face of $\mbox{Tor}(\Pi')$), see Figure \ref{Modif
  S0}b. The surfaces $S_0$ and $S_1$ intersect along a tropical curve $E$.

\begin{figure}[ht]
\centering
\begin{tabular}{cc}
\includegraphics[height=3cm, angle=0]{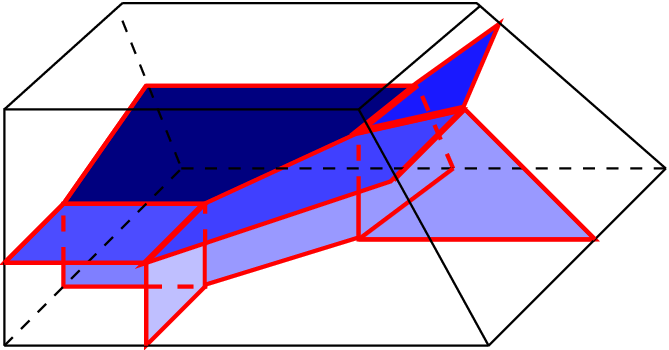} &
\includegraphics[height=3cm, angle=0]{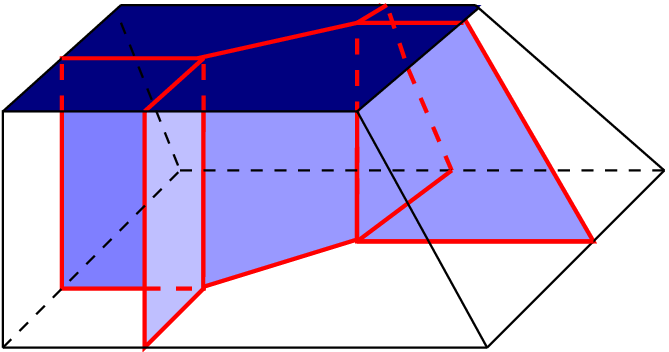}
\\
\\ a) $\Phi'^{-1}(t_0)$ with $t_0\ne -\infty$
& b) $\Phi'^{-1}(-\infty)$
\end{tabular}
\caption{The deformation process.}
\label{Modif S0}
\end{figure}

It is clear from $\Pi'$ that $S_1=\TT\Sigma_n$, so it remains to
prove that $S_0=\TT\Sigma_{n+2k}$. First, the morphisms $\pi$
restricts to a tropical morphism $S_0\to C$ whose fibers are $\TT
P^1$. It follows from elementary tropical intersection theory that
$C^2=n+2k=-E^2$ in $S_0$. 
For example, by adding a vertical edge to each tropical intersection
points of $C$ in $\TT\Sigma_n$, we see  that the
self-intersection of $C$ in $\TT\Sigma_n$ and $S_0$ are the same, i.e.
 $S_0=\TT\Sigma_{n+2k}$ (see Figure \ref{self int}).
Note that for the same reasons, we have $E^2=-n-2k$ in $S_0$.
\begin{figure}[ht]
\centering
\begin{tabular}{c}
\includegraphics[height=3cm, angle=0]{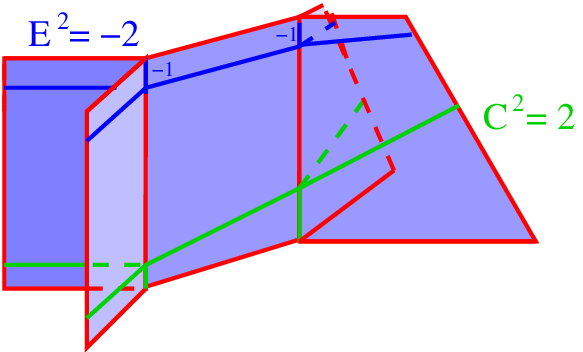} 
\end{tabular}
\caption{Self intersection of $C$ and $E$ in $S_0$.}
\label{self int}
\end{figure}
\end{proof}

\begin{example}
We depict in Figure \ref{Modif S0} 
this degeneration 
process when $n=0$, $k=1$,
and $C$ is the tropical curve depicted in Figure \ref{TSn}b.
\end{example}

In the next lemma, we use notation introduced in the proof of Theorem
\ref{def normal trop}. 
\begin{lemma}\label{blow down}
One can blow down $S_1$ to $E$ in $\TT X'$, i.e. there exists a non-singular
tropical variety $\TT X$ of dimension 3 and a tropical blow down
$\mbox{bl}:\TT X'\to\TT X$ which contracts the surface $S_1$ to a curve
isomorphic (up to tropical modifications) to $E$. 
\end{lemma}
\begin{proof}
This is an immediate consequence of the fact that one can blow down
the surface $S_1$ in $\mbox{Tor}(\Pi')$. Indeed, let $\Pi$ be the
polytope with vertices (see Figure \ref{blow}a)
$$(0,0,0),\ (2n+2k+1,0,0),\ (2k+1,2,0),\ (0,2,0), \ (0,0,1),\ 
(1,0,1).$$
The tropical 3-fold $\mbox{Tor}(\Pi)$ is non-singular, and the existence
of the desired blow-up $\mbox{Tor}(\Pi')\to \mbox{Tor}(\Pi)$ can be observed
directly at
the polytopes $\Pi$ and $\Pi'$ (see Figure \ref{blow}b). 
\begin{figure}[ht]
\centering
\begin{tabular}{cc}
\includegraphics[height=1.5cm, angle=0]{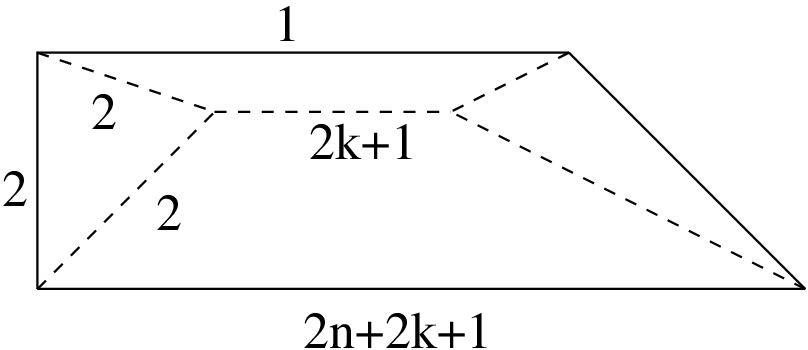} &
\includegraphics[height=1.5cm, angle=0]{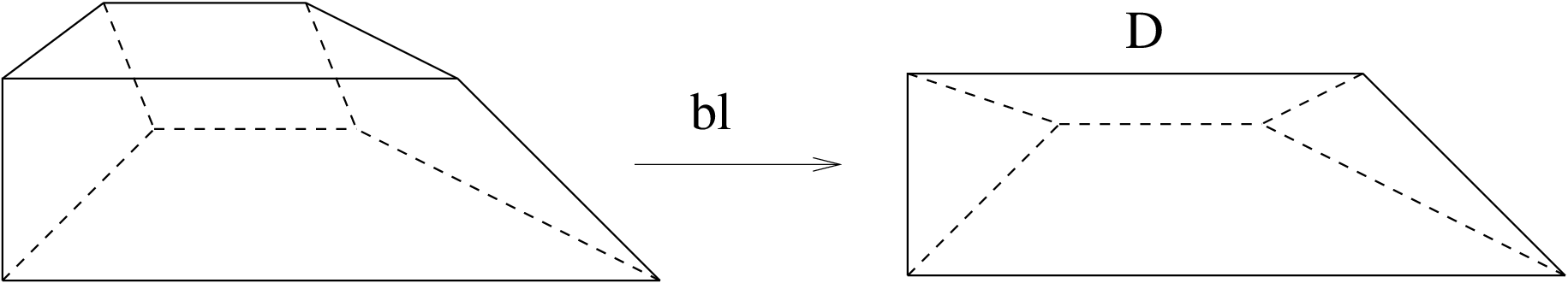}
\\
\\ a) $\Pi$
& b) Toric blow up of the curve $D$
\end{tabular}
\caption{The blow up $\mbox{Tor}(\Pi')\to \mbox{Tor}(\Pi)$.}
\label{blow}
\end{figure}
\end{proof}

As an immediate corollary, we obtain the tropical version of Kodaira
deformation of Hirzebruch surfaces.
\begin{corollary}
There exists a non-singular tropical variety $\TT X$ of dimension 3
and a tropical morphism $\Phi:\TT X\to \TT P^1$ such that
$$\forall t\ne -\infty,\  \Phi^{-1}(t)=\TT\Sigma_n, 
\quad \text{and}\quad
\Phi^{-1}(-\infty)=\TT\Sigma_{n+2k}.$$
\end{corollary}
\begin{proof}
Take $\TT X$ as in Lemma \ref{blow down}, and $\Phi:\TT X\to\TT P^1$
such that $\Phi'=\Phi\circ \mbox{bl}$.
\end{proof}

\begin{figure}[ht]
\centering
\begin{tabular}{cc}
\includegraphics[height=1cm, angle=0]{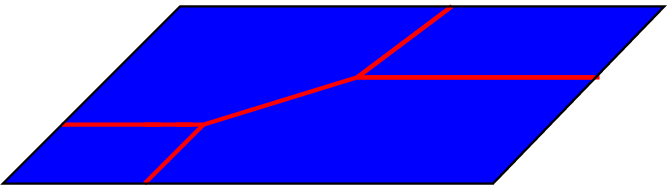} &
\includegraphics[height=3cm, angle=0]{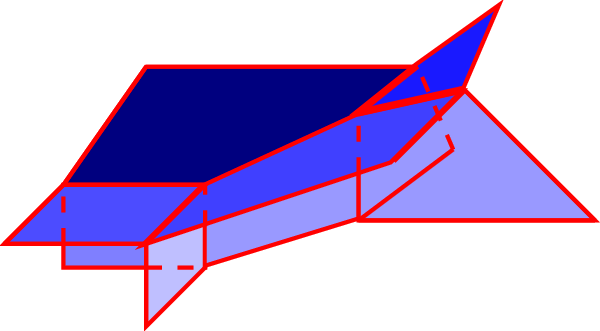}
\\
\\ a) $\Phi^{-1}(+\infty)$ 
& b) $\Phi^{-1}(t_0)$ with $t_0\ne \pm\infty$
\end{tabular}
\vspace{3ex}

\begin{tabular}{c}
\includegraphics[height=3cm, angle=0]{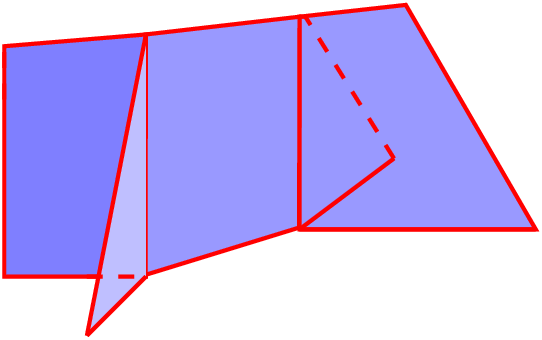} 
\\ $\Phi^{-1}(-\infty)$
\end{tabular}
\caption{Tropical Kodaira deformation.}
\label{Trop Kod}
\end{figure}

\begin{example}
We depict in Figure \ref{Trop Kod} this deformation when $n=0$, $k=1$,
and $C$ is the tropical curve depicted in Figure \ref{TSn}b.
\end{example}

\begin{remark}\label{rem:justification}
This section explains why we prove Theorem \ref{main formula} by
enumerating tropical curves in $X$. This latter surface (after the
change of coordinates $(x,y,z)\mapsto (x,-y,z)$) is the part
contained in the $\RR^3$-orbit of $\mbox{Tor}(\Pi')$ of the degeneration
process described in the proof of Theorem \ref{def normal trop}, where
$C$ is the tropical  curve of bidegree $(1,1)$
 in $\TT\Sigma_n$ defined by the tropical
polynomial $\tg 1+xy \td$. Since $X$ is the open part of a tropical
modification of $\TT\Sigma_n$, counting carefully tropical curves in
$X$ or  $\TT\Sigma_n$ should give the same numbers.

On the other hand, we choose the configuration $\omega$ inside the
face $\sigma_1$ of $X$ which degenerates to $\TT\Sigma_{n+2}$, replacing the
condition \emph{$t\to -\infty$} by the condition that \emph{points in
$\omega$  have very low $z$-coordinates}. In particular
 it is natural 
that the parts in $\sigma_1$ far from $L$ of the
 tropical curves we are counting look like curves in $\TT\Sigma_{n+2}$.
\end{remark}

\bibliographystyle {alpha}
\bibliography {bibliographie}

\end{document}

%% file: Deltadelta.pstex_t
\begin{picture}(0,0)%
\includegraphics{Deltadelta.pstex}%
\end{picture}%
\setlength{\unitlength}{3947sp}%
\begingroup\makeatletter\ifx\SetFigFont\undefined%
\gdef\SetFigFont#1#2#3#4#5{%
  \reset@font\fontsize{#1}{#2pt}%
  \fontfamily{#3}\fontseries{#4}\fontshape{#5}%
  \selectfont}%
\fi\endgroup%
\begin{picture}(4630,1661)(2086,-3448)
\put(2494,-2061){\makebox(0,0)[lb]{\smash{{\SetFigFont{9}{10.8}{\familydefault}{\mddefault}{\updefault}{\color[rgb]{0,0,0}$\binom{-\alpha_i}{\beta_i}$}%
}}}}
\put(3477,-3240){\makebox(0,0)[lb]{\smash{{\SetFigFont{9}{10.8}{\familydefault}{\mddefault}{\updefault}{\color[rgb]{0,0,0}$\binom{0}{-1}$}%
}}}}
\put(2625,-2323){\makebox(0,0)[lb]{\smash{{\SetFigFont{9}{10.8}{\familydefault}{\mddefault}{\updefault}{\color[rgb]{0,0,0}$\cdot$}%
}}}}
\put(2691,-2257){\makebox(0,0)[lb]{\smash{{\SetFigFont{9}{10.8}{\familydefault}{\mddefault}{\updefault}{\color[rgb]{0,0,0}$\cdot$}%
}}}}
\put(2769,-2201){\makebox(0,0)[lb]{\smash{{\SetFigFont{9}{10.8}{\familydefault}{\mddefault}{\updefault}{\color[rgb]{0,0,0}$\cdot$}%
}}}}
\put(3411,-1930){\makebox(0,0)[lb]{\smash{{\SetFigFont{9}{10.8}{\familydefault}{\mddefault}{\updefault}{\color[rgb]{0,0,0}$\binom{0}{1}$}%
}}}}
\put(4197,-2454){\makebox(0,0)[lb]{\smash{{\SetFigFont{9}{10.8}{\familydefault}{\mddefault}{\updefault}{\color[rgb]{0,0,0}$\binom{1}{n+2}$}%
}}}}
\put(2101,-2978){\makebox(0,0)[lb]{\smash{{\SetFigFont{9}{10.8}{\familydefault}{\mddefault}{\updefault}{\color[rgb]{0,0,0}$\binom{-1}{0}$}%
}}}}
\put(5638,-1995){\makebox(0,0)[lb]{\smash{{\SetFigFont{9}{10.8}{\familydefault}{\mddefault}{\updefault}{\color[rgb]{0,0,0}$(b-a,a)$}%
}}}}
\put(5114,-3043){\makebox(0,0)[lb]{\smash{{\SetFigFont{9}{10.8}{\familydefault}{\mddefault}{\updefault}{\color[rgb]{0,0,0}$(0,0)$}%
}}}}
\put(6424,-3043){\makebox(0,0)[lb]{\smash{{\SetFigFont{9}{10.8}{\familydefault}{\mddefault}{\updefault}{\color[rgb]{0,0,0}$(a(n+1)+b,0)$}%
}}}}
\put(4918,-1995){\makebox(0,0)[lb]{\smash{{\SetFigFont{9}{10.8}{\familydefault}{\mddefault}{\updefault}{\color[rgb]{0,0,0}$(0,a)$}%
}}}}
\end{picture}%

%% file: surfaceX.pstex_t
\begin{picture}(0,0)%
\includegraphics{surfaceX.pstex}%
\end{picture}%
\setlength{\unitlength}{3947sp}%
\begingroup\makeatletter\ifx\SetFigFont\undefined%
\gdef\SetFigFont#1#2#3#4#5{%
  \reset@font\fontsize{#1}{#2pt}%
  \fontfamily{#3}\fontseries{#4}\fontshape{#5}%
  \selectfont}%
\fi\endgroup%
\begin{picture}(1141,1656)(9150,-3456)
\put(10276,-2236){\makebox(0,0)[lb]{\smash{{\SetFigFont{10}{12.0}{\familydefault}{\mddefault}{\updefault}{\color[rgb]{0,0,0}$\sigma_2$}%
}}}}
\put(9976,-2986){\makebox(0,0)[lb]{\smash{{\SetFigFont{10}{12.0}{\familydefault}{\mddefault}{\updefault}{\color[rgb]{0,0,0}$\sigma_1$}%
}}}}
\put(9226,-2086){\makebox(0,0)[lb]{\smash{{\SetFigFont{10}{12.0}{\familydefault}{\mddefault}{\updefault}{\color[rgb]{0,0,0}$\sigma_3$}%
}}}}
\put(9901,-1936){\makebox(0,0)[lb]{\smash{{\SetFigFont{10}{12.0}{\familydefault}{\mddefault}{\updefault}{\color[rgb]{0,0,0}$L$}%
}}}}
\put(9165,-3396){\makebox(0,0)[lb]{\smash{{\SetFigFont{10}{12.0}{\familydefault}{\mddefault}{\updefault}{\color[rgb]{0,0,0}$X=V(\tg x+y+z\td)$}%
}}}}
\end{picture}%

%% file: Def.pstex_t
\begin{picture}(0,0)%
\includegraphics{Def.pstex}%
\end{picture}%
\setlength{\unitlength}{2368sp}%
\begingroup\makeatletter\ifx\SetFigFont\undefined%
\gdef\SetFigFont#1#2#3#4#5{%
  \reset@font\fontsize{#1}{#2pt}%
  \fontfamily{#3}\fontseries{#4}\fontshape{#5}%
  \selectfont}%
\fi\endgroup%
\begin{picture}(4619,3882)(5379,-3796)
\put(9676,-3286){\makebox(0,0)[lb]{\smash{{\SetFigFont{6}{7.2}{\familydefault}{\mddefault}{\updefault}{\color[rgb]{0,0,0}$\Sigma_{n+2}$}%
}}}}
\put(6151,-3286){\makebox(0,0)[lb]{\smash{{\SetFigFont{6}{7.2}{\familydefault}{\mddefault}{\updefault}{\color[rgb]{0,0,0}$\Sigma_{n+2}$}%
}}}}
\put(6376,-1861){\makebox(0,0)[lb]{\smash{{\SetFigFont{6}{7.2}{\familydefault}{\mddefault}{\updefault}{\color[rgb]{0,0,0}$V$}%
}}}}
\put(6376,-2161){\makebox(0,0)[lb]{\smash{{\SetFigFont{6}{7.2}{\familydefault}{\mddefault}{\updefault}{\color[rgb]{0,0,0}$V_\infty$}%
}}}}
\put(6076,-3736){\makebox(0,0)[lb]{\smash{{\SetFigFont{6}{7.2}{\familydefault}{\mddefault}{\updefault}{\color[rgb]{0,0,0}$\Sigma$}%
}}}}
\put(7726,-1561){\makebox(0,0)[lb]{\smash{{\SetFigFont{6}{7.2}{\familydefault}{\mddefault}{\updefault}{\color[rgb]{0,0,0}$bl$}%
}}}}
\put(6226,-61){\makebox(0,0)[lb]{\smash{{\SetFigFont{6}{7.2}{\familydefault}{\mddefault}{\updefault}{\color[rgb]{0,0,0}$\Sigma_n$}%
}}}}
\put(9076,-3736){\makebox(0,0)[lb]{\smash{{\SetFigFont{6}{7.2}{\familydefault}{\mddefault}{\updefault}{\color[rgb]{0,0,0}$\Sigma'$}%
}}}}
\put(9376,-61){\makebox(0,0)[lb]{\smash{{\SetFigFont{6}{7.2}{\familydefault}{\mddefault}{\updefault}{\color[rgb]{0,0,0}$\Sigma_n$}%
}}}}
\end{picture}%

%% file: Ex_mor.pstex_t
\begin{picture}(0,0)%
\includegraphics{Ex_mor.pstex}%
\end{picture}%
\setlength{\unitlength}{2368sp}%
\begingroup\makeatletter\ifx\SetFigFont\undefined%
\gdef\SetFigFont#1#2#3#4#5{%
  \reset@font\fontsize{#1}{#2pt}%
  \fontfamily{#3}\fontseries{#4}\fontshape{#5}%
  \selectfont}%
\fi\endgroup%
\begin{picture}(3627,2251)(6739,-3800)
\put(10351,-1711){\makebox(0,0)[lb]{\smash{{\SetFigFont{6}{7.2}{\familydefault}{\mddefault}{\updefault}{\color[rgb]{0,0,0}$\sigma_2$}%
}}}}
\put(9826,-1861){\makebox(0,0)[lb]{\smash{{\SetFigFont{6}{7.2}{\familydefault}{\mddefault}{\updefault}{\color[rgb]{0,0,0}$L$}%
}}}}
\put(9151,-3736){\makebox(0,0)[lb]{\smash{{\SetFigFont{6}{7.2}{\familydefault}{\mddefault}{\updefault}{\color[rgb]{0,0,0}$\sigma_1$}%
}}}}
\put(9301,-1786){\makebox(0,0)[lb]{\smash{{\SetFigFont{6}{7.2}{\familydefault}{\mddefault}{\updefault}{\color[rgb]{0,0,0}$\sigma_3$}%
}}}}
\put(10051,-3136){\makebox(0,0)[lb]{\smash{{\SetFigFont{6}{7.2}{\familydefault}{\mddefault}{\updefault}{\color[rgb]{0,0,0}$(1,1,-1)$}%
}}}}
\put(8476,-2686){\makebox(0,0)[lb]{\smash{{\SetFigFont{6}{7.2}{\familydefault}{\mddefault}{\updefault}{\color[rgb]{0,0,0}$(-2,0,0)$}%
}}}}
\put(7126,-1786){\makebox(0,0)[lb]{\smash{{\SetFigFont{6}{7.2}{\familydefault}{\mddefault}{\updefault}{\color[rgb]{0,0,0}$0$}%
}}}}
\put(10201,-2386){\makebox(0,0)[lb]{\smash{{\SetFigFont{6}{7.2}{\familydefault}{\mddefault}{\updefault}{\color[rgb]{0,0,0}$(1,-1,1)$}%
}}}}
\put(7201,-3136){\makebox(0,0)[lb]{\smash{{\SetFigFont{6}{7.2}{\familydefault}{\mddefault}{\updefault}{\color[rgb]{0,0,0}$1$}%
}}}}
\put(7726,-1861){\makebox(0,0)[lb]{\smash{{\SetFigFont{6}{7.2}{\familydefault}{\mddefault}{\updefault}{\color[rgb]{0,0,0}$h_1$}%
}}}}
\put(7651,-3586){\makebox(0,0)[lb]{\smash{{\SetFigFont{6}{7.2}{\familydefault}{\mddefault}{\updefault}{\color[rgb]{0,0,0}$h_2$}%
}}}}
\end{picture}%

%% file: Ex_mor2.pstex_t
\begin{picture}(0,0)%
\includegraphics{Ex_mor2.pstex}%
\end{picture}%
\setlength{\unitlength}{2368sp}%
\begingroup\makeatletter\ifx\SetFigFont\undefined%
\gdef\SetFigFont#1#2#3#4#5{%
  \reset@font\fontsize{#1}{#2pt}%
  \fontfamily{#3}\fontseries{#4}\fontshape{#5}%
  \selectfont}%
\fi\endgroup%
\begin{picture}(3627,2251)(6739,-3800)
\put(10351,-1711){\makebox(0,0)[lb]{\smash{{\SetFigFont{6}{7.2}{\familydefault}{\mddefault}{\updefault}{\color[rgb]{0,0,0}$\sigma_2$}%
}}}}
\put(9826,-1861){\makebox(0,0)[lb]{\smash{{\SetFigFont{6}{7.2}{\familydefault}{\mddefault}{\updefault}{\color[rgb]{0,0,0}$L$}%
}}}}
\put(9151,-3736){\makebox(0,0)[lb]{\smash{{\SetFigFont{6}{7.2}{\familydefault}{\mddefault}{\updefault}{\color[rgb]{0,0,0}$\sigma_1$}%
}}}}
\put(9301,-1786){\makebox(0,0)[lb]{\smash{{\SetFigFont{6}{7.2}{\familydefault}{\mddefault}{\updefault}{\color[rgb]{0,0,0}$\sigma_3$}%
}}}}
\put(10051,-3136){\makebox(0,0)[lb]{\smash{{\SetFigFont{6}{7.2}{\familydefault}{\mddefault}{\updefault}{\color[rgb]{0,0,0}$(1,1,-1)$}%
}}}}
\put(8476,-2686){\makebox(0,0)[lb]{\smash{{\SetFigFont{6}{7.2}{\familydefault}{\mddefault}{\updefault}{\color[rgb]{0,0,0}$(-1,0,0)$}%
}}}}
\put(7126,-1786){\makebox(0,0)[lb]{\smash{{\SetFigFont{6}{7.2}{\familydefault}{\mddefault}{\updefault}{\color[rgb]{0,0,0}$0$}%
}}}}
\put(10201,-2386){\makebox(0,0)[lb]{\smash{{\SetFigFont{6}{7.2}{\familydefault}{\mddefault}{\updefault}{\color[rgb]{0,0,0}$(1,-1,1)$}%
}}}}
\put(7726,-1861){\makebox(0,0)[lb]{\smash{{\SetFigFont{6}{7.2}{\familydefault}{\mddefault}{\updefault}{\color[rgb]{0,0,0}$h_3$}%
}}}}
\end{picture}%

%% file: specvert1.pstex_t
\begin{picture}(0,0)%
\includegraphics{specvert1.pstex}%
\end{picture}%
\setlength{\unitlength}{3947sp}%
\begingroup\makeatletter\ifx\SetFigFont\undefined%
\gdef\SetFigFont#1#2#3#4#5{%
  \reset@font\fontsize{#1}{#2pt}%
  \fontfamily{#3}\fontseries{#4}\fontshape{#5}%
  \selectfont}%
\fi\endgroup%
\begin{picture}(1680,1398)(3736,-4273)
\put(5401,-3286){\makebox(0,0)[lb]{\smash{{\SetFigFont{10}{12.0}{\familydefault}{\mddefault}{\updefault}{\color[rgb]{0,0,0}$L$}%
}}}}
\put(5326,-3586){\makebox(0,0)[lb]{\smash{{\SetFigFont{8}{9.6}{\familydefault}{\mddefault}{\updefault}{\color[rgb]{0,0,0}$\binom{\alpha_i}{\alpha_i}$}%
}}}}
\put(4576,-2986){\makebox(0,0)[lb]{\smash{{\SetFigFont{8}{9.6}{\familydefault}{\mddefault}{\updefault}{\color[rgb]{0,0,0}$l_i$}%
}}}}
\put(3751,-3736){\makebox(0,0)[lb]{\smash{{\SetFigFont{8}{9.6}{\familydefault}{\mddefault}{\updefault}{\color[rgb]{0,0,0}$k_i$}%
}}}}
\end{picture}%

%% file: specvert2.pstex_t
\begin{picture}(0,0)%
\includegraphics{specvert2.pstex}%
\end{picture}%
\setlength{\unitlength}{3947sp}%
\begingroup\makeatletter\ifx\SetFigFont\undefined%
\gdef\SetFigFont#1#2#3#4#5{%
  \reset@font\fontsize{#1}{#2pt}%
  \fontfamily{#3}\fontseries{#4}\fontshape{#5}%
  \selectfont}%
\fi\endgroup%
\begin{picture}(1680,1386)(5986,-4273)
\put(6001,-3736){\makebox(0,0)[lb]{\smash{{\SetFigFont{8}{9.6}{\familydefault}{\mddefault}{\updefault}{\color[rgb]{0,0,0}$k_i$}%
}}}}
\put(7651,-3286){\makebox(0,0)[lb]{\smash{{\SetFigFont{10}{12.0}{\familydefault}{\mddefault}{\updefault}{\color[rgb]{0,0,0}$L$}%
}}}}
\put(6826,-2986){\makebox(0,0)[lb]{\smash{{\SetFigFont{8}{9.6}{\familydefault}{\mddefault}{\updefault}{\color[rgb]{0,0,0}$l_i$}%
}}}}
\end{picture}%

%% file: Simple.pstex_t
\begin{picture}(0,0)%
\includegraphics{Simple.pstex}%
\end{picture}%
\setlength{\unitlength}{3947sp}%
\begingroup\makeatletter\ifx\SetFigFont\undefined%
\gdef\SetFigFont#1#2#3#4#5{%
  \reset@font\fontsize{#1}{#2pt}%
  \fontfamily{#3}\fontseries{#4}\fontshape{#5}%
  \selectfont}%
\fi\endgroup%
\begin{picture}(1824,3095)(8126,-3684)
\put(8326,-2761){\makebox(0,0)[lb]{\smash{{\SetFigFont{10}{12.0}{\familydefault}{\mddefault}{\updefault}{\color[rgb]{0,0,0}$L$}%
}}}}
\put(8551,-736){\makebox(0,0)[lb]{\smash{{\SetFigFont{10}{12.0}{\familydefault}{\mddefault}{\updefault}{\color[rgb]{0,0,0}$d_1-1$}%
}}}}
\put(9301,-736){\makebox(0,0)[lb]{\smash{{\SetFigFont{10}{12.0}{\familydefault}{\mddefault}{\updefault}{\color[rgb]{0,0,0}$d_2$}%
}}}}
\put(9676,-736){\makebox(0,0)[lb]{\smash{{\SetFigFont{10}{12.0}{\familydefault}{\mddefault}{\updefault}{\color[rgb]{0,0,0}$d_s$}%
}}}}
\put(9826,-1861){\makebox(0,0)[lb]{\smash{{\SetFigFont{10}{12.0}{\familydefault}{\mddefault}{\updefault}{\color[rgb]{0,0,0}$d_s$}%
}}}}
\put(9376,-1861){\makebox(0,0)[lb]{\smash{{\SetFigFont{10}{12.0}{\familydefault}{\mddefault}{\updefault}{\color[rgb]{0,0,0}$d_2$}%
}}}}
\end{picture}%

%% file: tildeC.pstex_t
\begin{picture}(0,0)%
\includegraphics{tildeC.pstex}%
\end{picture}%
\setlength{\unitlength}{3947sp}%
\begingroup\makeatletter\ifx\SetFigFont\undefined%
\gdef\SetFigFont#1#2#3#4#5{%
  \reset@font\fontsize{#1}{#2pt}%
  \fontfamily{#3}\fontseries{#4}\fontshape{#5}%
  \selectfont}%
\fi\endgroup%
\begin{picture}(3420,2799)(8464,-2398)
\put(11026, 14){\makebox(0,0)[lb]{\smash{{\SetFigFont{10}{12.0}{\familydefault}{\mddefault}{\updefault}{\color[rgb]{0,0,0}$L$}%
}}}}
\put(10395,-269){\makebox(0,0)[lb]{\smash{{\SetFigFont{10}{12.0}{\familydefault}{\mddefault}{\updefault}{\color[rgb]{0,0,0}$V_7$}%
}}}}
\put(9981,-663){\makebox(0,0)[lb]{\smash{{\SetFigFont{10}{12.0}{\familydefault}{\mddefault}{\updefault}{\color[rgb]{0,0,0}$V_5$}%
}}}}
\put(10154,-483){\makebox(0,0)[lb]{\smash{{\SetFigFont{10}{12.0}{\familydefault}{\mddefault}{\updefault}{\color[rgb]{0,0,0}$V_6$}%
}}}}
\put(9781,-822){\makebox(0,0)[lb]{\smash{{\SetFigFont{10}{12.0}{\familydefault}{\mddefault}{\updefault}{\color[rgb]{0,0,0}$V_4$}%
}}}}
\put(9521,-1076){\makebox(0,0)[lb]{\smash{{\SetFigFont{10}{12.0}{\familydefault}{\mddefault}{\updefault}{\color[rgb]{0,0,0}$V_3$}%
}}}}
\put(8601,-2009){\makebox(0,0)[lb]{\smash{{\SetFigFont{10}{12.0}{\familydefault}{\mddefault}{\updefault}{\color[rgb]{0,0,0}$V_1$}%
}}}}
\put(9074,-1536){\makebox(0,0)[lb]{\smash{{\SetFigFont{10}{12.0}{\familydefault}{\mddefault}{\updefault}{\color[rgb]{0,0,0}$V_2$}%
}}}}
\end{picture}%

%% file: Compute0.pstex_t
\begin{picture}(0,0)%
\includegraphics{Compute0.pstex}%
\end{picture}%
\setlength{\unitlength}{3947sp}%
\begingroup\makeatletter\ifx\SetFigFont\undefined%
\gdef\SetFigFont#1#2#3#4#5{%
  \reset@font\fontsize{#1}{#2pt}%
  \fontfamily{#3}\fontseries{#4}\fontshape{#5}%
  \selectfont}%
\fi\endgroup%
\begin{picture}(1374,967)(7864,-2816)
\put(8326,-2761){\makebox(0,0)[lb]{\smash{{\SetFigFont{10}{12.0}{\familydefault}{\mddefault}{\updefault}{\color[rgb]{0,0,0}$L$}%
}}}}
\end{picture}%

%% file: Compute1a.pstex_t
\begin{picture}(0,0)%
\includegraphics{Compute1a.pstex}%
\end{picture}%
\setlength{\unitlength}{3947sp}%
\begingroup\makeatletter\ifx\SetFigFont\undefined%
\gdef\SetFigFont#1#2#3#4#5{%
  \reset@font\fontsize{#1}{#2pt}%
  \fontfamily{#3}\fontseries{#4}\fontshape{#5}%
  \selectfont}%
\fi\endgroup%
\begin{picture}(1074,967)(8089,-2816)
\put(8326,-2761){\makebox(0,0)[lb]{\smash{{\SetFigFont{10}{12.0}{\familydefault}{\mddefault}{\updefault}{\color[rgb]{0,0,0}$L$}%
}}}}
\end{picture}%

%% file: Compute1b.pstex_t
\begin{picture}(0,0)%
\includegraphics{Compute1b.pstex}%
\end{picture}%
\setlength{\unitlength}{3947sp}%
\begingroup\makeatletter\ifx\SetFigFont\undefined%
\gdef\SetFigFont#1#2#3#4#5{%
  \reset@font\fontsize{#1}{#2pt}%
  \fontfamily{#3}\fontseries{#4}\fontshape{#5}%
  \selectfont}%
\fi\endgroup%
\begin{picture}(1074,967)(8089,-2816)
\put(8326,-2761){\makebox(0,0)[lb]{\smash{{\SetFigFont{10}{12.0}{\familydefault}{\mddefault}{\updefault}{\color[rgb]{0,0,0}$L$}%
}}}}
\end{picture}%

%% file: Compute1c.pstex_t
\begin{picture}(0,0)%
\includegraphics{Compute1c.pstex}%
\end{picture}%
\setlength{\unitlength}{3947sp}%
\begingroup\makeatletter\ifx\SetFigFont\undefined%
\gdef\SetFigFont#1#2#3#4#5{%
  \reset@font\fontsize{#1}{#2pt}%
  \fontfamily{#3}\fontseries{#4}\fontshape{#5}%
  \selectfont}%
\fi\endgroup%
\begin{picture}(924,967)(8089,-2816)
\put(8326,-2761){\makebox(0,0)[lb]{\smash{{\SetFigFont{10}{12.0}{\familydefault}{\mddefault}{\updefault}{\color[rgb]{0,0,0}$L$}%
}}}}
\end{picture}%

%% file: Compute2a.pstex_t
\begin{picture}(0,0)%
\includegraphics{Compute2a.pstex}%
\end{picture}%
\setlength{\unitlength}{3947sp}%
\begingroup\makeatletter\ifx\SetFigFont\undefined%
\gdef\SetFigFont#1#2#3#4#5{%
  \reset@font\fontsize{#1}{#2pt}%
  \fontfamily{#3}\fontseries{#4}\fontshape{#5}%
  \selectfont}%
\fi\endgroup%
\begin{picture}(924,967)(8089,-2816)
\put(8326,-2761){\makebox(0,0)[lb]{\smash{{\SetFigFont{10}{12.0}{\familydefault}{\mddefault}{\updefault}{\color[rgb]{0,0,0}$L$}%
}}}}
\end{picture}%

%% file: Compute2b.pstex_t
\begin{picture}(0,0)%
\includegraphics{Compute2b.pstex}%
\end{picture}%
\setlength{\unitlength}{3947sp}%
\begingroup\makeatletter\ifx\SetFigFont\undefined%
\gdef\SetFigFont#1#2#3#4#5{%
  \reset@font\fontsize{#1}{#2pt}%
  \fontfamily{#3}\fontseries{#4}\fontshape{#5}%
  \selectfont}%
\fi\endgroup%
\begin{picture}(924,1117)(8089,-2816)
\put(8326,-2761){\makebox(0,0)[lb]{\smash{{\SetFigFont{10}{12.0}{\familydefault}{\mddefault}{\updefault}{\color[rgb]{0,0,0}$L$}%
}}}}
\end{picture}%

%% file: Compute2c.pstex_t
\begin{picture}(0,0)%
\includegraphics{Compute2c.pstex}%
\end{picture}%
\setlength{\unitlength}{3947sp}%
\begingroup\makeatletter\ifx\SetFigFont\undefined%
\gdef\SetFigFont#1#2#3#4#5{%
  \reset@font\fontsize{#1}{#2pt}%
  \fontfamily{#3}\fontseries{#4}\fontshape{#5}%
  \selectfont}%
\fi\endgroup%
\begin{picture}(999,1117)(8089,-2816)
\put(8326,-2761){\makebox(0,0)[lb]{\smash{{\SetFigFont{10}{12.0}{\familydefault}{\mddefault}{\updefault}{\color[rgb]{0,0,0}$L$}%
}}}}
\end{picture}%

%% file: Compute2f.pstex_t
\begin{picture}(0,0)%
\includegraphics{Compute2f.pstex}%
\end{picture}%
\setlength{\unitlength}{3947sp}%
\begingroup\makeatletter\ifx\SetFigFont\undefined%
\gdef\SetFigFont#1#2#3#4#5{%
  \reset@font\fontsize{#1}{#2pt}%
  \fontfamily{#3}\fontseries{#4}\fontshape{#5}%
  \selectfont}%
\fi\endgroup%
\begin{picture}(1149,1074)(8089,-2773)
\put(8401,-2686){\makebox(0,0)[lb]{\smash{{\SetFigFont{10}{12.0}{\familydefault}{\mddefault}{\updefault}{\color[rgb]{0,0,0}$L$}%
}}}}
\end{picture}%

%% file: Compute2e.pstex_t
\begin{picture}(0,0)%
\includegraphics{Compute2e.pstex}%
\end{picture}%
\setlength{\unitlength}{3947sp}%
\begingroup\makeatletter\ifx\SetFigFont\undefined%
\gdef\SetFigFont#1#2#3#4#5{%
  \reset@font\fontsize{#1}{#2pt}%
  \fontfamily{#3}\fontseries{#4}\fontshape{#5}%
  \selectfont}%
\fi\endgroup%
\begin{picture}(999,1117)(8089,-2816)
\put(8326,-2761){\makebox(0,0)[lb]{\smash{{\SetFigFont{10}{12.0}{\familydefault}{\mddefault}{\updefault}{\color[rgb]{0,0,0}$L$}%
}}}}
\end{picture}%

%% file: Compute2d.pstex_t
\begin{picture}(0,0)%
\includegraphics{Compute2d.pstex}%
\end{picture}%
\setlength{\unitlength}{3947sp}%
\begingroup\makeatletter\ifx\SetFigFont\undefined%
\gdef\SetFigFont#1#2#3#4#5{%
  \reset@font\fontsize{#1}{#2pt}%
  \fontfamily{#3}\fontseries{#4}\fontshape{#5}%
  \selectfont}%
\fi\endgroup%
\begin{picture}(924,1117)(8089,-2816)
\put(8326,-2761){\makebox(0,0)[lb]{\smash{{\SetFigFont{10}{12.0}{\familydefault}{\mddefault}{\updefault}{\color[rgb]{0,0,0}$L$}%
}}}}
\end{picture}%

%% file: Cycle3.pstex_t
\begin{picture}(0,0)%
\includegraphics{Cycle3.pstex}%
\end{picture}%
\setlength{\unitlength}{2368sp}%
\begingroup\makeatletter\ifx\SetFigFont\undefined%
\gdef\SetFigFont#1#2#3#4#5{%
  \reset@font\fontsize{#1}{#2pt}%
  \fontfamily{#3}\fontseries{#4}\fontshape{#5}%
  \selectfont}%
\fi\endgroup%
\begin{picture}(1795,1599)(9093,-5923)
\end{picture}%

%% file: Cycle4.pstex_t
\begin{picture}(0,0)%
\includegraphics{Cycle4.pstex}%
\end{picture}%
\setlength{\unitlength}{2368sp}%
\begingroup\makeatletter\ifx\SetFigFont\undefined%
\gdef\SetFigFont#1#2#3#4#5{%
  \reset@font\fontsize{#1}{#2pt}%
  \fontfamily{#3}\fontseries{#4}\fontshape{#5}%
  \selectfont}%
\fi\endgroup%
\begin{picture}(1956,1599)(8932,-6073)
\end{picture}%

%% file: Cycle2.pstex_t
\begin{picture}(0,0)%
\includegraphics{Cycle2.pstex}%
\end{picture}%
\setlength{\unitlength}{2368sp}%
\begingroup\makeatletter\ifx\SetFigFont\undefined%
\gdef\SetFigFont#1#2#3#4#5{%
  \reset@font\fontsize{#1}{#2pt}%
  \fontfamily{#3}\fontseries{#4}\fontshape{#5}%
  \selectfont}%
\fi\endgroup%
\begin{picture}(1824,1205)(8914,-5473)
\end{picture}%

%% file: Cycle1.pstex_t
\begin{picture}(0,0)%
\includegraphics{Cycle1.pstex}%
\end{picture}%
\setlength{\unitlength}{2368sp}%
\begingroup\makeatletter\ifx\SetFigFont\undefined%
\gdef\SetFigFont#1#2#3#4#5{%
  \reset@font\fontsize{#1}{#2pt}%
  \fontfamily{#3}\fontseries{#4}\fontshape{#5}%
  \selectfont}%
\fi\endgroup%
\begin{picture}(1424,2090)(9059,-3802)
\put(9421,-1921){\makebox(0,0)[lb]{\smash{{\SetFigFont{5}{6.0}{\familydefault}{\mddefault}{\updefault}{\color[rgb]{0,0,0}$\sigma_3$}%
}}}}
\put(10401,-1851){\makebox(0,0)[lb]{\smash{{\SetFigFont{5}{6.0}{\familydefault}{\mddefault}{\updefault}{\color[rgb]{0,0,0}$\sigma_2$}%
}}}}
\put(9911,-1991){\makebox(0,0)[lb]{\smash{{\SetFigFont{5}{6.0}{\familydefault}{\mddefault}{\updefault}{\color[rgb]{0,0,0}$L$}%
}}}}
\put(9281,-3741){\makebox(0,0)[lb]{\smash{{\SetFigFont{5}{6.0}{\familydefault}{\mddefault}{\updefault}{\color[rgb]{0,0,0}$\sigma_1$}%
}}}}
\end{picture}%

%% file: vertextypeskleinA.pstex_t
\begin{picture}(0,0)%
\includegraphics{vertextypeskleinA.pstex}%
\end{picture}%
\setlength{\unitlength}{3947sp}%
\begingroup\makeatletter\ifx\SetFigFont\undefined%
\gdef\SetFigFont#1#2#3#4#5{%
  \reset@font\fontsize{#1}{#2pt}%
  \fontfamily{#3}\fontseries{#4}\fontshape{#5}%
  \selectfont}%
\fi\endgroup%
\begin{picture}(1687,2000)(8461,-4159)
\put(10133,-2306){\makebox(0,0)[lb]{\smash{{\SetFigFont{9}{10.8}{\familydefault}{\mddefault}{\updefault}{\color[rgb]{0,0,0}$\sigma_2$}%
}}}}
\put(9669,-2439){\makebox(0,0)[lb]{\smash{{\SetFigFont{9}{10.8}{\familydefault}{\mddefault}{\updefault}{\color[rgb]{0,0,0}$L$}%
}}}}
\put(10000,-2903){\makebox(0,0)[lb]{\smash{{\SetFigFont{9}{10.8}{\familydefault}{\mddefault}{\updefault}{\color[rgb]{0,0,0}$(0,-1,0)$}%
}}}}
\put(9868,-3565){\makebox(0,0)[lb]{\smash{{\SetFigFont{9}{10.8}{\familydefault}{\mddefault}{\updefault}{\color[rgb]{0,0,0}$(-1,-1,0)$}%
}}}}
\put(8476,-3168){\makebox(0,0)[lb]{\smash{{\SetFigFont{9}{10.8}{\familydefault}{\mddefault}{\updefault}{\color[rgb]{0,0,0}$(-1,0,0)$}%
}}}}
\put(9072,-4095){\makebox(0,0)[lb]{\smash{{\SetFigFont{9}{10.8}{\familydefault}{\mddefault}{\updefault}{\color[rgb]{0,0,0}$\sigma_1$}%
}}}}
\put(9205,-2372){\makebox(0,0)[lb]{\smash{{\SetFigFont{9}{10.8}{\familydefault}{\mddefault}{\updefault}{\color[rgb]{0,0,0}$\sigma_3$}%
}}}}
\end{picture}%

%% file: vertextypeskleinB.pstex_t
\begin{picture}(0,0)%
\includegraphics{vertextypeskleinB.pstex}%
\end{picture}%
\setlength{\unitlength}{3947sp}%
\begingroup\makeatletter\ifx\SetFigFont\undefined%
\gdef\SetFigFont#1#2#3#4#5{%
  \reset@font\fontsize{#1}{#2pt}%
  \fontfamily{#3}\fontseries{#4}\fontshape{#5}%
  \selectfont}%
\fi\endgroup%
\begin{picture}(1488,2000)(11310,-4159)
\put(12783,-2306){\makebox(0,0)[lb]{\smash{{\SetFigFont{9}{10.8}{\familydefault}{\mddefault}{\updefault}{\color[rgb]{0,0,0}$\sigma_2$}%
}}}}
\put(12320,-2439){\makebox(0,0)[lb]{\smash{{\SetFigFont{9}{10.8}{\familydefault}{\mddefault}{\updefault}{\color[rgb]{0,0,0}$L$}%
}}}}
\put(11723,-4095){\makebox(0,0)[lb]{\smash{{\SetFigFont{9}{10.8}{\familydefault}{\mddefault}{\updefault}{\color[rgb]{0,0,0}$\sigma_1$}%
}}}}
\put(11856,-2372){\makebox(0,0)[lb]{\smash{{\SetFigFont{9}{10.8}{\familydefault}{\mddefault}{\updefault}{\color[rgb]{0,0,0}$\sigma_3$}%
}}}}
\put(12651,-2770){\makebox(0,0)[lb]{\smash{{\SetFigFont{9}{10.8}{\familydefault}{\mddefault}{\updefault}{\color[rgb]{0,0,0}$(0,-1,0)$}%
}}}}
\put(12320,-3764){\makebox(0,0)[lb]{\smash{{\SetFigFont{9}{10.8}{\familydefault}{\mddefault}{\updefault}{\color[rgb]{0,0,0}$(0,0,1)$}%
}}}}
\put(11325,-2969){\makebox(0,0)[lb]{\smash{{\SetFigFont{9}{10.8}{\familydefault}{\mddefault}{\updefault}{\color[rgb]{0,0,0}$(0,1,1)$}%
}}}}
\end{picture}%

%% file: vertextypeskleinC.pstex_t
\begin{picture}(0,0)%
\includegraphics{vertextypeskleinC.pstex}%
\end{picture}%
\setlength{\unitlength}{3947sp}%
\begingroup\makeatletter\ifx\SetFigFont\undefined%
\gdef\SetFigFont#1#2#3#4#5{%
  \reset@font\fontsize{#1}{#2pt}%
  \fontfamily{#3}\fontseries{#4}\fontshape{#5}%
  \selectfont}%
\fi\endgroup%
\begin{picture}(1404,2001)(8744,-6612)
\put(10133,-4758){\makebox(0,0)[lb]{\smash{{\SetFigFont{9}{10.8}{\familydefault}{\mddefault}{\updefault}{\color[rgb]{0,0,0}$\sigma_2$}%
}}}}
\put(9669,-4891){\makebox(0,0)[lb]{\smash{{\SetFigFont{9}{10.8}{\familydefault}{\mddefault}{\updefault}{\color[rgb]{0,0,0}$L$}%
}}}}
\put(9072,-6548){\makebox(0,0)[lb]{\smash{{\SetFigFont{9}{10.8}{\familydefault}{\mddefault}{\updefault}{\color[rgb]{0,0,0}$\sigma_1$}%
}}}}
\put(9205,-4824){\makebox(0,0)[lb]{\smash{{\SetFigFont{9}{10.8}{\familydefault}{\mddefault}{\updefault}{\color[rgb]{0,0,0}$\sigma_3$}%
}}}}
\put(8899,-4995){\makebox(0,0)[lb]{\smash{{\SetFigFont{9}{10.8}{\familydefault}{\mddefault}{\updefault}{\color[rgb]{1,0,0}$n+1$}%
}}}}
\put(9573,-6194){\makebox(0,0)[lb]{\smash{{\SetFigFont{9}{10.8}{\familydefault}{\mddefault}{\updefault}{\color[rgb]{0,0,0}$(1,1,n+2)$}%
}}}}
\put(9883,-5361){\makebox(0,0)[lb]{\smash{{\SetFigFont{9}{10.8}{\familydefault}{\mddefault}{\updefault}{\color[rgb]{0,0,0}$(1,-n,1)$}%
}}}}
\put(8759,-5183){\makebox(0,0)[lb]{\smash{{\SetFigFont{9}{10.8}{\familydefault}{\mddefault}{\updefault}{\color[rgb]{0,0,0}$(0,1,1)$}%
}}}}
\end{picture}%

%% file: vertextypeskleinD.pstex_t
\begin{picture}(0,0)%
\includegraphics{vertextypeskleinD.pstex}%
\end{picture}%
\setlength{\unitlength}{3947sp}%
\begingroup\makeatletter\ifx\SetFigFont\undefined%
\gdef\SetFigFont#1#2#3#4#5{%
  \reset@font\fontsize{#1}{#2pt}%
  \fontfamily{#3}\fontseries{#4}\fontshape{#5}%
  \selectfont}%
\fi\endgroup%
\begin{picture}(2018,2001)(10847,-6612)
\put(12850,-4758){\makebox(0,0)[lb]{\smash{{\SetFigFont{9}{10.8}{\familydefault}{\mddefault}{\updefault}{\color[rgb]{0,0,0}$\sigma_2$}%
}}}}
\put(12386,-4891){\makebox(0,0)[lb]{\smash{{\SetFigFont{9}{10.8}{\familydefault}{\mddefault}{\updefault}{\color[rgb]{0,0,0}$L$}%
}}}}
\put(11789,-6548){\makebox(0,0)[lb]{\smash{{\SetFigFont{9}{10.8}{\familydefault}{\mddefault}{\updefault}{\color[rgb]{0,0,0}$\sigma_1$}%
}}}}
\put(11922,-4824){\makebox(0,0)[lb]{\smash{{\SetFigFont{9}{10.8}{\familydefault}{\mddefault}{\updefault}{\color[rgb]{0,0,0}$\sigma_3$}%
}}}}
\put(10862,-5819){\makebox(0,0)[lb]{\smash{{\SetFigFont{9}{10.8}{\familydefault}{\mddefault}{\updefault}{\color[rgb]{0,0,0}$(-1,0,0)$}%
}}}}
\put(10862,-5620){\makebox(0,0)[lb]{\smash{{\SetFigFont{9}{10.8}{\familydefault}{\mddefault}{\updefault}{\color[rgb]{1,0,0}$k$}%
}}}}
\put(11657,-5023){\makebox(0,0)[lb]{\smash{{\SetFigFont{9}{10.8}{\familydefault}{\mddefault}{\updefault}{\color[rgb]{1,0,0}$l$}%
}}}}
\put(11591,-5222){\makebox(0,0)[lb]{\smash{{\SetFigFont{9}{10.8}{\familydefault}{\mddefault}{\updefault}{\color[rgb]{0,0,0}$(0,1,1)$}%
}}}}
\put(12318,-5836){\makebox(0,0)[lb]{\smash{{\SetFigFont{9}{10.8}{\familydefault}{\mddefault}{\updefault}{\color[rgb]{0,0,0}$(-\alpha+k,-\alpha-l,-\alpha+k)$}%
}}}}
\put(12320,-6150){\makebox(0,0)[lb]{\smash{{\SetFigFont{9}{10.8}{\familydefault}{\mddefault}{\updefault}{\color[rgb]{0,0,0}$(-\alpha,-\alpha,\beta)$}%
}}}}
\end{picture}%

%% file: slopes.pstex_t
\begin{picture}(0,0)%
\includegraphics{slopes.pstex}%
\end{picture}%
\setlength{\unitlength}{3947sp}%
\begingroup\makeatletter\ifx\SetFigFont\undefined%
\gdef\SetFigFont#1#2#3#4#5{%
  \reset@font\fontsize{#1}{#2pt}%
  \fontfamily{#3}\fontseries{#4}\fontshape{#5}%
  \selectfont}%
\fi\endgroup%
\begin{picture}(3402,3324)(8386,-2923)
\put(11026, 14){\makebox(0,0)[lb]{\smash{{\SetFigFont{10}{12.0}{\familydefault}{\mddefault}{\updefault}{\color[rgb]{0,0,0}$L$}%
}}}}
\put(8401,-2761){\makebox(0,0)[lb]{\smash{{\SetFigFont{10}{12.0}{\familydefault}{\mddefault}{\updefault}{\color[rgb]{0,0,0}$\sigma_1$}%
}}}}
\end{picture}%

%% file: NP2a.pstex_t
\begin{picture}(0,0)%
\includegraphics{NP2a.pstex}%
\end{picture}%
\setlength{\unitlength}{3947sp}%
\begingroup\makeatletter\ifx\SetFigFont\undefined%
\gdef\SetFigFont#1#2#3#4#5{%
  \reset@font\fontsize{#1}{#2pt}%
  \fontfamily{#3}\fontseries{#4}\fontshape{#5}%
  \selectfont}%
\fi\endgroup%
\begin{picture}(1602,853)(7786,-3941)
\put(7801,-3511){\makebox(0,0)[lb]{\smash{{\SetFigFont{9}{10.8}{\familydefault}{\mddefault}{\updefault}{\color[rgb]{0,0,0}$1$}%
}}}}
\put(9301,-3511){\makebox(0,0)[lb]{\smash{{\SetFigFont{9}{10.8}{\familydefault}{\mddefault}{\updefault}{\color[rgb]{0,0,0}$1$}%
}}}}
\put(8476,-3211){\makebox(0,0)[lb]{\smash{{\SetFigFont{9}{10.8}{\familydefault}{\mddefault}{\updefault}{\color[rgb]{0,0,0}$b+n$}%
}}}}
\put(8401,-3886){\makebox(0,0)[lb]{\smash{{\SetFigFont{9}{10.8}{\familydefault}{\mddefault}{\updefault}{\color[rgb]{0,0,0}$2n+2+b$}%
}}}}
\end{picture}%

%% file: NP2b.pstex_t
\begin{picture}(0,0)%
\includegraphics{NP2b.pstex}%
\end{picture}%
\setlength{\unitlength}{3947sp}%
\begingroup\makeatletter\ifx\SetFigFont\undefined%
\gdef\SetFigFont#1#2#3#4#5{%
  \reset@font\fontsize{#1}{#2pt}%
  \fontfamily{#3}\fontseries{#4}\fontshape{#5}%
  \selectfont}%
\fi\endgroup%
\begin{picture}(1602,853)(7786,-3941)
\put(7801,-3511){\makebox(0,0)[lb]{\smash{{\SetFigFont{9}{10.8}{\familydefault}{\mddefault}{\updefault}{\color[rgb]{0,0,0}$1$}%
}}}}
\put(8476,-3211){\makebox(0,0)[lb]{\smash{{\SetFigFont{9}{10.8}{\familydefault}{\mddefault}{\updefault}{\color[rgb]{0,0,0}$c$}%
}}}}
\put(8401,-3886){\makebox(0,0)[lb]{\smash{{\SetFigFont{9}{10.8}{\familydefault}{\mddefault}{\updefault}{\color[rgb]{0,0,0}$n+2+c$}%
}}}}
\put(9301,-3511){\makebox(0,0)[lb]{\smash{{\SetFigFont{9}{10.8}{\familydefault}{\mddefault}{\updefault}{\color[rgb]{0,0,0}$1$}%
}}}}
\end{picture}%

%% file: FD.pstex_t
\begin{picture}(0,0)%
\includegraphics{FD.pstex}%
\end{picture}%
\setlength{\unitlength}{3947sp}%
\begingroup\makeatletter\ifx\SetFigFont\undefined%
\gdef\SetFigFont#1#2#3#4#5{%
  \reset@font\fontsize{#1}{#2pt}%
  \fontfamily{#3}\fontseries{#4}\fontshape{#5}%
  \selectfont}%
\fi\endgroup%
\begin{picture}(1827,1837)(10261,-4550)
\put(10610,-4472){\makebox(0,0)[lb]{\smash{{\SetFigFont{9}{10.8}{\familydefault}{\mddefault}{\updefault}{\color[rgb]{0,0,0}$b+n$}%
}}}}
\put(10276,-2836){\makebox(0,0)[lb]{\smash{{\SetFigFont{9}{10.8}{\familydefault}{\mddefault}{\updefault}{\color[rgb]{0,0,0}$b-1$}%
}}}}
\put(10876,-3436){\makebox(0,0)[lb]{\smash{{\SetFigFont{9}{10.8}{\familydefault}{\mddefault}{\updefault}{\color[rgb]{0,0,0}$l_1$}%
}}}}
\put(11251,-3436){\makebox(0,0)[lb]{\smash{{\SetFigFont{9}{10.8}{\familydefault}{\mddefault}{\updefault}{\color[rgb]{0,0,0}$l_{g+1}$}%
}}}}
\put(11551,-4486){\makebox(0,0)[lb]{\smash{{\SetFigFont{9}{10.8}{\familydefault}{\mddefault}{\updefault}{\color[rgb]{0,0,0}$n+1-\sum l_i$}%
}}}}
\end{picture}%

%% file: NP3.pstex_t
\begin{picture}(0,0)%
\includegraphics{NP3.pstex}%
\end{picture}%
\setlength{\unitlength}{3947sp}%
\begingroup\makeatletter\ifx\SetFigFont\undefined%
\gdef\SetFigFont#1#2#3#4#5{%
  \reset@font\fontsize{#1}{#2pt}%
  \fontfamily{#3}\fontseries{#4}\fontshape{#5}%
  \selectfont}%
\fi\endgroup%
\begin{picture}(1824,3095)(8126,-3684)
\put(8176,-1936){\makebox(0,0)[lb]{\smash{{\SetFigFont{10}{12.0}{\familydefault}{\mddefault}{\updefault}{\color[rgb]{0,0,0}$k$}%
}}}}
\put(9301,-736){\makebox(0,0)[lb]{\smash{{\SetFigFont{10}{12.0}{\familydefault}{\mddefault}{\updefault}{\color[rgb]{0,0,0}$d_2$}%
}}}}
\put(9676,-736){\makebox(0,0)[lb]{\smash{{\SetFigFont{10}{12.0}{\familydefault}{\mddefault}{\updefault}{\color[rgb]{0,0,0}$d_s$}%
}}}}
\put(9826,-1861){\makebox(0,0)[lb]{\smash{{\SetFigFont{10}{12.0}{\familydefault}{\mddefault}{\updefault}{\color[rgb]{0,0,0}$d_s$}%
}}}}
\put(9376,-1861){\makebox(0,0)[lb]{\smash{{\SetFigFont{10}{12.0}{\familydefault}{\mddefault}{\updefault}{\color[rgb]{0,0,0}$d_2$}%
}}}}
\put(8326,-2761){\makebox(0,0)[lb]{\smash{{\SetFigFont{10}{12.0}{\familydefault}{\mddefault}{\updefault}{\color[rgb]{0,0,0}$L$}%
}}}}
\put(8851,-736){\makebox(0,0)[lb]{\smash{{\SetFigFont{10}{12.0}{\familydefault}{\mddefault}{\updefault}{\color[rgb]{0,0,0}$l$}%
}}}}
\end{picture}%

%% file: proof.tex
Here we prove Theorems \ref{thm-corres} and \ref{thm:corres 2} 
relating the enumeration of complex
algebraic and tropical curves. Our proof is a mild adaptation of the
proof of {\cite[Theorem 1]{Mi03}} in the framework of tropical
morphisms and their approximation established
in {\cite[Section 6]{Br9}}. 
We begin this section by recalling this latter framework.
Then we reduce  Theorems \ref{thm-corres}
and \ref{thm:corres 2} to a Correspondence Theorem relating the
enumeration of tropical curves in $\RR^2$ to the enumeration of complex curves in
$(\CC^*)^2$ with a fixed Newton fan $\delta$ and having an ordinary
multiple point of maximal multiplicity at a fixed point on a toric
divisor of $\mbox{Tor}(\Pi_{\delta})$. 
We adapt
techniques from {\cite[Section 8]{Mi03}} to our situation, and  are
 eventually 
left to solve some easy local enumerative problems.

\subsection{Phase-tropical geometry}\label{sec:phase trop}
Here we recall definitions and results from 
{\cite[Section 6]{Br9}} 
that we  need later. We start to define the \emph{phase} of
a point, and of a tropical morphism. Roughly speaking, the phase of a
tropical variety  $Y$ is a choice, in a compatible way,
 of an algebraic variety $\mathcal
Y_p$ of dimension $\dim Y$ for each point $p\in Y$.

The notion of a phase-tropical limit is based on
the degeneration of the standard complex
structure on $(\CC^*)^n$ via the following self-diffeomorphism of  $(\CC^*)^n$:
$$\begin{array}{cccc}
H_t:&(\CC^*)^n&\longrightarrow& (\CC^*)^n
\\ &(z_i) &\longmapsto & (|z_i|^{\frac{1}{\log t}}\frac{z_i}{|z_i|})
 \end{array}.$$
We also define the two following maps:
$$\begin{array}{cccc}
\Log: & \ctorn &\longrightarrow & \RR^n
\\ & (z_i)&\longmapsto & (\log(|z_i|))\end{array}
 \quad \mbox{and}\quad
\begin{array}{cccc}
\Arg: & \ctorn &\longrightarrow & \sonen
\\ & (z_i)&\longmapsto & (\arg(z_i))\end{array}.$$

As usual, all definitions are particularly easy in the case of points.
\begin{definition}[{\cite[Definitions 6.1 and 6.4]{Br9}}]
Let $p\in\RR^n$ be a point.
A \emph{phase} of  $p$
 is the choice of a point $\phi(p)\in(S^1)^n$.

\vspace{1ex}
Let $(p_{t_j})$ be a sequence of points in  $\ctorn$ such that
$\lim\limits_{j\to+\infty} H_{\tj}(p_{\tj})$ exists as a point in
$\ctorn$.
The \emph{phase-tropical limit} of the sequence  $(p_{t_j})$
is defined as the point $p=\Log(\lim\limits_{j\to+\infty}
H_{\tj}(p_{\tj}))$
 enhanced with the phase $\Arg(\lim\limits_{j\to+\infty} H_{\tj}(p_{\tj}))$.
\end{definition}

\vspace{2ex}
Next, we define \emph{phase-tropical morphisms} and the \emph{phase-tropical limit} of a
sequence of algebraic maps from Riemann surfaces.
A \emph{pluriharmonic map $\Phi:S\to \sonen $} is a map from a
punctured Riemann surface $S$ which is the composition of an algebraic
map $\Phi_0:S\to \ctorn $
with the map $\Arg:\ctorn\to\sonen$. Such an algebraic map $\Phi_0$ is
called an \emph{algebraic lift} of $\Phi$. Clearly, two algebraic
lifts of $\Phi$ 
differ by a multiplicative translation in $\ctorn$ 
by a vector in $(\RR_{>0})^n$.

Given a pluriharmonic map $\Phi : S\to \sonen$ we can naturally
associate
 a map $\Phi^\varepsilon: S^1\to \sonen$ for
each of the punctures  $\varepsilon$ of $S$.
 Let us denote by $\overline S$ the compact Riemann
surface obtained from $S$ by performing a real blow-up at each
puncture of $S$. That is to say we replace each puncture $\varepsilon$
of $S$ with a \emph{boundary circle} $b_\varepsilon$ of
length $2\pi$, oriented as a boundary component of $\overline S$, the
metric on $b_\varepsilon$ being given by the conformal structure of $S$ at
$\varepsilon$. (see {\cite[Section 6.2]{MikhalkinOkounkov}} or
{\cite[Section 6.1]{Br9}}).  
The map $\Phi^\varepsilon: b_\varepsilon\to \sonen$ is defined
 as the limit of the map $\Phi_{|l}$ where $l$ is a small
 loop around $\varepsilon$ converging to $\varepsilon$.

The map $\Phi$ is proper at $\varepsilon$ if and only if
 $\Phi$ cannot be extended at $\varepsilon$ to a pluriharmonic map,
 i.e.\ an algebraic lift of $\Phi$ does not send any neighborhood of
 $\varepsilon$ into a compact subset of $\ctorn$.
In
 this case   the  map $\Phi^\varepsilon$ is a covering of some degree $
w\ge 1$
 of a geodesic on the flat torus $\sonen=(\RR/2\pi )^n$ 
(see {\cite[Section 6.2]{MikhalkinOkounkov}} or
{\cite[Section 6.1]{Br9}}). Hence it
is a
dilation of factor
$w$ if $\Phi^\varepsilon(b_\varepsilon)$ is equipped with the metric,
of total length $2\pi$,
induced by the natural flat metric on $(\RR/2\pi )^n$.

If $\Phi$ is not proper in a neighborhood of a puncture $\varepsilon$,
then it has a removable singularity, and $\Phi^\varepsilon$ maps the
whole circle $b_\varepsilon$ to a point. 

Recall that a stable  Riemann surface $S$
 is a, maybe reducible, nodal complex algebraic curve such that any of its
irreducible component is a punctured Riemann surface; the total number of nodes
 and punctures on a component $S_0$ of $S$ is at least 3 if
 $S_0$ is rational, and at least 1 if $S_0$ is elliptic. 
If  $\mathcal I$ is
the graph of intersection of the irreducible components of $S$ (i.e.\ the dual graph),
 the genus of $S$ is equal to $b_1(\mathcal I)+\sum_{S_0}g(S_0)$.
We denote by  $S^\circ$ the Riemann surface obtained from $S$
by removing all nodes.

In the following definition, we identify the group $H_1(\sonen,\ZZ)$ with
$\ZZ^n$ via the map $\Arg$ (recall that since $\CC^*$ is canonically
oriented, the group $H_1(\CC^*,\ZZ)$ is
canonically identified with $\ZZ$, and therefore $H_1(\ctorn,\ZZ)$ is
canonically identified with $\ZZ^n$).

Finally, we say that a tropical morphism $h:C\to\RR^n$
 is \emph{minimal} if $v(V,e)\ne 0$ for any edge $e$ of $C$ (see Definition \ref{def:trop morph}). The next
 definition only deals with minimal tropical morphisms. For a
 definition in a more general situation we refer to the forthcoming
 paper \cite{Mik08}.
\begin{definition}[{\cite[Definition 6.5]{Br9}}]\label{phase-curve}
Let $h:C\to\RR^n$ be a minimal tropical morphism.
A \emph{phase $\phi$ of $h$}  consists of the following the data

\begin{itemize}
\item for  each vertex
$V\in\Ve^0(C)$, a proper pluriharmonic map 
$$\Phi_{V}:S_{V}\to\sonen$$
where $S_V$ is a stable 
Riemann surface
of genus $g_V$
with $k$ punctures, equipped with  a one-to-one correspondence
$\varepsilon\leftrightarrow e$ between
the punctures of $S_V$ and the edges of $C$ adjacent to $V$, such that for each 
edge $e$, the homology class 
$[\Phi_V^e(b_e)]\in
  H_1(\sonen, \ZZ)=\ZZ^n$ satisfies
$$[\Phi_V^e(b_e)]=v(V,e)\in\ZZ^n;$$
\item for each edge $e\in\Ed^0(C)$ adjacent to
 $V,V'\in\Ve^0(C)$,
an orientation-reversing isometry 
\begin{equation}\label{rho-e}
\rho_e:b_e^V\to b^{V'}_e
\end{equation}
between the two boundary circles
of the punctures corresponding to $e$, such that
$\Phi^e_V=\Phi^e_{V'}\circ\rho_e$;

\item for  each vertex
$V\in\Ve^0(C)$ and each node $\kappa$ of $S_V$, 
an orientation-reversing isometry
\begin{equation}\label{rho-k}
\rho_\kappa:b_{\varepsilon'}\to b_{\varepsilon''}
\end{equation}
where $\varepsilon'$ and $\varepsilon''$ are the two punctures of 
$S_V^\circ$ corresponding to the node $\kappa$.
\end{itemize}
\end{definition}

Note that in this latter case, both boundary circles
$b_{\varepsilon'}$ and $b_{\varepsilon''}$
are mapped to the same point in $\sonen$, in particular one has 
 $\Phi^{\varepsilon'}_V=\Phi^{\varepsilon''}_V\circ\rho_\kappa.$
We denote by $(h,\phi)$ a \emph{phase-tropical morphism}, i.e. a
tropical morphism $h$ equipped with a phase $\phi$. 

\begin{remark}
There exist slight formal differences between
Definition \ref{phase-curve} and 
{\cite[Definition 6.5]{Br9}}. However these differences  only come
from minor differences in the presentation we chose in this
paper compared to the one chosen
in \cite{Br9}.
In this latter tropical curves
are not allowed to contain degenerate edges, and  pluriharmonic
maps $\Phi_V$ might be non-proper (in the language of \cite{Br9} $\Phi_V$
might have \emph{essential} boundary circles). In the present paper,
essential
boundary circles are replaced by degenerate edges, and the maps
$\Phi_V$ are now required to be proper. 
\end{remark}

By definition, a pluriharmonic map $\Phi_V$ is associated to each
vertex $V$ of a
phase-tropical morphism $(h:C\to\RR^n,\phi)$.
If  $S_V$ is smooth then we set $\hat S_V=S_V$,   and  $\hat \Phi_V=\Phi_V$.
If $S_V$ is nodal, then we denote by $\hat S_V$ the topological oriented surface obtained by replacing each node $\kappa$ with the corresponding boundary circle, i.e. either side of the isometry (\ref{rho-k}). By construction we have a natural continuous map $\hat S_V\to S_V $ contracting each boundary circle to the corresponding node of $S_V$. Furthermore the pluriharmonic map $\Phi_V$ naturally induces a continuous map $\hat \Phi_V: \hat S_V\to \sonen$.

If $W_c$ is a degenerate component of $C$, we denote by  $\hat S_{W_c}$  
the topological oriented surface obtained by gluing all surfaces $\hat S_V$,  with $V$ ranging over all vertices of $W_c$, along all boundary circles 
corresponding to edges of $W_c$ using isometries (\ref{rho-e}). 
We define $\widetilde S_{W_c}$ as the surface
 $\hat S_{W_c}$ with these boundary circles removed.
By construction 
the surface $\widetilde S_{W_c}$
 is the disjoint union of 
the surfaces  $\hat S_{V}$ over all vertices $V$ of $W_c$. 
Furthermore all pluriharmonic maps $\hat \Phi_V$ with $V$ a vertex of $W_c$ 
naturally induce a continuous map $\hat \Phi_{W_c}: \hat S_{W_c}\to \sonen$.

For
each non-degenerate edge $e$ of $C$ adjacent to $V$, we write $v(e,V)=(v_1,\ldots,v_n)$, 
and we  associate to $e$ 
the Riemann surface $S_e=\CC^*$ and the pluriharmonic map
$$\begin{array}{cccc}
\Phi_e :& \CC^* &\longrightarrow &\sonen
\\ & re^{i\theta}&\longmapsto & (v_1\theta,\ldots, v_n\theta)
\end{array}.$$
An algebraic lift of $\Phi_e$ is given by the map
$z\to (z^{v_1},\ldots, z^{v_n})$, 
which might be thought as the complexification of the map 
$\Phi_V^e:b_e^V\to (S^1)^n$:
 since 
$b_e^V$ 
(resp. each coordinate circle of $(S^1)^n$)
is
oriented, the tangent space of $b_e^V$ 
(resp.  $(S^1)^n$) can  naturally be identify with $\CC^* $
(resp. $\ctorn$). 
Note that $\Phi_e(\CC^*)=\Phi_V^e(b^V_e)$.

A connected open subset of $C$ is said to be \emph{admissible} if 
any  degenerate component of $C$ is either disjoint or contained in $W$, and if
$W$ contains two vertices of $C$ then these two vertices belong to the same
 degenerate component.
 In particular $W$ contains at most one degenerate component, and if not, at most one vertex.
Given $W$ such an admissible connected open subset of $C$,
we 
denote by $S_W=\widetilde S_W=\hat S_V$ and $\Phi_W=\hat \Phi_V$ if $V$ is the unique vertex of $C$ in $W$, by $S_W=\widetilde S_W=S_e$ and $\Phi_W=\Phi_e$
if $W$ is contained in the non-degenerate edge $e$ of $C$, and
by $S_W=\hat S_{W_c}$, $\widetilde S_W=\widetilde S_V$, 
and $\Phi_W=\hat \Phi_{W_c}$ if $W$ contains the degenerate component $W_c$.

Let $h:C\to\RR^n$ be a minimal tropical morphism.
A convex open subset $U\subset\RR^n$ is called \emph{$h$-admissible} if
$h(C)\cap U$ is connected, if $U$ contains
 at most one point $p$ such that $h^{-1}(p)$ contains a vertex of $C$,
 and if no vertex of $C$ is mapped to the boundary of $U$.

Finally, we call a \emph{positive multiplication} in $\ctorn$ a map of the form
$\tau(z_1,\ldots,z_n)= (a_1z_1,\ldots,a_nz_n)$ with $a_1,\ldots,a_n>0$.
In next definition we use the  standard flat metric on $\sonen$.
\begin{definition}[{\cite[Definition 6.5]{Br9}}]\label{wtroplimit}
Let 
$f_{\tj}:S_{\tj}\to\ctorn$ be a sequence of algebraic
 maps from  punctured Riemann surfaces $S_{\tj}$, with $t_j\to+\infty$
 when $j\to+\infty$.
We say that a phase-tropical morphism 
$(h:C\to \RR^n,\phi)$  is the
\emph{phase-tropical limit} of $f_{\tj}$
if for any choice of $h$-admissible open set $U\subset\RR^n$  and
all sufficiently large $\tj$ there is a 1-1 correspondence between
connected components
$W_{\tj}$ of $f_{\tj}^{-1}(\Log_{\tj}^{-1}(U))$ 
and connected components $W$ of $h^{-1}(U)$
with the following properties
of the corresponding components.
\begin{itemize}
\item
There exists an open embedding $\Xi^W_{\tj}:W_{\tj}\to S_W$ and, for any connected component $S$ of $\widetilde S_W$, an algebraic lift $\Psi_{S}:S\to \ctorn $  of $\Phi_W|_S$ and a sequence $(\tau_{t_j})$ of positive translations 
in $\ctorn$ such that for 
 any $z\in S$ 
$$\lim\limits_{t_j\to +\infty}\tau_{t_j}\circ f_{\tj}\circ (\Xi^W_{\tj})^{-1}(z) = \Psi_S(z).$$ 
In particular, we require that $z\in\Xi^W_{\tj}(W_{\tj})$ for large $\tj$.
\item
For any edge $e$ connecting vertices $V$ and $V'$, any point $z\in b^V_e$, any $\eta>0$
and a sufficiently large $\tj$
there exist 
\begin{enumerate}
\item a point $z_\eta\in S_V$ and a point $z'_\eta\in S_{V'}$;
\item a path $\gamma_\eta\subset\bar S_V$ connecting $z_\eta$ and $z$
and a path $\gamma'_\eta\subset\bar S_{V'}$ connecting $z'_\eta$ and $z'=\rho_e(z)$ (see Equation \eqref{rho-e} in Definition \ref{phase-curve})
such that the diameter of $\Phi_V(\gamma_\eta)\subset \sonen$ and that of
$\Phi_{V'}(\gamma'_\eta)\subset\sonen$ are less than $\eta$;
\item 
a path $\gamma_{\tj}\subset S_{\tj}$ connecting $(\Xi^v_{\tj})^{-1}(z_\eta)$ and $(\Xi^{v'}_{\tj})^{-1}(z'_\eta)$
such that the diameter of $\Arg(f_{\tj}(\gamma_{\tj}))\subset \sonen$ is less than $\eta$.
\end{enumerate}

\end{itemize}

\end{definition}

It follows from Definition \ref{wtroplimit} that if
$(h:C\to \RR^n,\phi)$  is the 
phase-tropical limit of $f_{\tj}:S_{\tj}\to(\CC^*)^n$, then $h(C)$ is
the limit, in the Hausdorff metric on compact sets of $\RR^n$, of
$\Log_{t_j}(f_{\tj}(S_{\tj}))$. 
In the case of curves mapped to $\ctor$, each curve $f_{\tj}(S_{\tj})$
has a polynomial equation $P_{t_j}(x,y)=0$, and
 the 
equation of the tropical limit $h(C)\subset \RR^2$ 
can be deduced from the sequence of
polynomials
$P_{t_j}(x,y)$.
\begin{proposition}[{\cite[Theorem 6.4]{Mik12}}]\label{lem:NP}
Let $P_{t_j}(x,y)=\sum_{i,j}a_{i,j,\tj}x^iy^j$ be a sequence of complex
polynomials, such that one of the sequence $a_{i_0,j_0,\tj}$ is the
constant sequence equal to 1, and
 the absolute value of any 
coefficients $a_{i,j,\tj}$ satisfies
$|a_{i,j,\tj}|=O_{j\to +\infty}(1) $. Suppose in addition that
$\Log_{t_j}(\{P_{t_j}(x,y)=0\})$ converges, in the Hausdorff metric on
compact sets of $\RR^2$, to a tropical curve $C\subset \RR^2$. Then $C$
is given by the tropical polynomial
$$\tg \sum_{i,j}\lambda_{i,j}x^iy^j \td $$
where $\lambda_{i,j}=\inf\{\lambda\ | \ a_{i,j,\tj}=o_{j\to +\infty}(\tj^\lambda) \}$.
\end{proposition}

\begin{example}
We consider the following family of algebraic maps
$$\begin{array}{cccc}
f_t: & \CC^*\setminus\{1,\frac{t}{t-1}\} &\longrightarrow & \ctor
\\ & z &\longmapsto & \left(\frac{z}{1-z}, \frac{t}{(t-1)z -t}  \right)
\end{array}. $$
The image of the map $f_t$ is given by the polynomial
$P_t(x,y)=1+x+y+t^{-1}xy$ in $\ctor$, 
so it converges to the embedded tropical curve $C'$ given by the
polynomial $\tg 0+ 
x + y + (-1)xy \td$. 
This tropical curve has two vertices $(0,0)$ and
$(1,1)$, one bounded edge of direction $(1,1)$, and four ends of
direction $(-1,0),(0,-1),(1,0),$ and $(0,1)$ (see Figure \ref{fig:ex
conv}). Now one computes easily that the family $(f_t)$ converges
tropically to the phase-tropical morphism $(h:C\to \RR^2,\phi)$ where
\begin{figure}[ht]
\centering
\begin{tabular}{c}
\input{Ex_mor3.pstex_t}
\end{tabular}
\caption{A tropical limit.}
\label{fig:ex conv}
\end{figure}
\begin{itemize}
\item $h$ is the unique parameterization of $C'$ from a rational tropical
curve 
$C$  with four ends and with all weights equal to $1$;
\item if $V_0$ is the vertex of $C$ mapped to $(0,0)$, then 
$S_{V_0}$ is the curve in $\ctor$ given by the equation
 $1+x+y=P_{+\infty}(x,y)$
 and 
$\Phi_{V_0}$ is the restriction of the argument map to $S_{V_0}$;

\item if $V_1$ is the vertex of $C$ mapped to $(1,1)$, then 
$S_{V_1}$ is the curve in $\ctor$ given by the equation 
$x+y+xy=\lim_{t\to+\infty} \frac{1}{t}P_t(tx,ty) $ and 
$\Phi_{V_1}$ is the restriction of the argument map to $S_{V_1}$.
\end{itemize}
Note that there is a unique possibility for the map $\rho_e$
corresponding to the bounded edge $e$ of $C$.
\end{example}

\begin{example}
Proposition  \ref{lem:NP} implies that a phase tropical morphism
$(h,\phi)$ which is the phase-tropical limit of
 a sequence of algebraic maps $f_{\tj}:S_{\tj}\to\ctorn$ with a fixed
 Newton fan $\delta$ might have a Newton fan 
different from 
 $\delta$. 
For example, the the image of map $f_t(z)=(-1-e^{-t}z,z)$ in $(\CC^*)^2$ satisfies the
equation
$1+x+e^{-t}y=0$, and the family $(f_t)$ 
has a tropical limit with Newton polygon the segment
$[(0,0),(1,0)]\subset\RR^2$. 
\end{example}

We define the \emph{Euler characteristic of a phase-tropical 
morphism $(h,C\to\RR^n,\phi)$} by
$$\chi(h)=\sum_{V\in\Ve^0(C)}\chi(S_V).$$
In particular  one has $2\chi_{\mbox{trop}}(C)\le \chi(h)
+\#\Ed^\infty(C)$, with equality if and only if all surfaces $S_V$ are smooth. 
The next lemma is an immediate consequences of Definition \ref{wtroplimit}.
\begin{lemma}\label{lem:irr}
Let $f_{\tj}:S_{\tj}\to\ctorn$ be a sequence of non-constant algebraic maps which
converges, as a tropical limit, to a phase-tropical morphism
$(h:C\to\RR^n,\phi)$. Suppose that the Riemann surfaces $S_{t_j}$ are
connected and of constant Euler characteristic $\chi$ 
for  $t_j$ large enough. Then $\chi(h)\ge \chi$, and if equality
holds then $C$ is connected.
\end{lemma}

The following compactness result is fundamental
for proving Correspondence Theorems in the phase-tropical  framework.
 \begin{proposition}[{\cite[Proposition 6.8]{Br9}}]\label{trop-compact}
Let $f_t:S_t\to\ctorn$ 
be  a family of algebraic maps with a fixed genus and Newton fan, 
defined for all sufficiently large positive parameter $t>>1$.
Then there exists a phase-tropical morphism $(f,\phi)$ 
and a sequence $\tj\to +\infty$ such that $f_{\tj}$ converges to
$(f,\phi)$
in the sense of Definition \ref{wtroplimit}.
\end{proposition}

Now we prove that the tropical limit of a sequence of 
algebraic maps to the surface $\mathcal X$ in $(\CC^*)^3$
with equation $x+y+z=0$ is a tropical morphism to $X$ in the sense of Definition \ref{def:morphX}.
It is well known that $X$ is
the limit, for the Hausdorff metric on compact sets of $\RR^n$,
of $\Log_t(\mathcal X)$ when $t\to +\infty$ (see for example \cite{Mik8}). Note that $\mathcal X$ is
invariant under any multiplicative translation by  an element of
$\ctorn$ of the form $(\lambda,\ldots, \lambda)$.

The next proposition is a
particular case of a result proved by Mikhalkin and the first author
in the forthcoming
paper \cite{Br12}. Since the latter is not available yet, we provide a
proof for the sake of completeness.

We use the following notation in the proof of
Proposition \ref{prop:RH}.
Let us consider a sequence of algebraic maps $f_j: S_{j}\to(\CC^*)^3$
from a family of punctured Riemann surfaces
such that the sequence of compactified maps
$\overline f_{j}:\overline S_{j}\to\CC P^3$ converges 
 to some map $\overline f_{\infty}:\overline S_\infty\to\CC P^3$. Note that
the Riemann surface $\overline S_\infty$ might be reducible, and some of its
connected components might be mapped to the toric boundary divisors of
$\CC P^3$. 
We denote by $f_\infty: S_\infty\to (\CC^*)^3$ the 
 restriction of $\overline f_\infty$ to $(\CC^*)^3$.
We emphasize that the Riemann surface $S_\infty$ might be disconnected.

\begin{proposition}[\cite{Br12}]\label{prop:RH}
Let $f_{\tj}:S_{\tj}\to(\CC^*)^3$ be a sequence of algebraic maps which
converges, as a tropical limit, to a phase-tropical morphism
$(h:C\to\RR^3,\phi)$. Suppose in addition that 
 $f_{\tj}(S_{\tj})\subset \mathcal X$  for 
$\tj$ large enough.
Then $h$ is
a tropical morphism to $X$.
\end{proposition}
\begin{proof}
Since $h(C)$ and $X$ are respectively 
the limit, in the Hausdorff metric on compact sets of $\RR^3$, of
$\Log_{t_j}(f_{\tj}(S_{\tj}))$ and $\Log_{t_j}(\mathcal X)$, we have $h(C)\subset X$.

Let $V$ be a vertex of $C$.
If $\tau=(\tau_{t_j})$ is a sequence of positive translation in $(\CC^*)^3$
we denote by $f^\tau_{\tj}$ the sequence of maps
$\tau_{t_j}\circ
 f_{\tj}$.
By
Definition \ref{wtroplimit} there exists
a sequence $\tau_{t_j}(x,y,z)=(a_{t_j,1}x,a_{t_j,2}y,a_{t_j,3}z)$ of
positive translations  
in $(\CC^*)^3$ for which
there exists a connected Riemann surface $S_V\subset
S_\infty$  such that $f^\tau_{\infty \ |S_V}: S_V\to\ctorn$ is an algebraic
lift of the phase $\Phi_V:S_V\to (S^1)^3$ of
$h:C\to\RR^3$ at the vertex $V$.
Note that $\lim \Log_{t_j} (a_{t_j,1},a_{t_j,2},a_{t_j,3})=h(V)$.

There are only three possibilities, up to a positive translation, for
the limit $\mathcal X^\tau_\infty$ of $\tau_{t_j}(\mathcal X)$: it has equation either $x+z=0$,
$x+y=0$, $y+z=0$, or $x+y+z=0$. Moreover this latter arises if and
only if $f(V)=(u,u,u)$ and the sequence
$(t_j^{-u}a_{t_j,1},t_j^{-u}a_{t_j,2},t_j^{-u}a_{t_j,3})$ is contained
in a compact set of $(\CC^*)^3$.
Moreover, since $f_{\tj}(S_{\tj})\subset \mathcal X$  for all $\tj$
large enough, we have  $f^\tau_{\infty \ |S_V}(S_V)\subset \mathcal
X^\tau_\infty$.
In  particular this implies that $h:C\to X$ is a tropical premorphism, 
i.e.\  $h$ satisfies the second condition of Definition \ref{def:premorphism}.

Let $V$ be a vertex of $C$ with $d_V>0$. In particular, $h(V)=(u,u,u)$
with $u\in\RR$ and the corresponding sequence
$(t_j^{-u}a_{t_j,1},t_j^{-u}a_{t_j,2},t_j^{-u}a_{t_j,3})$ is contained
in a compact set of $(\CC^*)^3$. After composing the map $f_{\tj}$ with the
multiplicative translation by $(t_j^{-u}a_{t_j,1},t_j^{-u}a_{t_j,2},t_j^{-u}a_{t_j,3})$, and 
the map $h$ with the additive translation by $(-u,\ldots,-u)$, we may
suppose that $u=0$ and $(a_{t_j,1},a_{t_j,2},a_{t_j,3})=(1,1,1)$.
From what we said above, there exists an algebraic lift 
$f_{\infty \ |S_V}: S_V\to\ctorn$ of the phase $\Phi_V:S_V\to (S^1)^3$
such that $f_{\infty}( S_V)\subset \mathcal X$.

Let us consider the plane $\mathcal P\simeq \ctor\subset (\CC^*)^3$
with equation $z=1$. Then the curve $S_0=\mathcal P\cap \mathcal X$ is
the Riemann sphere punctured in three points, and the map
$(x,y,z)\mapsto ((\frac{x}{z},\frac{y}{z},1),z)$ provides an algebraic
isomorphism between $\mathcal X$ and $S_0\times \CC^*$. 
The restriction on
$f_{\infty}(S_V)$ of 
 the projection to the first factor $\mathcal X\to S_0$ gives
a  holomorphic map $\pi:S_V\to S_0$, which
can be extended to all  punctures of $S_V$ 
 corresponding to  an edge of $C$
entirely mapped to $L$. In this way we obtain a proper holomorphic map 
$\overline \pi:\overline S_V\to S_0$, which has 
 degree
$d_V$ by {\cite[Theorem 1.1]{ST07}}. Smoothing  nodes of $\overline
S_V$ if there are any, we may even further suppose that 
$\overline S_V$ is a non-singular Riemann surface of genus $g_V$ with
$k_V$ punctures, where $k_V$ is the number of edges adjacent to $V$
which are not entirely mapped to $L$.
Hence
it follows from the Riemann-Hurwitz formula that
$$\chi(S_V)= d_{V}\chi(S_0) -\iota$$
where $\iota\ge 0$ is the sum of ramification indices of
$\overline\pi$ over all points of $\overline S_V$. Since we have 
$\chi(\overline S_V)=2-2g_V -k_V$ and $\chi(S_0)=-1$, we obtain
$$k_V-d_V-2+2g_V\ge 0.$$
which means exactly that $f$ is a tropical morphism to $X$.
\end{proof}

\subsection{The proof of Theorem \ref{thm-corres}}
Let us first recast  notations from previous sections.
We denote by $X$ (resp. $\mathcal X$) the tropical hypersurface in
$\RR^3$ (resp. in $(\CC^*)^3$)
defined by the tropical polynomial $\tg x+y+z \td$ (resp. the equation
 $x+y+z=0$).  The tropical surface $X$ is made of three half-planes
 $\sigma_1=\{x=y\geq z\}$, $\sigma_2=\{x=z\geq
y\}$ and $\sigma_3=\{y=z\geq x\}$ meeting along the line
$L=\RR(1,1,1)$.
We have fixed the following Newton fan
$$\Delta= \{(1,-n,1)^a, (0,1,1)^{an+b}, (-1,0,0)^a, (0,-1,0)^b,
(0,0,-1)^{a(n+1)+b} \}$$
as well as an integer $\chi\in\ZZ$ and a generic configuration
$\omega$ of
$\#\Delta_2+\#\Delta_3 -\chi$ points in 
$\sigma_1\cup \sigma_2 $, where 
 $\Delta_i\subset \Delta$ consists of
elements of $\Delta$ contained in $\sigma_i$.

Let us fix the following  two additional Newton
fans
$$\delta_z=\{(1,-n)^a, (0,1)^{an+b}, (-1,0)^a, (0,-1)^b\},$$ 
and
$$\delta_x=\{(-n,1)^a, (1,1)^{an+b}, (-1,0)^b, (0,-1)^{a(n+1)+b}\}.$$ 
Note that $\delta_z$ is the image of $\delta_0$ (see
Theorem \ref{thm-corres}) under the map $(x,y)\mapsto (x,-y)$, in
particular
$N_{2\chi}^{irr}(\delta_0)=N_{2\chi}^{irr}(\delta_z) $.
We define the projections
$$\begin{array}{cccc}
\pi_z: & (\CC^*)^3 &\longrightarrow &\ctor
\\ & (x,y,z)&\longmapsto &(x,y)
\end{array}\ \mbox{and}\ 
\begin{array}{cccc}
\pi_x: & (\CC^*)^3 &\longrightarrow &\ctor
\\ & (x,y,z)&\longmapsto &(y,z)
\end{array}. $$
Note that $\delta_z=\pi_z(\Delta)$, 
$\delta_x=\pi_x(\Delta)$, and $\#\delta_z =\#\Delta_2 +\#\Delta_3$.
Furthermore, we choose  a configuration $\omega^\CC$ of 
$\#\delta_z-\chi$ points in $\X$. We denote by
 $\S(\Delta,\omega^\CC)$ the set of all irreducible 
algebraic curves in $\X$ of Euler characteristic $2\chi-\#\Delta$, with
Newton fan $\Delta$, and passing through $\omega^\CC$. We define
$N_{2\chi}^{irr}(\Delta)=\#\S(\Delta,\omega^\CC)$.

The strategy to prove the Correspondence Theorem is as follows: we first ``put the algebraic enumerative problem into three-space'', i.e.\ we prove in Lemma \ref{lem:hirzebruch} that $N_{2\chi}^{irr}(\Delta)=N_{2\chi}^{irr}(\delta_z)$. Note that this 
implies in particular that $N_{2\chi}^{irr}(\Delta)$  does not
depend on $\omega$.
As usual, we have to prove two statements for a Correspondence Theorem: first we have to show that the algebraic curves considered in our enumerative problem (i.e.\ now the curves in $\S(\Delta,\omega^\CC)$) degenerate to the curves in the tropical enumerative problem, i.e.\ to elements in $\TT\mathcal{S} (\omega)$. Second, we have to show that the number of algebraic curves degenerating to a fixed tropical curve equals the tropical multiplicity. 
The first part is Lemma \ref{lem:S to TS}.
For the second part, we use the projection $\pi_x$ to ``put the algebraic enumerative problem back into the plane'', but differently. Our intention to use this different projections is that then we are able to apply known techniques for Correspondence Theorems for plane curves.

Contrary to the situation where we project with $\pi_z$ however, we do not obtain an enumerative problem that involves only simple point conditions. Instead, we obtain curves with a multiple point (see Lemma \ref{lem:multiple}).
We then relate this new plane algebraic enumerative problem to our
tropical curves in $X$
and 
their projections 
by
$\pi_x$. Compared to the existing Correspondence Theorems for plane curves, our situation differs since we have vertices to which several edges of the same direction are adjacent (the projections of vertices in $L$). We compute the number of algebraic preimages for these vertices locally in Lemma \ref{lem:local comp}.

\begin{lemma}\label{lem:hirzebruch}
For a generic configuration of points $\omega^\CC$, the set
$\S_{2\chi}(\Delta,\omega^\CC)$ is finite and
$$N_{2\chi}^{irr}(\Delta)=N_{2\chi}^{irr}(\delta_z). $$
\end{lemma}
\begin{proof}
There clearly exists a bijection between curves in $\S(\Delta,\omega^\CC)$
 and irreducible complex
algebraic
curves in
$(\CC^*)^2$  of Euler characteristic $2\chi-\#\delta_z$, with
Newton fan $\delta_z=\pi_z(\Delta)$, and passing through $\pi_z(\omega^\CC)$. 
\end{proof}

Recall that we chose a tropically generic configuration $\omega$ of 
$\#\delta_z-\chi$ points 
in $\sigma_1\cup\sigma_2\subset X$. By Proposition \ref{prop-curves}
this implies in particular that for any point $p\in\omega$ and any
morphism $h:C\to X$ in $\TT\mathcal S(\omega)$ there exists 
a unique edge $e_p$ of $C$
such that $p\in h(e_p)$.
Let us equip $\omega$  with a phase structure $\phi_\omega$ in
$\Arg(\X)$, 
i.e. 
we equip each point $p$ in $\omega$ with a phase 
$\phi_p \in\Arg(\X)$.
Let us choose an approximation $(\omega_t^\CC)_{t>0}$ of $(\omega,\phi_\omega)$ by
generic configurations in $\X$. That is to say for each 
$t>0$, 
we choose a point $p_t\in\X$  for each $p\in\omega$ 
in such a way that the configuration 
$\omega_t^\CC$ formed by those points is generic, and that
$(p,\phi_p)$ is the tropical limit of $(p_t)$.

\begin{proposition}\label{lem:S to TS}
Let $(h,\phi)$ be an accumulation point of the sequence of sets
$\S(\Delta,\omega_{t}^\CC)$, in the sense of a tropical limit.
Then $h$ is an
element of $\TT\S(\omega)$.
Moreover for any  $p\in\omega$,  one has $\phi_p\in\Phi_{e_p}(S_{e_p})$. 
\end{proposition}
\begin{proof}
Let $(h:C\to \RR^3,\phi)$ be such an accumulation point. Since
$f_t(S_t)\subset\X$ for any element $f_t:S_t\to\X$ of
$\S(\Delta,\omega_{t}^\CC)$, it follows from Proposition \ref{prop:RH} that $h$
is a tropical morphism to $X$. Since in addition $\omega_{t}^\CC\subset
f_t(S_t)$, we clearly have $\omega\subset h(C)$ and 
$\phi_p\in\Phi_{e_p}(S_{e_p})$.

\vspace{1ex}
Next we show that the Newton fan of $h$ is equal to $\Delta$. Let us
first look at the projection to the $(x,y)$-coordinates. According to Lemma \ref{lem:hirzebruch}, each curve
$\pi_z(f_{\tj}(S))\subset\ctor$ is given by an equation
$P_{\tj}(x,y)=0$. Moreover, up to rescaling
the coefficients of $P_{t_j}$,
 we may suppose that the
biggest absolute value of the coefficients of $P_{\tj}$ is equal to
one.
We denote by $|\delta_z|$ the linear system on
$\mbox{Tor}(\Pi_{\delta_z})=\Sigma_n$ defined by the fan $\delta_z$. It is
naturally a projective space of dimension $N=\ZZ^2\cap \Pi_{\delta_z}-1$.
We denote by $\mathcal V_{\delta_z,2\chi}$ 
 the closure in $\mathbb C P^N$ of the set of all nodal complex
irreducible algebraic curves in $\ctor$ with Newton fan
$\delta_z$, and whose 
normalization has Euler characteristic $2\chi$. 

Passing through a point in $\ctor$ imposes a linear condition on 
curves in $|\delta_z|$. 
Hence curves in $|\delta_z|$ which pass through all points in
$\omega_{\tj}^\CC$ form a linear subspace $\L_{\tj}\subset
|\delta_z|$. 
By construction all points $\omega_{\tj}^\CC$ have a tropical limit in
$(\CC^*) ^3$,
so the coefficients of the equations defining
$\L_{\tj}$ may also be chosen to have a tropical limit in
$\CC^*$, i.e. they are
all equivalent to some function $ct_j^\lambda$. 
For each $\tj$, elements of $\S(\delta_z,\omega_{\tj}^\CC)$ correspond precisely
to  intersections of $\mathcal V_{\delta_z,2\chi}$ with the linear space
$\L_{tj}$. Since equations defining $\mathcal V_{\delta_z,2\chi}$ do not
depend on $t$, it follows from the analytic dependency of the root of
a polynomial with respect to its coefficients that all coefficients of $P_{\tj}$
also have a tropical limit in $\CC^*$.
Hence according to Lemma \ref{lem:NP}, the tropical curve
$\lim\Log_{tj}(\{P_{\tj}=0\})$ has Newton polygon $\Pi_{\delta_z}$.

\vspace{1ex}
It follows from what we just proved
that the Newton fan of $C$ has less elements than $\Delta$,
since $C$ might have ends of weight at least 2.
We also have $\chi(C)\le 2\chi$ according to Lemma \ref{lem:irr}.
Since $\omega$ is generic,
it follows from Proposition \ref{prop-curves} that  $\chi(C)= 2\chi$, and that
the Newton fan of $C$ has as many elements than $\Delta$, and so is
equal to $\Delta$.
\end{proof}

The projection $\pi_x$ relates the number
$N_{2\chi}^{irr}(\Delta)$ to another enumerative invariant of some
toric surface. 
This relation
 will allow us to relate the multiplicity of a
tropical curve in $\TT\S(\omega)$ 
to an actual number of complex curves in
$\S(\Delta,\omega^\CC)$. Together with Lemma \ref{lem:hirzebruch} this
will imply Theorem \ref{thm-corres}.
Let us denote by $\L$ the compactification in $\mbox{Tor}(\Pi_{\delta_x})$ of the
line in $\ctor$ defined by the equation $y+z=0$, and let us denote by
$\E$ the toric divisor of $\mbox{Tor}(\Pi_{\delta_x})$ corresponding to the
elements $(1,1)$ of $\delta_x$.
Let us denote by $\S(\delta_x,\omega^\CC)$ the set of all irreducible 
algebraic curves in $\ctor$ with Newton fan $\delta_x$, of Euler
characteristic $2\chi-\#\delta_x$, passing through all points in
$\pi_x(\omega^\CC)$, and whose closure in $\mbox{Tor}(\Pi_{\delta_x})$ has an ordinary
multiple point of multiplicity $an+b$ at the point $\L\cap\E$.

\begin{lemma}\label{lem:gen mult}
If $\omega^\CC$ is generic, then the set $\S(\delta_x,\omega^\CC)$ is
finite. Moreover if $S\in\S(\delta_x,\omega^\CC)$ and $D$ is a branch
of $S$ such that the closure of $D$
in $\mbox{Tor}(\Pi_{\delta_x})$ 
intersects $\L$, then this intersection is transverse.
\end{lemma}
\begin{proof}
We define $N=\#\left( \Pi_{\delta_x}\cap\ZZ^2 \right)-1$. 
The space $\mathcal V_{\delta_x,2\chi}$ of irreducible curves in $\ctor$ with
Newton fan $\delta_x$ 
and Euler characteristic $2\chi-\#\delta_x$ has dimension $\#\delta_x
-\chi$ and is naturally a quasiprojective variety in the linear system
$|\Pi_{\delta_x}|=\CC P^{N}$ on $\mbox{Tor}(\Pi_{\delta_x})$.
 Passing through a generic configuration of points $\omega^\CC$
imposes $\#\delta_z -\chi$ linearly independent conditions on curves in 
$|\Pi_{\delta_x}|$, and having a point of
multiplicity $an+b$ at $\L\cap\E$ imposes $an+b$ extra linearly independent
conditions. Hence the set $\S(\delta_x,\omega^\CC)$ is the
intersection of $\mathcal V_{\delta_x,2\chi}$ with a linear space of
complementary codimension. For such a generic linear space, the
intersection will be finite and transverse.
So for a generic configuration $\omega^\CC$, the set
$\S(\delta_x,\omega^\CC)$ is finite, and any of its elements  
cannot satisfy any
 further independent condition, like for example  having $\L$ as a tangent.
\end{proof}

\begin{lemma}\label{lem:multiple}
If $\omega^\CC$ is generic, then 
the projection $\pi_x$ establishes a bijection between the sets
$\S(\Delta,\omega^\CC)$ and $\S(\delta_x,\omega^\CC)$.
\end{lemma}
\begin{proof}
The surface $\X$ is the image of the map
$$\begin{array}{cccc}
\iota:& (\CC^*)^2\setminus\L &\longrightarrow & (\CC^*)^3
\\ &(y,z)&\longmapsto & (-y-z,y,z)
\end{array}. $$
In particular, any irreducible curve $S$ in $(\CC^*)^2$ intersecting
$\L$ in finitely many points
has a birational lift $\iota(S\setminus \L)$ in $\X$. 
Suppose now that $S$ has Newton fan $\delta_x$.
Then a puncture $p$ of $S$ corresponding to a vector $(1,1)$ of its
Newton fan will 
be a puncture of $\iota(S\setminus \L)\subset(\CC^*)^3$ corresponding to
a vector $(1-l,1,1)$ where $l$ is the order of contact of the closure
of $S$ and $\L$ at the point $\L\cap\E$ in $\mbox{Tor}(\Pi_{\delta_x})$. 
Hence if $\omega^\CC$ is generic, it follows from Lemma \ref{lem:gen mult} that
$\iota(S\setminus \L)\in \S(\Delta,\omega^\CC)$ if and only 
if $S\in \S(\delta_x,\omega^\CC)$.
\end{proof}

Hence it remains to compute, given an element
$h\in\TT\S(\omega)$, how many elements of
$\S(\Delta,\omega_{t}^\CC)$ have their amoeba $\Log_t(f_{t}(S_{t}))$
contained in a small neighborhood of $h(C)$.
To do so, we work with the projection $\pi_x$ in order to apply 
results from {\cite[Section 8.2]{Mi03}}.
In this latter paper, everything is stated in terms of curves in
$\ctor$ given by an equation. However, the generalization to maps from abstract curves to
$\ctor$ is straightforward.

Given a tropical morphism $h:C:\to\RR^2$, we denote by $w(e)$ the
weight of the edge $e$ of $C$ for $h$ (see Definition \ref{def:trop morph}).
\begin{proposition}\label{prop:numb cplx}
Let $h:C\to X$ be a tropical morphism in
$\TT\S(\omega)$
equipped with a phase $\phi$ such that 
$\phi_p\in \Phi_{e_p}(S_{e_p})$ for all points $p$ in $\omega$.
Then there exists a sequence $t_j\to +\infty$ such  that exactly 
$$\prod_{p\in \omega}w(e_p)$$
elements of
$\S(\Delta,\omega_{t_j}^\CC)$ converge tropically to $(h,\phi)$.
\end{proposition}
\begin{proof}
The tropical morphism $h$ composed with the projection $\pi_x$ induces
a tropical morphism $\overline h:\overline C\to\RR^2$, where
$\overline C$ is obtained from $C$ by contracting all ends $e$ with
$v(e)=(-1,0,0)$ (i.e. $C$ is an open tropical modification
of $\overline C$, see {\cite[Section 2.1]{Br9}}). 
Moreover, the
weight of an edge of $\overline C$ is the same as the weight of the
corresponding edge of $C$. Hence
according to Lemma \ref{lem:multiple}, it is equivalent to prove that
there exist a sequence $t_j\to +\infty$ such  that exactly 
$\prod_{p\in \pi_x(\omega)}w(e_p)$
elements of
$\S(\delta_x,\omega_{t_j}^\CC)$ converge tropically to
$(\overline h,\pi_x(\phi))$.

The proof of this latter statement
follows exactly along the lines of {\cite[Section 8.2]{Mi03}}. The only
point to check is that the tropical
morphism $\overline h:\overline C\to \RR^2$ is \emph{regular} in the sense of
{\cite[Definition 2.22]{Mi03}}. 
Let $\widetilde h:\widetilde C\to \RR^2$ be the tropical morphism to
$\RR^2$ such that $\widetilde h(\widetilde C)= \overline h(\overline
C)$ (set theoretically), where $\widetilde C$ is
obtained from $\overline C$ by identifying all ends of $\overline C$ mapped to
$\pi_x(L)$ and adjacent to a common
 vertex of $\overline C$. Clearly, the morphism 
$\overline h:\overline C\to \RR^2$ uniquely determines
$\widetilde h:\widetilde C\to \RR^2$, and their spaces of deformation
are canonically isomorphic.

It follows from Proposition \ref{prop-curves} that
the tropical curve $\widetilde C$ is 3-valent and $\widetilde h$ is an
immersion, so by {\cite[Proposition 2.23]{Mi03}}, the
dimension of the space of deformations of $\widetilde h$ is equal to
$$\Ed^{\infty}(\widetilde C) +g(\widetilde C) -1
= \Ed^{\infty}(\overline C) +g(\overline C) -1
- \sum_{V\in\Ve^0(\overline  C)}(\val(V)-3).$$
Thus $\overline h$ is regular, and all proofs in  {\cite[Section
8.2]{Mi03}} apply literally. The fact that we require  the
complex morphisms 
to pass through a (unique) point $an+b$ times instead
of requiring to pass through   $an+b$ distinct points in general
position does not make any difference 
since these conditions provides
independent equations. 
\end{proof}

Now we have to compute, given a tropical morphism $f:C\to \RR^2$, 
the number of phases $\phi$ we can endow $f$ with, such that 
$\phi_p\in \Phi_{e_p}(S_{e_p})$ for all points $p$ in $\omega$.
Again, this can be
done by a straightforward adaptation of {\cite[Section 8.2]{Mi03}}.
There is only one local computation needed here which is not covered
by
\cite{Mi03}, and that we perform now.

We first fix some notation.
Let  $T_0$ be
the triangle with vertices $(0,0),(l_1,0),$ and $(wl_2,wl_3)$ where $l_1,l_2,l_3$
and $w$ are four positive integers such that $\gcd(l_2,l_3)=1$. 
We denote by $\E_1$ (resp. $\E_2$ and $\E_3$)  the toric divisor of $\mbox{Tor}(T_0)$
corresponding to $[(0,0);(l_1,0)]$ (resp. $[(0,0);(wl_2,wl_3)]$ and 
$[(l_1,0);(wl_2,wl_3)]$).  Finally we choose a point $p_1$
 on $\E_1\setminus\left(\E_2\cup \E_3 \right)$,  and a point $p_2$
 on $\E_2\setminus\left(\E_1\cup \E_3 \right)$.

Note that the restriction to the case of $T_0$ in the following lemma
does not cause any loss of generality. Indeed let $T$ be a triangle with
vertices in $\ZZ^2$, and choose any two points $q_1$ and $q_2$ on two
different toric divisors of $\mbox{Tor}(T)$ such that neither $q_1$ or $q_2$
is the intersection point of two toric divisors.
Then there exists a unique choice of
$l_1$, $l_2$, $l_3$ and $w$ and  a unique toric
isomorphism $\mbox{Tor}(T)\to \mbox{Tor}(T_0)$ mapping $q_1$ to $p_1$ and $q_2$ to
$p_2$, and sending the linear system defined by $T$ on $\mbox{Tor}(T)$ to the linear
system defined by $T_0$ on $\mbox{Tor}(T_0)$.
\begin{lemma}\label{lem:local comp}
With the above notation, up to re-parameterization of $\CC P^1$, there
exist exactly $\frac{1}{w}\binom{wl_3}{l_1}$ algebraic maps
$f:\CC P^1\to \mbox{Tor}(T)$ such that 
\begin{itemize}
\item $f^{-1}(\E_1)=f^{-1}(p_1)$ and consists of $l_1$ distinct points
on $\CC P^1$;
\item $f^{-1}(\E_2)=f^{-1}(p_2)$ and consists of a single point;
\item $f^{-1}(\E_3)$ consists of a single point.
\end{itemize}
\end{lemma}
\begin{proof}
Since $\gcd(l_2,l_3)=1$ there exist two integers $u$ and $v$ such that $vl_2 -
ul_3=1$. We fix on $(\CC^*)^2$ the
unique choice of coordinates $(x,y)$  such that
$p_1=\psi(p_2)=(1,0)$ where $\psi$ is the automorphism of $(\CC^*)^2$
given by $\psi(x,y)= (x^{l_2}y^{l_3}, x^uy^v)$.
We also choose a coordinate system on $\CC P^1$ such that 
$f^{-1}(p_2)=\{0\}$,  $f^{-1}(\E_3)=\{\infty\}$, and $f(1)=p_1$.

Hence all curves we are counting may be parameterized by a 
map of the form
$$\begin{array}{cccc}
f: & \CC^*\setminus\{1\}\setminus\{z\  | \ Q_2(z)=0\} &\longrightarrow& \ctor
\\& z &\longmapsto & (z^{wl_3}, \frac{(z-1)Q_2(z)}{z^{wl_2}})
\end{array} $$
where $Q_2(z)$ is a  complex polynomial of degree $l_1-1$. 
From the conditions imposed on $f$, 
the roots of $Q_2$ must be a subset  of $l_1-1$ elements of the set of
all $(wl_3)$-th roots of unity distinct from $1$. Hence we have 
$\binom{wl_3-1}{l_1-1}$ distinct possibilities for the roots of
$Q_2$, and any such choice  determines $Q_2$ up to a multiplicative
constant.
The equation that this constant  satisfies is given by the
condition $f(0)=p_2$. 
Since we have 
$\psi\circ f (z)=\left((z-1)^{l_3}Q_2(z)^{l_3},
z^{w}(z-1)^vQ_2(z)^v  \right)$ this condition translates to 
$(-1)^{l_3}Q(0)^{l_3}=1$. Hence there exists exactly $l_3$ distinct
polynomials $Q_2(z)$ once the set of its roots is chosen.

In conclusion there exist exactly $l_3\binom{wl_3-1}{l_1-1}$
admissible functions $f$ with the chosen coordinates system on $\CC
P^1$. Since there exist $l_1$ possibilities to choose
a point $q$ in $f^{-1}(p_1)$ in order to fix the coordinate system on
$\CC P^1$, the number of  maps $f$ up to re-parameterization of $\CC
P^1$ is $\frac{l_3}{l_1}\binom{wl_3-1}{l_1-1}=\frac{1}{w}\binom{wl_3}{l_1} $.
\end{proof}

In the next proposition we apply this
computation to vertices of the tropical curve mapped to $L$. After the projection with $\pi_x$, these vertices are adjacent to $l$ edges of the same direction, corresponding to the $l_1$ preimages of $p_1$ here. From this computation, we obtain the binomial coefficients $\binom{k+l}{l}$ that appear in the definition of tropical multiplicity of a vertex mapped to $L$.

\begin{proposition}\label{prop:numb phase}
Let $h:C\to X$ be a tropical morphism in
$\TT\S(\omega)$. Then there exist exactly 
$$\frac{\mu_h}{\prod_{p\in e}w(e)} $$
possibilities to choose a phase $\phi$ for $h$ such that
$\phi_p\in \Phi_{e_p}(S_{e_p})$ for all points $p$ in $\omega$.
\end{proposition}
\begin{proof}
As in the proof of Proposition \ref{prop:numb cplx}, we can equivalently compute 
how many possibilities there are to phase $\overline h:\overline
C\to\RR^2$ 
in a coherent way with
$\pi_x(\phi_\omega)$.
Again we follow the lines of \cite{Mi03}. As in {\cite[Lemma
4.20]{Mi03}}, any connected component of $C\setminus\P$ is a tree
containing exactly one end not mapped to $L$, otherwise the set
$\TT\S(\omega)$ would not be finite. Hence we may reconstruct
all possible phases step by step as in  {\cite[Section
8.2]{Mi03}}, each step consisting of solving a simple enumerative
problem in a toric surface. 
The only local computation we need here
 which is not covered by \cite{Mi03} is done in Lemma \ref{lem:local comp}. 
\end{proof}

Now Theorem \ref{thm-corres} follows immediately from
Propositions   \ref{lem:S
to TS}, \ref{prop:numb cplx}, and \ref{prop:numb phase}.

\subsection{The proof of Theorem \ref{thm:corres 2}}\label{sec:inter trop numbers corres}

The proof of Theorem \ref{thm:corres 2} follows exactly the same strategy
as the proof of Theorem \ref{thm-corres}: we  reformulate the
definition of the
numbers $\mathcal N(u,n,d,\alpha)$ in terms of enumeration of curves
in some toric surface having prescribed intersections with the toric
divisors; then we prove a Correspondence Theorem for this latter
enumerative problem by
applying the methods from \cite{Mi03} and \cite{Shu12}.
As in the proof of Theorem \ref{thm-corres}, we first ``lift'' our enumerative problem to three-space and then project it back to the plane using the different projection $\pi_x$. In this way, we can make use of known techniques.

We fix a coordinate system  $(x,y)$ on $\ctor$ such that
$(x,\frac{1}{y})$ is a standard coordinate system on $\Sigma_n$.
Let $S_0$ be the curve of bidegree $(1,1)$ in $\Sigma_n$ with equation
$x+y=0$, and let $F_0$ be the curve with equation $x=\infty$.
We define
$$\begin{array}{cccc}
\iota': & \ctor\setminus \{x+y=0\} &\longrightarrow & (\CC^*)^3
\\& (x,y) &\longmapsto & (x,y,-x-y)
\end{array} .$$
The next lemma is a straightforward computation.
\begin{lemma}\label{lem:newtonfan}
Let $C$ be an algebraic curve in
$\Sigma_n$ of bidegree $(a,0)$, and let
$D$ be a
  branch of $C$ such that 
$(D\circ S_0)_{p_0}=d+\alpha(n+1)$ and $(D\circ
  F_0)_{p_0}=\alpha$. 
If $C'$ denotes the restriction of $C$ to the torus orbit
corresponding to the coordinate system $(x,y)$ fixed above, then
the element of the Newton fan of $\iota'(C'\setminus \{x+y=0\})$
corresponding to $D$ is $(\alpha, \alpha, \alpha-d)$.
\end{lemma}

Let $\delta$ be the Newton fan
\begin{align*}\delta= \{ (1,1)^{un-\sum \alpha_i}, (-n, 1)^u,
  (\alpha_1, &\alpha_1-d_1),\ldots,  (\alpha_s,
  \alpha_s-d_s),\\& (0,-d_{r+1}),\ldots,  (0,-d_{r+s})\}.\end{align*}
We denote by $\mathcal L$ the closure in $\mbox{Tor}(\Pi_\delta)$  of the
line in $\ctor$ with equation $y+z=0$, and by $\E$ the toric divisor of
$\mbox{Tor}(\Pi_\delta)$ corresponding to  elements $(1,1)$ of $\delta$.
Given an algebraic curve $C$ in $\ctor$, we denote by 
$\overline C$ its closure in $\mbox{Tor}(\Pi_\delta)$.

Using the projection $\delta_x$ and the  techniques from
Lemma \ref{lem:multiple} and Lemma \ref{lem:newtonfan}, 
we obtain the following corollary. 
\begin{corollary}\label{cor:other def of N}
The number $\mathcal N(u,n,d,\alpha)$ 
is equal to the number of algebraic curves $C$ in $\ctor$ with Newton
  fan $\delta$ and such that
\begin{itemize}
\item $\overline C$ has $u$ connected components, whose normalization are 
all rational;
\item for each element $(\alpha_i, \alpha_i-d_i)$ or $(0,-d_i)$ of
$\delta$, the intersection of $\overline C$ with the corresponding
toric divisor is fixed;
\item $\overline C$ has an ordinary multiple
point of multiplicity $un-\sum \alpha_i$ at the point $\mathcal L\cap \E$.
\end{itemize}
\end{corollary}

Now the proof of Theorem \ref{thm-corres} extends literally to Theorem 
\ref{thm:corres 2}, using refinements of \cite{Mi03} from
  \cite{Shu12} to enumerate
curves with prescribed intersections with toric divisors.

\section{Concluding remarks}\label{sec:rem}
We discuss some of the
 possible extensions of  results and methods
presented in this 
paper.

\begin{enumerate}
\item Although Theorem \ref{thm-corres}  assumes that configurations
$\omega$ are contained in the two faces $\sigma_1$ and $\sigma_2$, 
Theorem \ref{main formula} is obtained just by considering 
configurations
$\omega$ contained in $\sigma_1$. 
It should be possible to  generalize Theorem \ref{main formula} for any
configuration $\omega\subset \sigma_1\cup\sigma_2$. This would also
require  to enlarge the family of $(1,1)$-relative invariants 
considered here. 

\item It would also be  interesting to relate enumerative
invariants of $\Sigma_n$ and $\Sigma_{n+2k}$ when $k\ge 2$. According
to Appendix \ref{hirz}, one possible way would be to study enumerative
geometry of the tropical surface $X_k$ in $\RR^3$ given by the polynomial
$\tg x^{k} + y + z\td$. In this case the assumption we made throughout
Section \ref{sec:inter trop numbers}, i.e. that $d=\#\Delta_i$ for
some $i$, fails.  In particular the study of enumerative geometry of $X_k$
requires more care for $k\ge 2$.

\item Related to the previous remark is the question of determining the
multiplicity of a tropical morphism to $X$. In general, the multiplicity
of a vertex 
tropically
mapped to  the line $L$ should be expressed in terms of
\emph{triple Hurwitz numbers} weighted by some binomials coefficients.
Although all those numbers are in principle computable, 
 no nice general formula is known yet. 
In the particular case treated in this paper, the corresponding
Hurwitz numbers are very simple: it is the number of rational maps 
$\CC P^1\to \CC P^1$ of
degree $d$ with a prescribed pole and zero of maximal order. In
particular we could perform easily all computations keeping hidden the
Hurwitz numbers aspect. However for more general enumerative problems
in $X$, these Hurwitz numbers will show up naturally.

\item  More generally, the study of enumerative geometry of general
tropical surfaces, or even tropical varieties of any dimension, is of
great interest. So far, little is known about this problem. 
In this case  all Hurwitz numbers will come into the game, not only the
triple ones mentioned above.
A generalization of
Proposition \ref{prop:numb cplx} in nice cases should be a consequence of 
general results about the realization  of regular phase-tropical curves in
the
forthcoming paper \cite{Mik08}.

\end{enumerate}

%% file: Ex_mor3.pstex_t
\begin{picture}(0,0)%
\includegraphics{Ex_mor3.pstex}%
\end{picture}%
\setlength{\unitlength}{2368sp}%
\begingroup\makeatletter\ifx\SetFigFont\undefined%
\gdef\SetFigFont#1#2#3#4#5{%
  \reset@font\fontsize{#1}{#2pt}%
  \fontfamily{#3}\fontseries{#4}\fontshape{#5}%
  \selectfont}%
\fi\endgroup%
\begin{picture}(3999,1524)(6364,-2923)
\put(8176,-1861){\makebox(0,0)[lb]{\smash{{\SetFigFont{6}{7.2}{\familydefault}{\mddefault}{\updefault}{\color[rgb]{0,0,0}$h$}%
}}}}
\put(7126,-1786){\makebox(0,0)[lb]{\smash{{\SetFigFont{6}{7.2}{\familydefault}{\mddefault}{\updefault}{\color[rgb]{0,0,0}$V_1$}%
}}}}
\put(6976,-2611){\makebox(0,0)[lb]{\smash{{\SetFigFont{6}{7.2}{\familydefault}{\mddefault}{\updefault}{\color[rgb]{0,0,0}$V_0$}%
}}}}
\put(9526,-2461){\makebox(0,0)[lb]{\smash{{\SetFigFont{6}{7.2}{\familydefault}{\mddefault}{\updefault}{\color[rgb]{0,0,0}$(0,0)$}%
}}}}
\put(9376,-1861){\makebox(0,0)[lb]{\smash{{\SetFigFont{6}{7.2}{\familydefault}{\mddefault}{\updefault}{\color[rgb]{0,0,0}$(1,1)$}%
}}}}
\end{picture}%